\documentclass[a4paper]{article}
\usepackage{lipsum}
\usepackage{amsfonts}
\usepackage{graphicx}
\usepackage{epstopdf}
\usepackage{algorithm}
\usepackage{algorithmic}
\ifpdf
  \DeclareGraphicsExtensions{.eps,.pdf,.png,.jpg}
\else
  \DeclareGraphicsExtensions{.eps}
\fi

\usepackage{amsmath, amsfonts, amssymb}
\usepackage{color, xcolor}
\usepackage{siunitx}
\usepackage{subfig}
\usepackage{multirow}
\usepackage{hyperref}
\usepackage{cleveref}
\usepackage{booktabs}
\usepackage{geometry}
\usepackage[todonotes={textsize=tiny}]{changes}
\definechangesauthor[name={Zhicheng HU}, color=red]{zc}
\setdeletedmarkup{{\color{blue} \sout{#1}}}
\usepackage{tikz}
\usetikzlibrary{calc} 
\usetikzlibrary {decorations.fractals,spy}
\graphicspath{{images/}}
\crefname{equation}{}{}

\newcommand\dd{\,\mathrm{d}}

\newcommand\bbR{\mathbb{R}}
\newcommand\bbN{\mathbb{N}}

\newcommand\He{\mathit{He}}
\newcommand\Kn{\mathit{Kn}}

\newcommand\mM{\mathcal{M}}
\newcommand\mF{\mathcal{F}}

\newcommand\mH{\mathcal{H}}

\newcommand\mR{\mathcal{R}}
\newcommand\mO{\mathcal{O}}

\newcommand\bA{{\bf{A}}}

\newcommand\bS{\boldsymbol{S}}
\newcommand{\bbf}{\boldsymbol{f}}

\newcommand\bR{\boldsymbol{R}}

\newcommand\bxi{\boldsymbol{\xi}}
\newcommand\bx{\boldsymbol{x}}

\newcommand\bg{\boldsymbol{g}}
\newcommand\bq{\boldsymbol{q}}
\newcommand\bu{\boldsymbol{u}}
\newcommand\bv{\boldsymbol{v}}

\newcommand\bn{\boldsymbol{n}}
\newcommand\br{\boldsymbol{r}}

\newcommand\bLambda{{\boldsymbol{\Lambda}}}

\newcommand\bvpi{\boldsymbol{\varpi}}
\newcommand\bvrho{\boldsymbol{\varrho}}

\newcommand\pd[2]{\dfrac{\partial {#1}}{\partial {#2}}}

\newcommand\opd[2]{\dfrac{\dd {#1}}{\dd {#2}}}

\title{An efficient nonlinear multigrid solver for the simulation of 
  rarefied gas cavity flow\thanks{This work was partially supported by the National Natural
      Science Foundation of China, No. 12171240, and the Fundamental
      Research Funds for the Central Universities, China,
      No. NS2021054. The computational resources were supported by
      High Performance Computing Platform of Nanjing University of
      Aeronautics and Astronautics, China.}}

\author{Zhicheng Hu\thanks{
    Department of Mathematics, College of Science,
    Nanjing University of Aeronautics and Astronautics,
    Nanjing 211106, China; Key Laboratory of Mathematical Modelling
    and High Performance Computing of Air Vehicles (NUAA), MIIT,
    Nanjing 211106, China. Email: huzhicheng@nuaa.edu.cn (Z. Hu),
    liguanghan@nuaa.edu.cn (G. Li).} \and Guanghan
  Li\footnotemark[2]}

\date{}
\begin{document}
\maketitle

\begin{abstract}
  We study efficient simulation of steady state for rarefied gas flow,
  which is modeled by the Boltzmann equation with BGK-type collision
  term. A nonlinear multigrid solver is proposed to resolve the
  efficiency issue by the following approaches. The unified framework
  of numerical regularized moment method is first adopted to derive
  the high-quality discretization of the underlying problem. A fast
  sweeping iteration is introduced to solve the derived discrete
  problem more efficiently than the usual time-integration scheme on a
  single level grid. Taking it as the smoother, the nonlinear
  multigrid solver is then established to significantly improve the
  convergence rate. The OpenMP-based parallelization is applied in the
  implementation to further accelerate the computation. Numerical
  experiments for two lid-driven cavity flows and a bottom-heated
  cavity flow are carried out to investigate the performance of the
  resulting nonlinear multigrid solver. All results show the
  efficiency and robustness of the solver for both first- and
  second-order spatial discretization.
\end{abstract}


\section{Introduction}
\label{sec:intro}

Rarefied gas dynamics is a classical discipline of fluid mechanics,
which studies gas flows in the situation when the molecular mean free
path, compared with the characteristic length of problems, can not be
negligible. Over the past several decades, it has attracted more and
more attention in the high-tech fields including astronautics and
micro-electro-mechanical systems, due to the rapid development of
aerospace and semiconductor industries, etc. To describe the rarefied
gas flows correctly, the Boltzmann equation, rather than the
traditional continuum models such as Euler equations and Navier-Stokes
equations, should be taken into account. Clearly, numerical methods
are necessary to solve the Boltzmann equation. In practice, for a
variety of applications, the original quadratic collision term of the
Boltzmann equation can be approximated by a simplified relaxation
model, e.g., the Bhatnagar-Gross-Krook (BGK) model \cite{BGK}, to get
a considerable reduction of computational effort without loss of
accuracy. However, because of the intrinsic high dimensional nature,
it still encounters a big challenge to develop an accurate and
efficient solver for the simplified Boltzmann equation, especially for
a great of important purposes \cite{nabapure2021dsmc,Zhu2021,Yang2021}
when the steady state of problems has to be under investigation.

This work is concerned with efficient simulation of steady state for
rarefied gas flows modeled by the Boltzmann equation with BGK-type
relaxation models. To this end, we are interested in the framework of
numerical regularized moment method originated in \cite{NRxx} and then
developed in \cite{Cai, Li, Microflows1D, framework} for
discretization of the underlying problem. As one of the most powerful
methods, the moment method derives a system of continuous equations,
which can be viewed not only as an extended hydrodynamic model but
also a high-order velocity discretization of the Boltzmann
equation. The numerical regularized moment method enhances this system
with a slight modification so that some desired properties, such as
global hyperbolicity \cite{Fan_new} and convergence \cite{Cai2020,
  cai2021moment}, could be well guaranteed for the revised moment
system up to arbitrary order. Moreover, in the framework of numerical
regularized moment method, the full discretization would be obtained
formally by using the finite volume method for spatial discretization
without explicitly writing out the moment system. Such a unified
framework makes the practical application of high-order moment system
much easier. Additionally, both first- and second-order spatial
discretization are of interest in the present paper.

To efficiently solve the derived discrete steady-state problem, a fast
sweeping iteration, based on the forward Euler scheme and the
cell-by-cell Gauss-Seidel iteration with four alternating direction
sweepings in two-dimensional (2D) case, is first proposed on a single
level grid. In comparison to the forward Euler scheme, the present
fast sweeping iteration would converge much faster, while takes for
each iteration almost the same computational cost as four steps of the
forward Euler scheme. It turns out that the fast sweeping iteration is
usually more efficient than the forward Euler scheme. Indeed, it can
also be observed that the fast sweeping iteration is more robust than
the forward Euler scheme especially when the second-order spatial
discretization is under consideration.

In view of that the multigrid method
\cite{brandt2011book,hackbusch1985book} is one of the most popular
acceleration techniques for steady-state computation, we then
concentrate on the multigrid acceleration for our discrete
problem. Actually, a nonlinear multigrid (NMG) iteration, which takes
an SGS-Newton iteration as the smoother, has been successfully
developed for the hyperbolic moment models in one-dimensional case
\cite{hu2014nmg}. In the current work, we generalize this NMG
iteration to 2D case by using the previous fast sweeping iteration
instead of the SGS-Newton iteration, and providing an appropriate
construction of 2D restriction and prolongation operators. The
resulting NMG solver would always converge within dozens of
iterations, leading to a significant improvement in efficiency in
comparison to the fast sweeping iteration. To further accelerate the
computation, the OpenMP-based parallelization \cite{dagum1998openmp}
is employed in the implementation of the NMG solver. Almost all
operations, including \emph{for} loop to traverse grid cells of a
given grid, e.g., smoothing, restriction and prolongation, can be
parallelized trivially without appreciably affecting the results,
although it should be pointed out that the parallel fast sweeping
iteration is not exactly equivalent to the serial fast sweeping
iteration.

Finally, a number of numerical experiments for 2D cavity flows such as
two lid-driven flows and a bottom-heated flow are presented to explore
the behavior of the resulting NMG solver. The efficiency and
robustness of it are well validated by numerical results for both
first- and second-order spatial discretization.

The rest of this paper is organized as follows. In \cref{sec:model}, a
brief review of the Boltzmann equation and the full discretization as
well as the basic time-integration scheme are introduced. The details
of three acceleration methods, namely, fast sweeping iteration,
nonlinear multigrid method, and OpenMP-based parallelization, are
described in \cref{sec:solver} to establish the final steady-state
solver. Numerical experiments are presented in \cref{sec:num-ex}, and
a brief summary is given in the last section.

\section{Model equations and basic numerical method}
\label{sec:model}

In this section, we will give a brief review of the governing model
equations for rarefied gas dynamics, and introduce a unified framework
of numerical method on a given single level grid by using the
numerical regularized moment method.

\subsection{Boltzmann equation with BGK-type collision term}
\label{sec:model-bolt}

Let $f(t,\bx,\bxi)$ be the distribution function such that
$f(t,\bx,\bxi)\dd \bx \dd \bxi$ gives the number of particles that
locate in the infinitesimal cell $\dd \bx$ about $\bx$ with velocities
within the infinitesimal cell $\dd \bxi$ about $\bxi$ at time $t$. For
the rarefied gas in a 2D cavity, the evolution of $f$ can be well
described by the Boltzmann equation
\begin{align}
    \label{eq:boltzmann}
  \pd{f}{t} + \bxi \cdot \nabla_{\bx} f  = Q[f,f], \quad t \in \bbR^{+}, \quad \bx \in \Omega \subset \bbR^{2}, \quad \bxi \in \bbR^{3},
\end{align}
subject to some appropriate boundary conditions and the assumption
$\frac{\partial f}{\partial x_{3}} \equiv 0$. Here the right-hand side
$Q[f,f]$ is the collision term representing the interaction between
gas molecules. It usually takes the quadratic form
\begin{align}
  \label{eq:binary-collision}
  Q[f,f](t,\bx,\bxi) = \int_{\bbR^{3}} \int_{\mathbb{S}_{+}^{2}}
  B(\vert \bxi-\bxi_*\vert, \bn ) (f'f_*'-ff_*) \dd \bn \dd \bxi_*,
\end{align}
where the collision kernel $B(\cdot,\cdot)$ is a non-negative function
determined by the potential between gas molecules,
$\bn \in \mathbb{S}_{+}^{2}$ is a unit vector directed along the line
joining the centers of two colliding particles with pre-collision
velocities $\bxi$ and $\bxi_{*}$, and $f_{*}=f(t,\bx,\bxi_*)$,
$f' = f(t,\bx, \bxi')$, $f_{*}'= f(t, \bx, \bxi_*')$, in which $\bxi'$
and $\bxi_*'$, given by
\begin{align*}
  \bxi' = \bxi - [(\bxi-\bxi_{*})\cdot \bn]\bn, \quad  \bxi_{*}' = \bxi_{*} - [(\bxi_{*}-\bxi)\cdot \bn]\bn,  
\end{align*}
are the corresponding post-collision velocities of the two colliding
particles.

While the distribution function $f$ gives a detailed information of
the state of the gas, in most practical applications, we are mainly
concerned about the macroscopic physical quantities, such as the mass
density $\rho$, mean velocity $\bu$, and temperature $\theta$, rather
than the distribution function itself. For these quantities of
interest, they are related to the moments of the distribution function
as follows
\begin{align}
  \label{eq:rho-u-theta}
\begin{aligned}
  & \rho(t,\bx) = m_* \int_{\bbR^3} f(t,\bx,\bxi) \dd \bxi, \\ &
  \rho(t,\bx) \bu(t,\bx) = m_* \int_{\bbR^3} \bxi f(t,\bx,\bxi) \dd \bxi, \\
  & \rho(t,\bx) \vert \bu(t,\bx) \vert^2 + 3 \rho(t,\bx) \theta(t,\bx) = m_*
  \int_{\bbR^3} \vert \bxi \vert^2 f(t,\bx, \bxi) \dd \bxi,
\end{aligned}
\end{align}
where $m_*$ is the mass of a single gas molecule. Similarly, the
stress tensor $\sigma_{ij}$, $i,j=1,2,3$, and the heat flux $\bq$ are
defined by
\begin{align}
  \label{eq:sigma-heat}
\begin{aligned}
  & \sigma_{ij}(t,\bx) = m_* \int_{\bbR^3} (\xi_i - u_i(t,\bx))(\xi_j -
  u_j(t,\bx)) f(t,\bx,\bxi) \dd \bxi - \rho(t,\bx) \theta(t,\bx) \delta_{ij},
  \\ & \bq(t,\bx) = \frac{m_*}{2} \int_{\bbR^3} \vert
  \bxi-\bu(t,\bx) \vert^2 (\bxi-\bu(t,\bx)) f(t,\bx,\bxi) \dd \bxi,
\end{aligned}
\end{align}
where $\delta_{ij}$ is the Kronecker delta symbol.

It is obvious that the quadratic collision term
\cref{eq:binary-collision} is the most complicated part in the
Boltzmann equation, which leads to a great challenge for efficient
numerical simulation. In recent years, some fast numerical methods
have been successfully proposed to partially resolve the efficiency
issue in the computation of \cref{eq:binary-collision}, see e.g.,
\cite{Wu2013,Gamba2017} for Fourier spectral method, \cite{Wang2019}
for Hermite spectral method and \cite{Cai2019} for Burnett spectral
method. However, lots of researchers are still very interested in the
approach to approximate the quadratic collision term
\cref{eq:binary-collision} by a simpler collision model, which is able
to predict the major physical features of interest with a much smaller
computational cost, in a number of situations at moderate Knudsen
numbers. The most famous class of simplified collision models are the
BGK-type relaxation models, which have the general form of
\begin{align}
  \label{eq:col-relaxation}
  Q[f,f]\approx Q(f) = \nu (f^{\text{E}} - f), 
\end{align}
where $\nu$ is the average collision frequency assumed independent of
the molecular velocity, and $f^{\text{E}}$ is the local equilibrium
distribution function depending on the specific model selected. In
particular, we have
\begin{itemize}
\item For the simplest BGK model \cite{BGK}, $f^{\text{E}}$ is the
  local Maxwellian given by
  \begin{align}
    \label{eq:BGK-f^E}
    f^{\text{E}}(t,\bx,\bxi)=\mM_{f}(t,\bx,\bxi) = \frac{\rho(t, \bx)}{m_* [2\pi
      \theta(t,\bx)]^{3/2}} \exp\left(-\frac{\vert \bxi -
        \bu(t,\bx) \vert^2}{2\theta(t,\bx)} \right).
  \end{align}
\item For the ellipsoidal statistical BGK (ES-BGK) model 
  \cite{Holway}, $f^{\text{E}}$ is an anisotropic Gaussian 
  distribution defined by
  \begin{equation}
    \label{eq:ES-f^E}
    f^{\text{E}}(t,\bx,\bxi) = \frac{\rho(t,\bx)}{m_* \sqrt{ \det[2\pi \bLambda]}
    } \exp\left(-\frac{1}{2} ( \bxi - \bu(t,\bx) )^T \bLambda^{-1}
      (\bxi-\bu(t,\bx)) \right),
  \end{equation}
  where $\bLambda=(\lambda_{ij})$ is a $3\times 3$ matrix given by
  \begin{align*}
    \lambda_{ij}(t,\bx) = \theta(t,\bx) \delta_{ij} + \left( 1-
      \frac{1}{\Pr} \right) \frac{\sigma_{ij}(t,\bx)}{\rho(t,\bx)}, \quad i,
    j = 1,2,3,
  \end{align*}
  in which $\Pr$ is the Prandtl number.
\item For the Shakhov model \cite{Shakhov}, $f^{\text{E}}$ is the
  product of the local Maxwellian and a cubic polynomial of $\bxi$ as
  \begin{equation}
    \label{eq:Shakhov-f^E}
    f^{\text{E}}(t,\bx, \bxi) = \left[ 1+\frac{(1-\Pr)(\bxi-\bu(t,
        \bx))\cdot \bq(t, \bx) }{5\rho(t,\bx) [\theta(t,\bx)]^2} \left(
        \frac{ \vert \bxi-\bu(t,\bx) \vert^2}{\theta(t,\bx)} - 5
      \right)\right] \mM_{f}(t,
    \bx, \bxi).
  \end{equation}
\end{itemize}
It is apparent from \cref{eq:rho-u-theta,eq:sigma-heat} that
$f^{\text{E}}$ is nonlinearly dependent on $f$ via the macroscopic
physical quantities. In addition, both the ES-BGK model and the
Shakhov model reduce to the BGK model in the case when $\Pr=1$.

In view of the enormous computational cost of the quadratic collision
term \cref{eq:binary-collision}, we shall mainly use the BGK-type
relaxation models as examples in this work, to illustrate the
framework of the present numerical method and to reveal the
main features of the resulting solver.

\subsection{Velocity discretization and moment system}
\label{sec:model-vdis}

In order to discretize the Boltzmann equation \cref{eq:boltzmann} in
velocity space, we would like to use the Hermite spectral method,
which first expands the distribution function into a series of Hermite
functions as
\begin{align}
  \label{eq:dis-expansion}
  f(t, \bx, \bxi) = \sum_{\alpha \in \bbN^{3}} f_\alpha(t, \bx)
  \mH_{\alpha}^{[\bvpi, \vartheta]} (\bxi), 
\end{align}
where $\mH_{\alpha}^{[\bvpi, \vartheta]}(\cdot)$ is the $\alpha$th
basis function defined by
\begin{align*}
  \mH_{\alpha}^{[\bvpi, \vartheta]}(\bxi) =
  \frac{1}{m_*(2\pi\vartheta)^{3/2}\vartheta^{|\alpha|/2}} \prod\limits_{d=1}^3
  \He_{\alpha_d}({v}_d)\exp \left(-\frac{{v}_d^2}{2} \right),
  \quad {\bv} = \frac{\bxi-\bvpi}{
    \sqrt{\vartheta}},~~\forall \bxi\in\bbR^3,
\end{align*}
and $f_{\alpha}(t,\bx)$ is the corresponding coefficient that is
independent of $\bxi$. In the above expression, $|\alpha|$ is the sum
of all components of the multi-dimensional index $\alpha$, i.e.,
$|\alpha| = \alpha_1 + \alpha_2 + \alpha_3$, and $\He_n(\cdot)$ is the
Hermite polynomial of degree $n$ given by
\begin{equation*}
  \He_n(x) = (-1)^n\exp \left(\frac{x^{2}}{2} \right) \frac{\dd^n}{\dd x^n} 
  \exp \left(-\frac{x^{2}}{2}\right).
\end{equation*}
The parameters $\bvpi\in \bbR^{3}$ and $\vartheta \in \bbR^{+}$ in the
basis function could be either constants or variables depending on $t$
and $\bx$. Besides to set them as constants for the most common
Hermite spectral method, a popular approach, originated by Grad
\cite{Grad}, is to adaptively choose $\bvpi$ and $\vartheta$ as the
local mean velocity $\bu$ and temperature $\theta$, respectively,
according to the distribution function $f$ itself via
\cref{eq:rho-u-theta}.
For the Grad method, we can easily deduce from
\cref{eq:rho-u-theta,eq:sigma-heat} that
\begin{align}\label{eq:moments-relation}
  \begin{aligned}
    &f_0 = \rho, \qquad f_{e_1} = f_{e_2} = f_{e_3} = 0, \qquad
    \sum_{d=1}^3 f_{2e_d} = 0,\\
    &\sigma_{ij} = (1+\delta_{ij}) f_{e_i+e_j}, \quad q_i = 2
    f_{3e_i} + \sum_{d=1}^3 f_{2e_d+e_i}, \qquad i,j=1,2,3,
  \end{aligned}
\end{align}
where $e_1$, $e_2$, and $e_3$ denote the multi-dimensional indices
$(1,0,0)$, $(0,1,0)$, and $(0,0,1)$, respectively. Moreover, it can be
expected that the series \cref{eq:dis-expansion} converges fast when
the distribution function is smooth and not too far away from the
Maxwellian. Actually, the first term of the series
\cref{eq:dis-expansion} already gives the local Maxwellian associated
with the distribution function $f$. For the general choice of $\bvpi$
and $\vartheta$, the macroscopic physical quantities of interest can
also be related directly to these two parameters and the low-order
coefficients, i.e., $f_{\alpha}$ with $|\alpha|\leq 3$. We refer to
\cite{Hu2019} for the detailed relations in this general
case. Additionally, as pointed out in \cite{hu2020burnett}, the
convergence of the series \cref{eq:dis-expansion} could usually be
expected for a smooth distribution function if $\theta < 2 \vartheta$.

Due to the importance of the Grad method, we will restrict ourselves
to the case of $\bvpi=\bu$ and $\vartheta=\theta$ for the
approximation of every distribution function in this
paper. Numerically, the distribution function is approximated by a
truncated series of order $M$, that is,
\begin{align}
  \label{eq:truncated-dis}
  f(t, \bx, \bxi) \approx \sum_{|\alpha| \leq M} f_\alpha(t, \bx)
  \mH_{\alpha}^{[\bu(t, \bx), \theta(t,\bx)]} (\bxi) \in \mF_{M}^{[\bu,\theta]},
\end{align}
where $\mF_{M}^{[\bu,\theta]}$ is the finite-dimensional linear space
spanned by $\mH_{\alpha}^{[\bu,\theta]}(\bxi)$ for all $\alpha$ with
$|\alpha|\leq M$. In view of the relation \cref{eq:moments-relation},
we have that the unknown variables in the above approximation are
$\bu$, $\theta$, and all coefficients $f_{\alpha}$ with
$|\alpha|\leq M$ except the coefficients with $|\alpha|=1$ and one of
the coefficients with $|\alpha|=2$. Let $\bvrho$ represent the vector
composed of these unknown quantities. Following the framework of
deriving the regularized moment model from the Boltzmann equation
\cref{eq:boltzmann} that was proposed in \cite{Fan_new,framework}, we
can get a globally hyperbolic system for $\bvrho$, which can be
written in a quasi-linear form as
\begin{align}
  \label{eq:moment-system}
  \pd{\bvrho}{t} + \bA_{1}(\bvrho) \pd{\bvrho}{x_{1}} + \bA_{2}(\bvrho) \pd{\bvrho}{x_{2}} = \bS(\bvrho),
\end{align}
where the left-hand side is derived from the left-hand side of the
Boltzmann equation, and the right-hand side corresponds to the
collision term. The detailed expression of the above system can be
found for the Boltzmann equation with the ES-BGK model in
\cite{Microflows1D} and with the Shakhov model in \cite{Li}. Here, we
omit it since it has no effect on the description of the solver
presented in this work.

Nevertheless, it is worth mentioning that the moment system
\cref{eq:moment-system} of order $M$ with $M\geq 2$ contains the
classical hydrodynamic equations. In particular, the equations with
the multi-dimensional indices $\alpha=0$ and $\alpha=e_{i}$,
$i=1,2,3$, are the conservation laws for mass and momentum
respectively. Therefore, the moment system of order $M$ is also known
as the extended hydrodynamic equations in the literature.

\subsection{Spatial finite volume discretization}
\label{sec:model-xdis}

It is tedious to discretize the moment system in spatial domain
directly based on the form \cref{eq:moment-system} for
$\bvrho$. Alternatively, in the framework of the numerical regularized
moment method originated in \cite{NRxx}, the spatial discretization of
\cref{eq:moment-system} is designed in a unified approach for the
moment system of arbitrary order, based on the underlying Boltzmann
equation \cref{eq:boltzmann} with the truncated approximation
\cref{eq:truncated-dis}. Our discretization given below is obtained
following this approach.

In the rest of the paper, the spatial coordinates $x_{1}$ and $x_{2}$
as well as the associated notations will be replaced by $x$ and $y$,
respectively, for convenience. Let the spatial domain $\Omega$ be
simply a rectangle with the length of $L_{x}$ and the height of
$L_{y}$. Suppose $\Omega$ is divided into a rectangular grid with
$N_{x}\times N_{y}$ cells, for which the grid points are given by
$(x_{i+1/2},y_{j+1/2})$, $i=0,1,\ldots,N_{x}$,
$j=0,1,\ldots,N_{y}$. Let $f_{ij}(t,\bxi)$ represent the average
distribution function over the $(i,j)$th grid cell
$[x_{i-1/2},x_{i+1/2}]\times [y_{j-1/2},y_{j+1/2}]$,
$i=1,2,\ldots,N_{x}$, $j=1,2,\ldots,N_{y}$. Then by applying the
finite volume method to the Boltzmann equation \cref{eq:boltzmann}, we
can obtain an ODE for the evolution of the average distribution
function over each grid cell as
\begin{align}
  \label{eq:ode-system}
  \opd{f_{ij}}{t} = - \left[ \frac{F_{i+\frac{1}{2},j} - F_{i-\frac{1}{2},j}}{\Delta x_{i}} + \frac{F_{i,j+\frac{1}{2}} - F_{i,j-\frac{1}{2}}}{\Delta y_{j}}\right] + Q(f_{ij}) =: \mR_{ij}(\bbf),
\end{align}
where $\Delta x_{i} = x_{i+1/2}-x_{i-1/2}$ and
$\Delta y_{j} = y_{j+1/2}-y_{j-1/2}$ are the length and height,
respectively, of the $(i,j)$th cell, $Q(f_{ij})$ is the average of
collision term over the $(i,j)$th cell, $\bbf$ is the vector composed
of all $f_{ij}$, and $F_{i+1/2,j}$ represents the numerical flux
defined at the boundary between the $(i,j)$th and $(i+1,j)$th
cells. Other numerical fluxes are defined similarly, so they are
omitted below.

With the assumption that all
$f_{ij}(t,\bxi) \in \mF_{M}^{[\bu_{ij},\theta_{ij}]}$, where
$\bu_{ij}(t)$ and $\theta_{ij}(t)$ are the averages of mean velocity
and temperature, respectively, on the $(i,j)$th cell and the relation
\cref{eq:moments-relation} holds for the expansion coefficients
$f_{ij,\alpha}(t)$, let us investigate the approximation of
\cref{eq:ode-system} over the $(i,j)$th cell in
$\mF_{M}^{[\bu_{ij},\theta_{ij}]}$. To this end, the four numerical
fluxes and $Q(f_{ij})$ are first projected into
$\mF_{M}^{[\bu_{ij},\theta_{ij}]}$, so that, for examples, we have
formally that
\begin{align}
  \label{eq:flux-force-expansion}
  \begin{aligned}
    & F_{i+\frac{1}{2},j}\approx \sum_{|\alpha|\leq M} F_{\alpha}(f_{i+\frac{1}{2},j}^{-}, f_{i+\frac{1}{2},j}^{+}) \mH_{\alpha}^{[\bu_{ij},\theta_{ij}]}(\bxi),\\
    & Q(f_{ij}) \approx \sum_{|\alpha|\leq M} Q_{ij,\alpha} \mH_{\alpha}^{[\bu_{ij},\theta_{ij}]}(\bxi),
  \end{aligned}
\end{align}
where the coefficients $F_{\alpha}$ depend on the left and right
limits of distribution function, i.e.,
$f_{i+1/2,j}^{\pm} \in \mF_{M}^{[\bu_{i+1/2,j}^{\pm},
  \theta_{i+1/2,j}^{\pm}]}$, at the right boundary of the $(i,j)$th
cell, and the coefficients $Q_{ij,\alpha}$ depend only on the
$(i,j)$th distribution function $f_{ij}$. If the numerical flux is
designed specially to reflect the hyperbolicity of the moment system
\cref{eq:moment-system}, then after matching the coefficients in both
sides of \cref{eq:ode-system} for the same basis function
$\mH_{\alpha}^{[\bu_{ij},\theta_{ij}]}(\bxi)$, we will get a system
which can equivalently be viewed as a spatial discretization of the
moment system \cref{eq:moment-system} on the $(i,j)$th
cell.

In our numerical experiments, the numerical flux introduced in
\cite{Microflows1D} is adopted. The detailed expression will not be
given here for brevity. However, it is noted that in order to calculate
the approximation of the numerical flux in
$\mF_{M}^{[\bu_{ij},\theta_{ij}]}$ as \cref{eq:flux-force-expansion},
a transformation to project a function from one space
$\mF_{M}^{[\bvpi,\vartheta]}$ into another space
$\mF_{M}^{[\tilde{\bvpi},\tilde{\vartheta}]}$ would be heavily
involved, since the expansion of the distribution function $f_{ij}$ on each cell
usually uses the basis functions with different parameters, and so are
the expansions of $f_{i+1/2,j}^{\pm}$. Such a transformation is in
fact one of cores of the unified discretization for the moment system
of arbitrary order. A fast algorithm with the time complexity of
$\mO(M^{3})$ for it has been implemented in \cite{NRxx, Qiao}.
For other parts of the present solver, this transformation will also
be employed frequently, and may not be explicitly pointed out.

It remains to determine the distribution functions $f_{i+1/2,j}^{\pm}$
at all cell boundaries. Obviously, it is equivalent to directly
determine the parameters $\bu_{i+1/2,j}^{\pm}$,
$\theta_{i+1/2,j}^{\pm}$, and all expansion coefficients
$f_{i+1/2,j,\alpha}^{\pm}$ with $|\alpha|\leq M$. Taking the
distribution functions on the left and right boundaries of the
$(i,j)$th cell, i.e., $f_{i-1/2,j}^{+}$ and $f_{i+1/2,j}^{-}$,
respectively, as examples, let us consider the reconstruction of them
providing that the average distribution functions
$f_{\iota j}\in \mF_{M}^{[\bu_{\iota j},\theta_{\iota j}]}$ on the
$(\iota,j)$th cell, $\iota=i, i\pm 1$, are given. Using the same
strategy as in \cite{hu2019efficient}, the corresponding parameters
and expansion coefficients of $f_{i-1/2,j}^{+}$ and $f_{i+1/2,j}^{-}$
could be computed efficiently by avoiding the transformation mentioned
above. To be specific, we have
\begin{align}
  \label{eq:moments-recon}
  \begin{aligned}
  \bu_{i-1/2,j}^{+}= \bu_{ij} - \frac{\Delta x_{i}}{2} \bg_{ij}, & \qquad 
\bu_{i+1/2,j}^{-} = \bu_{ij} + \frac{\Delta x_{i}}{2} \bg_{ij}, \\ 
  \theta_{i-1/2,j}^{+}= \theta_{ij} - \frac{\Delta x_{i}}{2} g_{ij}, & \qquad
\theta_{i+1/2,j}^{-} = \theta_{ij} + \frac{\Delta x_{i}}{2} g_{ij}, \\
  f_{i-1/2,j,\alpha}^{+}= f_{ij,\alpha} - \frac{\Delta x_{i}}{2} g_{ij,\alpha}, & \qquad
f_{i+1/2,j,\alpha}^{-} = f_{ij,\alpha} + \frac{\Delta x_{i}}{2} g_{ij,\alpha}, 
\end{aligned}
\end{align}
where $\bg_{ij}$, $g_{ij}$ and $g_{ij, \alpha}$ are reconstructed
slopes of the corresponding variables in the horizontal direction of
the $(i,j)$th cell. If all slopes are set to $0$, then the resulting
discretization would have the first-order spatial
accuracy. Accordingly, a second-order spatial discretization can be
obtained by computing these slopes as follows
\begin{align*}
\bg_{ij} = \frac{\bu_{i+1,j} - \bu_{i-1,j}}{x_{i+1} - x_{i-1}},  \qquad
g_{ij} = \frac{\theta_{i+1,j} - \theta_{i-1,j}}{x_{i+1} - x_{i-1}},  \qquad
g_{ij,\alpha} = \frac{f_{i+1,j,\alpha} - f_{i-1,j,\alpha}}{x_{i+1} - x_{i-1}}, 
\end{align*}
where $x_{i} = (x_{i+1/2} + x_{i-1/2}) / 2$ is the horizontal
coordinate of the center of the $(i,j)$th cell.

Now it is time to compute the coefficients $Q_{ij, \alpha}$ in
\cref{eq:flux-force-expansion} for the given distribution function
$f_{ij}$. For the BGK-type relaxation models, $Q_{ij,\alpha}$ can be
evaluated analytically with the time complexity of $\mO(M^{3})$, as shown
in \cite{Microflows1D} for the ES-BGK model and \cite{Li} for the
Shakhov model. By contrast, it is noted that for the quadratic
collision term \cref{eq:binary-collision}, an accurate computation of
$Q_{ij,\alpha}$, proposed in \cite{Wang2019}, takes the time
complexity of $\mO(M^{9})$.

Finally, proper boundary conditions are necessary for the spatial
discretization of the moment system on $\partial \Omega$. In our
experiments, the boundary conditions derived from the Maxwell boundary
condition would be employed. We refer to \cite{Li,Microflows1D} for
more details on such boundary conditions.

\subsection{Temporal discretization and time-integration scheme}
\label{sec:model-tdis}

The system of ODEs \cref{eq:ode-system} together with the boundary
conditions can be solved by the general time-integration methods. The
simplest time-integration scheme is the forward Euler scheme, which
can be formally written as
\begin{align}
  \label{eq:euler-scheme}
  f_{ij}^{n+1} = f_{ij}^{n} + \Delta t \mR_{ij}(\bbf^{n}), 
\end{align}
for all $i$ and $j$, where the superscript $n$ is used to denote the
approximation of the variables at time $t_{n}$, and
$\Delta t= t_{n+1}-t_{n}$ is the time step size.
According to the CFL condition, the time step size is chosen as
\begin{align}
  \label{eq:global-CFL}
  \Delta t = \min_{ij} \Delta t_{ij},
\end{align}
where the local time step size $\Delta t_{ij}$ must satisfy
\begin{align}
  \label{eq:local-CFL}
  \Delta t_{ij} \left( \frac{|u_{ij,1}^{n}|+C_{M+1} \sqrt{\theta_{ij}^{n}}}{\Delta x_{i}} + \frac{|u_{ij,2}^{n}|+C_{M+1} \sqrt{\theta_{ij}^{n}}}{\Delta y_{j}}\right) < 1,
\end{align}
in which $u_{ij,1}^{n}$ and $u_{ij,2}^{n}$ are the first two
components of the vector $\bu_{ij}^{n}$, and $C_{M+1}$ is the maximal
root of the Hermite polynomial of degree $M+1$.

As pointed out in \cite{hu2016acceleration}, the scheme
\cref{eq:euler-scheme} numerically consists of two steps. First, find
the approximation of the right-hand side in
$\mF_{M}^{[\bu_{ij}^{n},\theta_{ij}^{n}]}$ as an intermediate
distribution function $f_{ij}^{*}$. Next, evaluate the mean velocity
$\bu_{ij}^{n+1}$ and temperature $\theta_{ij}^{n+1}$ from $f_{ij}^{*}$
upon \cref{eq:rho-u-theta}, and then project $f_{ij}^{*}$ into
$\mF_{M}^{[\bu_{ij}^{n+1}, \theta_{ij}^{n+1}]}$ to get the approximate
distribution function $f_{ij}^{n+1}$ at time $t_{n+1}$.

It is well known that the forward Euler scheme \cref{eq:euler-scheme}
has first-order accuracy in time. When the second-order spatial
discretization is taken into account, the time-integration scheme has
to be modified accordingly into a second-order method to match the
spatial accuracy and to improve the stability of the method. As an
example, Heun's method can be employed for this purpose. Actually,
Heun's method has been introduced in \cite{hu2019efficient} to solve
the system of ODEs \cref{eq:ode-system} with second-order spatial
discretization. It is also noted that the computational cost in each
time step will be doubled for Heun's method in comparison to the
forward Euler scheme.

In the current work, we are interested in the steady state of the
rarefied gas cavity flows, when the time $t\to \infty$. Consequently,
the time-integration schemes always take a long time simulation to
achieve the steady state. In order to improve the efficiency of the
steady-state computation, we may give up time accuracy and try to
develop some fast convergent iteration methods using the acceleration
techniques introduced in the next section.

\section{Acceleration methods and steady-state solver}
\label{sec:solver}

We focus on the steady-state computation of cavity flows directly in
this section. At first, the discrete steady-state problem is rewritten as
\begin{align}
  \label{eq:local-steady-discretization}
  \mR_{ij}(\bbf) = 0, 
\end{align}
for all $i$ and $j$, by dropping the derivatives with respect to $t$
in \cref{eq:ode-system}. In this situation, the forward Euler scheme
\cref{eq:euler-scheme} becomes a simplest iteration to find the
steady-state solution. Based on this iteration, below three
acceleration methods, including fast sweeping iteration, nonlinear
multigrid method, and OpenMP-based parallelization, will be described
in detail to establish the efficient steady-state solver.

\subsection{Fast sweeping iteration}
\label{sec:solver-fs}

The forward Euler scheme \cref{eq:euler-scheme} is also referred to as
the Richardson iteration for steady-state problem
\cref{eq:local-steady-discretization}. It is essentially a Jacobi-type
iteration. A natural but often effective strategy to improve the
convergence rate is to modify the Jacobi-type iteration to a
Gauss-Seidel sweeping iteration. Precisely speaking, in the
Gauss-Seidel sweeping iteration, the grid cells are swept
cell-by-cell, and the newest approximation of $f_{ij}$ would be
utilized in the computation of the local residual $\mR_{ij}(\bbf)$ as
soon as they are available, in contrast to the forward Euler scheme
\cref{eq:euler-scheme}, which uses $f_{ij}^{n}$ to compute all local
residuals. Consequently, the resulting scheme to update the
approximation of $f_{ij}$ on the $(i,j)$th cell reads
\begin{align}
  \label{eq:gs-scheme}
  f_{ij}^{n+1} = f_{ij}^{n} + \Delta t_{ij} \mR_{ij}(\bbf^{*}),
\end{align}
where the vector $\bbf^{*}$ initially equals $\bbf^{n}$, and its
$(i,j)$th component will be replaced immediately by $f_{ij}^{n+1}$
after obtaining it. In the above scheme, the local time step size
$\Delta t_{ij}$ is adopted to replace the global time step size
$\Delta t$, since the former one is evidently more suitable than the
latter one for the Gauss-Seidel sweeping iteration. Moreover, such a
slight substitution of the time step size might speed up the
steady-state computation significantly, especially when the
non-uniform grid is used \cite{imamura2005acceleration}.

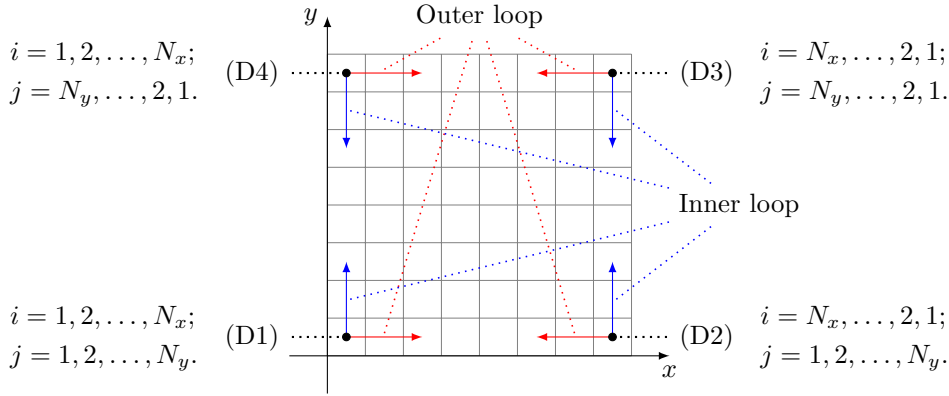
\begin{figure}[!htbp]
  \centering
    \begin{tikzpicture}[scale=1.]
    \draw[step=0.5, help lines] (0.5,0.5) grid(4.5,4.5);
    \draw[-latex] (0.5,0) -- (0.5,5) node[left] {$y$};
    \draw[-latex] (0,0.5) -- (5,0.5) node[below] {$x$};
    \draw[-latex,red] (0.75,0.75) -- (1.75,0.75);
    \draw[-latex,blue] (0.75,0.75) -- (0.75,1.75);
    \draw[dotted,thick] (0.75,0.75) -- (0,0.75) node[left]{$
      \begin{aligned}
        & i = 1,2,\ldots,N_{x};\\
        & j = 1,2,\ldots,N_{y}.
      \end{aligned}
      \quad (\text{D}1)$};
    \draw [fill] (0.75,0.75) circle [radius=0.05];
    \draw[-latex,red] (4.25,0.75) -- (3.25,0.75);
    \draw[-latex,blue] (4.25,0.75) -- (4.25,1.75);
    \draw[dotted,thick] (4.25,0.75) -- (5,0.75) node[right]{$(\text{D}2) \quad
      \begin{aligned}
        & i = N_{x},\ldots,2,1; \\ &j = 1,2,\ldots,N_{y}.
      \end{aligned}$};
    \draw [fill] (4.25,0.75) circle [radius=0.05];
    \draw[-latex,red] (4.25,4.25) -- (3.25,4.25);
    \draw[-latex,blue] (4.25,4.25) -- (4.25,3.25);
    \draw[dotted,thick] (4.25,4.25) -- (5,4.25) node[right]{$(\text{D}3) \quad
      \begin{aligned}
        & i = N_{x},\ldots,2,1; \\ & j = N_{y},\ldots,2,1.
      \end{aligned}$};
    \draw [fill] (4.25,4.25) circle [radius=0.05];
    \draw[-latex,red] (0.75,4.25) -- (1.75,4.25);
    \draw[-latex,blue] (0.75,4.25) -- (0.75,3.25);
    \draw[dotted,thick] (0.75,4.25) -- (0,4.25) node[left]{$
      \begin{aligned}
        & i = 1,2,\ldots,N_{x};\\ & j = N_{y},\ldots,2,1.
      \end{aligned}
      \quad (\text{D}4)$};
    \draw [fill] (0.75,4.25) circle [radius=0.05];
    \draw (2.5,5) node[text centered] (oloop) {Outer loop};
    \draw[dotted,red,semithick] (1.25,0.8) -- (oloop);
    \draw[dotted,red,semithick] (3.75,0.8) -- (oloop);
    \draw[dotted,red,semithick] (1.25,4.3) -- (oloop);
    \draw[dotted,red,semithick] (3.75,4.3) -- (oloop);
    \draw (5,2.5) node[right] (iloop) {Inner loop};
    \draw[dotted,blue,semithick] (0.8,1.25) -- (iloop);
    \draw[dotted,blue,semithick] (4.3,1.25) -- (iloop);
    \draw[dotted,blue,semithick] (0.8,3.75) -- (iloop);
    \draw[dotted,blue,semithick] (4.3,3.75) -- (iloop);
  \end{tikzpicture}
  \caption{The sketch of alternating direction sweepings.}
  \label{fig:sweeping}
\end{figure}

The sweeping direction is also important for the Gauss-Seidel sweeping
iteration. In one-dimensional case, the symmetric Gauss-Seidel
iteration always works well for the steady-state moment system, as can
be seen in \cite{hu2016acceleration,hu2015}. For multi-dimensional
cases, alternating direction sweepings are frequently adopted in
conjunction with the Gauss-Seidel iteration to give a desirable
iteration. In this paper, four alternating direction sweepings, that
are the same as in \cite{Zhang2016}, are taken into consideration to
construct the fast sweeping iteration. These four alternating
direction sweepings can be simply represented by
\begin{align*}
  \begin{aligned}
    & (\text{D}1): \quad i = 1,2,\ldots,N_{x} \text{ (outer loop)}; \quad j = 1,2,\ldots,N_{y} \text{ (inner loop)}; \\
    & (\text{D}2): \quad i = N_{x},\ldots,2,1 \text{ (outer loop)}; \quad j = 1,2,\ldots,N_{y} \text{ (inner loop)}; \\
    & (\text{D}3): \quad i = N_{x},\ldots,2,1 \text{ (outer loop)}; \quad j = N_{y},\ldots,2,1 \text{ (inner loop)}; \\
    & (\text{D}4): \quad i = 1,2,\ldots,N_{x} \text{ (outer loop)}; \quad j = N_{y},\ldots,2,1 \text{ (inner loop)}. \\
    \end{aligned}  
\end{align*}
The sketch of them is shown in \cref{fig:sweeping}.

Applying the above four alternating direction sweepings sequentially
to the Gauss-Seidel sweeping scheme \cref{eq:gs-scheme}, we get one
step of the fast sweeping iteration, which is summarized in
\Cref{alg:fs}, and denoted by $\bbf^{n+1}= \text{FS}(\bbf^{n},\br)$
throughout the paper. Here, $\br$ is introduced to represent the
vector of the right-hand side of
\cref{eq:local-steady-discretization}. Its default value is $0$
everywhere. Yet $\br$ might be nonzero, when the fast sweeping
iteration is utilized in the multigrid solver given in the next
subsection.

\begin{algorithm}[!htbp]
  \renewcommand{\algorithmicrequire}{\textbf{Input:}}
  \renewcommand{\algorithmicensure}{\textbf{Output:}}
  \caption{One step of the fast sweeping iteration}\label{alg:fs}
  \begin{algorithmic}[1]
    \REQUIRE Grid, $\bbf^{n}$, and $\br$ (the vector of the right-hand side of \cref{eq:local-steady-discretization}, default: $\br=0$)
    \ENSURE The new approximation $\bbf^{n+1}$ denoted by $\bbf^{n+1}=\text{FS}(\bbf^{n},\br)$
    \vspace*{0.7em}
    \STATE $\tilde{\bbf} \leftarrow \bbf^{n}$;
    \FOR {$k$ from $1$ \TO $4$}
    \FOR[$(\text{D}k)$ is the sweeping directions shown in \cref{fig:sweeping}] {$i,j \in (\text{D}k)$}
    \STATE Determine $\Delta t_{ij}$ from \cref{eq:local-CFL} with $\tilde{f}_{ij}$;
    \STATE $\tilde{f}_{ij} \leftarrow \tilde{f}_{ij} + \Delta t_{ij} \left(\mR_{ij}(\tilde{\bbf})-r_{ij}\right)$; \COMMENT{Modified from \cref{eq:gs-scheme}}
    \ENDFOR
    \ENDFOR
    \RETURN $\bbf^{n+1} \leftarrow \tilde{\bbf}$;
  \end{algorithmic}
\end{algorithm}

As mentioned in \cite{Zhang2016}, we are able to observe that the
present fast sweeping iteration is several times faster than the
forward Euler scheme. Indeed, it is found in our numerical experiments
that the fast sweeping iteration takes drastically fewer number of
iterations to achieve the steady state in almost all cases. Apart from
this, the fast sweeping iteration turns out to be more stable than the
time-integration schemes. The CFL number to determine the time step
size could be larger than that in time-integration schemes to further
improve the convergence. Especially when a second-order spatial
discretization is applied, it is easy to see the situation that the fast
sweeping iteration converges while the forward Euler scheme fails to
converge, for a relatively large CFL number.

We conclude this subsection by pointing out that the convergence
behavior of the fast sweeping iteration with respect to the number of
grid cells is almost the same as the time-integration schemes, that
is, the number of iterations to achieve the steady state will be
doubled as both $N_{x}$ and $N_{y}$ are doubled. Aiming to improve
this behavior, we are going to consider the multigrid method in the
following subsection.

\subsection{Nonlinear  multigrid solver}
\label{sec:solver-mg}

The system \cref{eq:local-steady-discretization} that we intend to
solve is evidently a complicated nonlinear system. For a nonlinear
problem, there are mainly two approaches to develop multigrid method
\cite{hackbusch1985book}. One is the global-linearization-based method
such as the Newton-multigrid method, the other is the intrinsic
nonlinear multigrid (NMG) method, known also as the full approximation
storage (FAS) method
\cite{brandt1977multi-level,brandt2011book}. Based on our experience,
we recommend using the NMG method here, because efficient
implementation of global linearization usually depends on the specific
numerical flux, and is nontrivial especially for the system
\cref{eq:local-steady-discretization} with a relatively large order
$M$. Moreover, a great advantage of the NMG method is that the
previous fast sweeping iteration as well as the time-integration
schemes can be reused directly as the smoother, in the framework of
the NMG method.

Besides the smoother, the remaining key ingredients of the NMG method,
that need to be further specified for the given problem, include the
coarse grid correction and the transfer operators between two adjacent
levels of grids, namely, restriction and prolongation. In this
subsection, we first utilize two levels of grids to illustrate these
components, and then present the complete NMG solver by recursion.

\subsubsection{A nonlinear two-grid  iteration}
\label{sec:solver-mg-coarse}

For convenience, let us denote operators and variables related to the
fine and coarse grids, respectively, by subscripts $h$ and $H$. Then,
the underlying problem on the fine grid can be
written simply as
\begin{align}
  \label{eq:fine-grid-problem}
  \mR_{h}(\bbf_{h}) = \br_{h}, 
\end{align}
where the left-hand side $\mR_{h}$ is the global form of the
discretization operator $\mR_{ij}$ on the fine grid, and the
right-hand side $\br_{h}$ is a known vector independent of the unknown
solution $\bbf_{h}$. In particular, we have $\br_{h}\equiv 0$ for the
original steady-state problem \cref{eq:local-steady-discretization}.

According to the framework of the NMG method, the fine grid problem
\cref{eq:fine-grid-problem} is basically solved by simple iterative
methods, which are able to damp the high-frequency error of a given
approximation of the solution quickly in a few number of
iterations. This procedure is called smoothing and the corresponding
iteration is known as the smoother. Usually, the previous fast
sweeping iteration as well as the time-integration schemes fulfills
our requirement. Hence we utilize the fast sweeping iteration as the
smoother in the current work.

Given an initial approximation of the solution for the fine grid
problem \cref{eq:fine-grid-problem}, suppose after several
pre-smoothing steps we get the numerical solution denoted by
$\bar{\bbf}_{h}$, in which the $(i,j)$th component
$\bar{f}_{h,ij}(\bxi) \in \mF_M^{[\bar{\bu}_{h,ij},
  \bar{\theta}_{h,ij}]}$, and $\bar{\bu}_{h,ij}$,
$\bar{\theta}_{h,ij}$ are the corresponding velocity and temperature
determined from $\bar{f}_{h,ij}$. In order to efficiently damp the
low-frequency components of the error, both the numerical solution
$\bar{\bbf}_{h}$ and the residual
$\bar{\bR}_{h} = \br_{h} - \mR_{h}(\bar{\bbf}_{h})$ on the fine grid
would be then restricted into the coarse grid, to formulate the coarse
grid problem as
\begin{align}
  \label{eq:coarse-grid-problem}
  \mR_{H}(\bbf_{H}) = \br_{H} := \mR_{H}(I_{h}^{H}\bar{\bbf}_{h}) + I_{h}^{H} \bar{\bR}_{h}, 
\end{align}
where $\mR_{H}$ is the discretization operator on the coarse grid
defined analogously to the fine grid counterpart $\mR_{h}$, and
$I_{h}^{H}$ is the restriction operator that will be illustrated in
\cref{sec:solver-operations}. It is obvious that the coarse grid
problem \cref{eq:coarse-grid-problem} could be solved by the same
strategy as the fine grid problem \cref{eq:fine-grid-problem}. In
practice, we adopt $I_{h}^{H}\bar{\bbf}_{h}$ as the initial guess for
the solution of the coarse grid problem \cref{eq:coarse-grid-problem}.

When a new approximation of the solution for the coarse grid problem
\cref{eq:coarse-grid-problem}, denoted by $\tilde{\bbf}_{H}$, is
obtained, we can calculate the correction on the coarse grid and update
the fine grid solution $\bar{\bbf}_{h}$ to $\hat{\bbf}_{h}$ by
\begin{align}
  \label{eq:update-fine-grid-solution}
  \hat{\bbf}_h = \bar{\bbf}_h + I_H^h\left( \tilde{\bbf}_H - I_h^H \bar{\bbf}_h \right),
\end{align}
where $I_H^h$ is the prolongation operator transferring functions from
the coarse grid to the fine grid, and will be given in the next
subsection. Finally, taking $\hat{\bbf}_{h}$ as the initial value,
several post-smoothing steps would be applied for the fine grid
problem \cref{eq:fine-grid-problem}. This completes a single step of
the nonlinear two-grid iteration.

\subsubsection{Restriction and prolongation operators}
\label{sec:solver-operations}
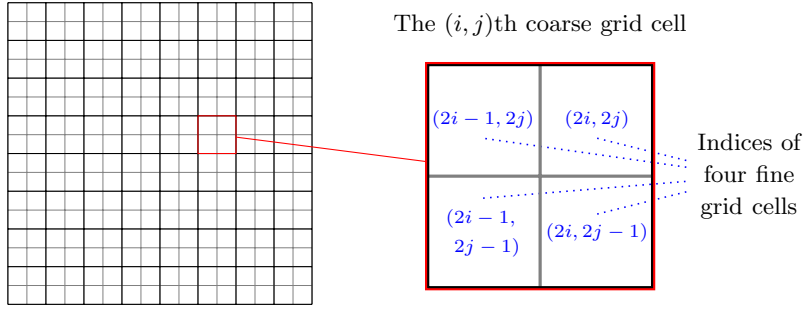
\begin{figure}[!htbp]
  \centering
    \begin{tikzpicture}
      [spy using outlines={rectangle, red, magnification=6, size=3cm, connect spies}]
      \draw[step=0.25, help lines] (0,0) grid(4,4);
      \draw[step=0.5] (0,0) grid(4,4);
      
      \spy on (2.75,2.25) in node (zoom) at (7,1.7);
      \draw (7,1.7) node (zc) {};
      \draw ($(zc)+(0,2.)$) node[text centered] {\small The $(i,j)$th coarse grid cell};
      \draw ($(zc)+(0.75,0.75)$) node[blue,text centered] (f1) {\scriptsize $(2i,2j)$};
      \draw ($(zc)+(-0.75,0.75)$) node[blue,text centered] (f2) {\scriptsize $(2i-1,2j)$};
      \draw ($(zc)+(-0.75,-0.75)$) node[blue,text centered] (f3) {\scriptsize $
        \begin{aligned}
          (& 2i-1,\\ \scriptsize &2j-1)
        \end{aligned}$};
      \draw ($(zc)+(0.75,-0.75)$) node[blue,text centered] (f4) {\scriptsize $(2i,2j-1)$};
      \draw ($(zc)+(2.75,0)$) node[text centered, text width=40] (fine) {\small Indices of four fine grid cells};
      \draw[dotted,blue,semithick] (f1.south) -- (fine) -- (f2.south); \draw[dotted,blue,semithick] (f3.north) -- (fine) -- (f4.north);
    \end{tikzpicture}
    \caption{Geometric relationship between two adjacent levels of grids.}
    \label{fig:2-grid}
\end{figure}

To construct the restriction operator $I_{h}^{H}$ and the prolongation
operator $I_{H}^{h}$ in detail, it is helpful to first clarify the
generation of the multi-levels of grids. Recalling that in this paper
the spatial domain $\Omega$ is assumed to be a rectangle and be
discretized by a rectangular grid, the coarse grid can be easily
generated in a standard way by merging cells of the given fine
grid. As shown in \cref{fig:2-grid}, we have that the $(i,j)$th coarse
grid cell is precisely composed of four fine grid cells with indices
given by $(2i-1,2j-1)$, $(2i,2j-1)$, $(2i-1,2j)$ and $(2i,2j)$,
respectively. With this geometric relationship, both the restriction
operator and the prolongation operator can be constructed locally
between the $(i,j)$th coarse grid cell and the corresponding four fine
grid cells for the cell-centered function that adopted in our
discretization. In particular, they are established here
by following the same strategy as employed in \cite{hu2014nmg} for the
implementation of one-dimensional versions of these operators. For
simplicity of presentation, the four fine grid indices $(2i-1,2j-1)$,
$(2i,2j-1)$, $(2i-1,2j)$ and $(2i,2j)$ shall be replaced by the local
indices $1$, $2$, $3$, $4$, respectively, below without causing
confusion.

We can see from the right-hand side of \cref{eq:coarse-grid-problem}
that there are two fine grid variables, i.e., the numerical solution
$\bar{\bbf}_{h}$ and the corresponding residual $\bar{\bR}_{h}$, which
are required to be restricted into the coarse grid. The resulting
restrictions are given by $\bar{\bbf}_{H} = I_{h}^{H}\bar{\bbf}_{h}$
and $\bar\bR_{H} = I_{h}^{H} \bar{\bR}_{h}$, respectively. Before
giving their components $\bar{f}_{H,ij}$ and $\bar{R}_{H,ij}$ on the
$(i,j)$th coarse grid cell, it is worth reminding that on the
$\iota$th fine grid cell, $\iota=1,2,3,4$, the numerical solution
$\bar{f}_{h,\iota} \in \mF_M^{[\bar{\bu}_{h,\iota},
  \bar{\theta}_{h,\iota}]}$ is a linear combination of the basis
functions
$\mH_{\alpha}^{[\bar{\bu}_{h,\iota}, \bar{\theta}_{h,\iota}]}(\bxi)$
with $|\alpha|\leq M$. Thereby, the corresponding residual
$\bar{R}_{h,\iota}$ is also calculated in
$\mF_M^{[\bar{\bu}_{h,\iota}, \bar{\theta}_{h,\iota}]}$ in our
implementation, which means $\bar{R}_{h,\iota}$ is expressed as a
linear combination of the same basis functions as $\bar{f}_{h,\iota}$
too.

In order to get the expressions of $\bar{f}_{H,ij}$ and
$\bar{R}_{H,ij}$, which are expected to belong to the same function
space $\mF_M^{[\bar{\bu}_{H,ij}, \bar{\theta}_{H,ij}]}$, we have first
to determine the parameters $\bar\bu_{H,ij}$ and
$\bar{\theta}_{H,ij}$.
Similar to \cite{hu2014nmg}, these two quantities can be computed from
the equations for conservation of mass, momentum, and energy, i.e.,
\begin{align}
  \label{eq:uH-thetaH-computation}
  \begin{aligned}
    & \bar\rho_{H, ij} \Delta s_{H,ij} = \sum_{\iota=1}^{4} \bar\rho_{h,\iota} \Delta s_{h,\iota},\\
    & \bar\rho_{H,ij}\bar\bu_{H,ij} \Delta s_{H,ij} = \sum_{\iota=1}^{4} \bar\rho_{h,\iota} \bar\bu_{h,\iota} \Delta s_{h,\iota},\\
    & \left(\bar\rho_{H,ij}\bar\bu_{H,ij}^{2} + 3\bar\rho_{H,ij}\bar\theta_{H,ij}\right) \Delta s_{H,ij} = \sum_{\iota=1}^{4} \left(\bar\rho_{h,\iota} \bar\bu_{h,\iota}^{2}+ 3\bar\rho_{h,\iota}\bar\theta_{h,\iota}\right) \Delta s_{h,\iota},
  \end{aligned}
\end{align}
where $\Delta s_{H,ij}$ is the area of the $(i,j)$th coarse grid cell,
$\Delta s_{h,\iota}$ is the area of the $\iota$th fine grid cell, and
$\bar\rho_{h,\iota}$ is the $\iota$th fine grid density extracted from
$\bar{f}_{h,\iota}$.

Now let us project the fine grid solution $\bar{f}_{h,\iota}$ into
$\mF_M^{[\bar{\bu}_{H,ij}, \bar{\theta}_{H,ij}]}$ to get a linear
combination of the basis functions
$\mH_{\alpha}^{[\bar{\bu}_{H,ij}, \bar{\theta}_{H,ij}]}(\bxi)$ with
$|\alpha|\leq M$. Denoting the resulting coefficients by
$\bar{f}^{*}_{h,\iota,\alpha}$, the coefficients of $\bar{f}_{H,ij}$
in the linear combination of the same basis functions can then be
evaluated by
\begin{align}
  \label{eq:fH}
  \bar{f}_{H,ij,\alpha} = \frac{1}{\Delta s_{H,ij}} \sum_{\iota=1}^{4} \bar{f}^{*}_{h,\iota,\alpha} \Delta s_{h,\iota},
\end{align}
which preserves the property of conservation well as the
one-dimensional case in \cite{hu2014nmg}. This completes the
construction of the restriction $\bar{f}_{H,ij}$. Similarly, we can
get the restriction $\bar{R}_{H,ij}$.

For the prolongation operator $I_{H}^{h}$, the simplest identity
operator is employed, that is, we have
$\left(I_{H}^{h}\bg_{H}\right)_{h,\iota} = g_{H,ij}$ for any coarse
grid quantity $\bg_{H}$. Consequently, the update
\cref{eq:update-fine-grid-solution} can be rewritten as 
\begin{align}
  \label{eq:rewritten-update}
  \hat{\bbf}_h = \bar{\bbf}_h - I_h^H \bar{\bbf}_h + \tilde{\bbf}_H.
\end{align}
It is pointed out that the above formula implicitly includes the
transformation between function spaces with different parameters,
since the three terms of the right-hand side actually have expressions
under different basis functions.

\subsubsection{Complete multigrid algorithm}
\label{sec:solver-mg-alg}

It is natural to extend the nonlinear two-grid iteration to a
nonlinear multigrid iteration by recursively applying the two-grid
strategy to the coarse grid problem \cref{eq:coarse-grid-problem},
until it can be solved efficiently by the single level method.

To complete the NMG algorithm on a given grid, a sequence of coarse
grids would first be generated by merging grid cells level by
level. Suppose the total levels of grids is $K+1$, and let us
introduce subscripts $h_{k}$, $k=0,1,\ldots,K$, to denote operators
and variables related to the $k$th-level grid, where $k=0$ and $K$
correspond to the coarsest and the finest grid, respectively. Then a
$(k+1)$-level NMG iteration for the $k$th-level problem
\begin{align}
  \label{eq:kth-level-problem}
  \mR_{h_{k}}(\bbf_{h_{k}}) = \br_{h_{k}},
\end{align}
to produce the new approximation of the solution from a given
approximation $\bbf_{h_{k}}^{n}$, can be summarized in \cref{alg:nmg}.

\begin{algorithm}[!htbp]
  \renewcommand{\algorithmicrequire}{\textbf{Input:}}
  \renewcommand{\algorithmicensure}{\textbf{Output:}}
  \renewcommand{\algorithmicloop}{\textbf{begin}}
  \renewcommand{\algorithmicendloop}{\algorithmicend}
  \caption{One step of $(k+1)$-level NMG iteration for \cref{eq:kth-level-problem}}\label{alg:nmg}
  \begin{algorithmic}[1]
    \REQUIRE Level index $k$, right-hand side $\br_{h_{k}}$, initial approximation $\bbf_{h_{k}}^{n}$
    \ENSURE The new approximation $\bbf_{h_{k}}^{n+1}$ denoted by $\bbf_{h_{k}}^{n+1}=\text{NMG}_{k}(\bbf_{h_{k}}^{n},\br_{h_{k}})$
    \vspace*{0.7em}
    \IF {$k=0$}
    \STATE Call the coarsest grid solver to obtain the new approximation $\bbf_{h_{0}}^{n+1}$;
    \ELSE
    \STATE \emph{Pre-smoothing}: perform $s_{1}$ steps of the fast sweeping iteration (\cref{alg:fs}) to obtain a new approximation $\bar\bbf_{h_{k}}$ by $\bar\bbf_{h_{k}} = \text{FS}^{s_{1}}(\bbf_{h_{k}}^{n}, \br_{h_{k}})$;
    \LOOP[\emph{Coarse grid correction}]
    \STATE Compute the fine grid residual as $\bar{\bR}_{h_{k}} = \br_{h_{k}} - \mR_{h_{k}}(\bar\bbf_{h_{k}})$;
    \STATE Construct the initial approximation of the coarse grid solution as\\
    \centerline{$\bar\bbf_{h_{k-1}} = I_{h_{k}}^{h_{k-1}} \bar\bbf_{h_{k}}$;}
    \STATE Calculate the right-hand side of the coarse grid problem \cref{eq:coarse-grid-problem} as\\
    \centerline{$\br_{h_{k-1}} = \mR_{h_{k-1}}(\bar\bbf_{h_{k-1}}) + I_{h_{k}}^{h_{k-1}} \bar\bR_{h_{k}}$;}
    \STATE Recursively call $\gamma$ steps of the $k$-level NMG iteration to obtain the new approximation of the coarse grid problem as\\
    \centerline{$\tilde\bbf_{h_{k-1}}=\text{NMG}_{k-1}^{\gamma}(\bar\bbf_{h_{k-1}},\br_{h_{k-1}})$;}
    \STATE Update the fine grid solution by \cref{eq:rewritten-update} as $\hat{\bbf}_{h_{k}} = \bar{\bbf}_{h_{k}} - \bar{\bbf}_{h_{k-1}} + \tilde{\bbf}_{h_{k-1}}$;
    \ENDLOOP
    \STATE \emph{Post-smoothing}: perform $s_{2}$ steps of the fast sweeping iteration to get the final approximation $\bbf_{h_{k}}^{n+1}$, i.e., $\bbf_{h_{k}}^{n+1} = \text{FS}^{s_{2}}(\hat\bbf_{h_{k}}, \br_{h_{k}})$;
    \ENDIF
  \end{algorithmic}
\end{algorithm}

It remains to give the coarsest grid solver in \cref{alg:nmg}. Since
the coarsest grid problem is analogous to the problem defined on the
other levels of grids, the fast sweeping iteration is applied again
for the coarsest grid solver. Moreover, for the sake of efficiency,
the fast sweeping iteration is performed at most $s_{3}$ steps instead
of being performed until convergence, in each calling of the coarsest
grid solver. Here, $s_{3}$ is a positive integer close to the
smoothing steps $s_{1} + s_{2}$. 

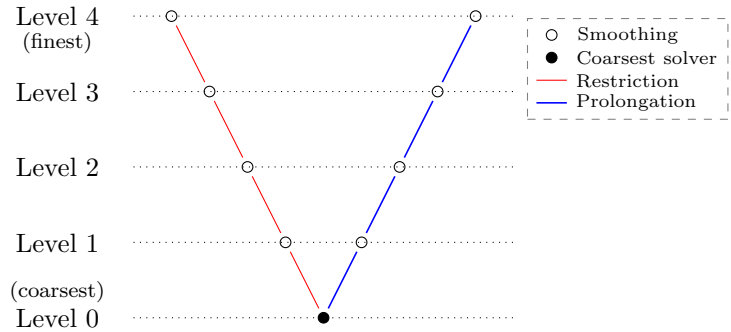
\begin{figure}[!htbp]
  \centering
      \begin{tikzpicture}
        \foreach \y/\ytext in {0, 1, 2, 3, 4}
        {
          \draw[dotted] (0,\y)--(5,\y);
          \draw (-1,\y) node [text centered] (l\y) {Level $\y$};
        }
        \draw ($(l0.north)+(0,0.1)$) node[text centered] {\footnotesize (coarsest)};
        \draw ($(l4.south)+(0,-0.1)$) node[text centered] {\footnotesize (finest)};
        \draw (2.5,0) node (c1) {}; \draw[fill] (c1) circle [radius=2pt]; 
        \foreach \x/\xtext in {2/1,3/2}
        {\draw (\x,1) node (c1x\xtext) {}; \draw (c1x\xtext) circle [radius=2pt];}
        \foreach \x/\xtext in {1.5/1,3.5/2}
        {\draw (\x,2) node (c2x\xtext) {}; \draw (c2x\xtext) circle [radius=2pt];}
        \foreach \x/\xtext in {1/1,4/2}
        {\draw (\x,3) node (c3x\xtext) {}; \draw (c3x\xtext) circle [radius=2pt];}
        \foreach \x/\xtext in {0.5/1,4.5/2}
        {\draw (\x,4) node (c4x\xtext) {}; \draw (c4x\xtext) circle [radius=2pt];}
        \draw[red] (c1) -- (c1x1) -- (c2x1) -- (c3x1) -- (c4x1);
        \draw[blue,semithick] (c1) -- (c1x2) -- (c2x2) -- (c3x2) -- (c4x2);
        \draw (6.5,3.3) node[gray,dashed,draw,minimum width=75pt,minimum height=38pt,rectangle] (legend) {};
      \draw ($(legend)+(-1.0,0.45)$) circle[radius=2pt]; \draw ($(legend)+(0.25,0.45)$) node {\scriptsize Smoothing\hspace*{1.8em}};
      \draw[fill] ($(legend)+(-1.0,0.15)$) circle[radius=2pt]; \draw ($(legend)+(0.25,0.15)$) node {\scriptsize Coarsest solver};
      \draw[red] ($(legend)+(-1.2,-0.15)$) -- ($(legend)+(-0.8,-0.15)$); \draw ($(legend)+(0.25,-0.15)$) node {\scriptsize Restriction\hspace*{1.8em}};      
      \draw[blue,semithick] ($(legend)+(-1.2,-0.45)$) -- ($(legend)+(-0.8,-0.45)$); \draw ($(legend)+(0.25,-0.45)$) node {\scriptsize Prolongation\hspace*{1.0em}};
    \end{tikzpicture}
    \caption{Diagram of a $V$-cycle NMG iteration.}
    \label{fig:nmg-v}
\end{figure}

Besides the coarsest grid solver, it is noted that the parameter
$\gamma$ in \cref{alg:nmg} is usually taken as $1$ or $2$,
corresponding to the so-called $V$-cycle or $W$-cycle NMG method,
respectively. Although, in practice, the $W$-cycle NMG method may have
better convergence rate than the $V$-cycle NMG method, the former one
takes much more computational cost for each iteration than the latter
one. Accordingly, we only report the numerical results of the
$V$-cycle NMG method in the current work. A diagram of a $V$-cycle
$5$-level NMG iteration is given in \cref{fig:nmg-v}.

We can now get a $(K+1)$-level NMG solver for the original
steady-state problem \cref{eq:local-steady-discretization} on a given
grid, by repeatedly performing the $(K+1)$-level NMG iteration, i.e.,
$\bbf_{h_{K}}^{n+1} = \text{NMG}_{K}(\bbf_{h_{K}}^{n}, 0)$, until the
steady state has been achieved. In particular, the one-level NMG
solver reduces to the single level solver of fast sweeping
iteration. Incidentally, we consider the steady state is achieved for
a given numerical solution $\bbf_{h_{K}}^{n}$ in our experiments, if
the resulting residual relative to the initial residual is less than a
given tolerance $\mathit{Tol}$, that is,
\begin{align}
  \label{eq:relative-residual}
  \frac{\Vert \mR_{h_{K}}(\bbf_{h_{K}}^{n}) \Vert}{\Vert \mR_{h_{K}}(\bbf_{h_{K}}^{0}) \Vert} \leq \mathit{Tol},
\end{align}
where the norm of a discrete residual $\mR(\bbf)$ is defined by
\begin{align*}
  \left\Vert \mR(\bbf) \right\Vert = \sqrt{\frac{1}{L_{x}L_{y}} \sum_{i=1}^{N_{x}}\sum_{j=1}^{N_{y}}\left\Vert \mR_{ij}(\bbf) \right\Vert^{2} \Delta s_{ij}},
\end{align*}
in which $\Delta s_{ij}$ is the area of the $(i,j)$th grid cell, and
the weighted $L^{2}$ norm of the space
$\mF_{M}^{[\bu_{ij},\theta_{ij}]}$ is employed to compute the local norm
$\left\Vert \mR_{ij}(\bbf) \right\Vert$ as \cite{hu2014nmg,hu2019efficient}.

\subsection{Parallelization of the solver with OpenMP}
\label{sec:solver-openmp}

The parallelization of an algorithm is able to fully utilize the
hardware resources of any computer with multi-core processors, so that
the simulation could be further accelerated. The OpenMP
\cite{dagum1998openmp} is a portable, scalable model that gives
programmers a simple and flexible interface for developing parallel
applications. It supports multi-platform shared-memory multiprocessing
programming in C/C++ and Fortran on many platforms ranging from the
standard desktop computer to the supercomputer. Therefore, we adopt
the OpenMP to parallelize the previous NMG solver, which is
implemented in C++ in our work.

In particular, the OpenMP executable directive is mainly applied to the
\emph{for} loop that traverses grid cells of a given grid. For almost
all operations in the NMG solver, such as smoothing, restriction,
prolongation, and residual computation, they can be parallelized
trivially in this way, no matter whether there is data racing.
Meanwhile, the parallel computational results shall keep the same as
the serial computational results except for the fast sweeping
iteration which is used in the modules of smoothing and coarsest grid
solver. It should be pointed out that the fast sweeping iteration,
i.e., \cref{alg:fs}, is essentially a serial algorithm due to the data
dependence. When applying the OpenMP-based parallelization to the
outer \emph{for} loop of \cref{alg:fs}, the grid is actually divided
into several continuous blocks, and the sweeping iteration is executed
in each continuous block. As a result, the parallel fast sweeping
iteration would be not exactly equivalent to the original serial fast
sweeping iteration. Nevertheless, it can be observed in all our
simulations that the convergence would be almost unaffected by the
parallel fast sweeping iteration.

\section{Numerical experiments}
\label{sec:num-ex}

A number of numerical experiments for three types of 2D square cavity
flows, i.e., single lid-driven flow, four-sided lid-driven flow, and
bottom-heated flow, are carried out to explore the main features of
the proposed NMG solver such as efficiency and robustness.  Numerical
results for the first-order spatial discretization are first
investigated in detail. Then the results of single lid-driven cavity
flow for the second-order spatial discretization are presented.

Throughout the experiments, we assume that all sides of the cavity are
completely diffusive, and the cavity is filled with the argon
gas. Thus we have the molecular mass $m_* = \SI{6.63e-26}{\kg}$, and
the Prandtl number $\Pr = 2/3$. The average collision frequency $\nu$
is taken to be
\begin{align}
  \label{eq:nu-shakhov}
  \nu = \beta \sqrt{\frac{\pi}{2}} \frac{1}{\Kn} \rho \theta^{1-w},
\end{align}
where $\Kn$ is the Knudsen number, $w$ is the viscosity index given by
$0.81$, and $\beta$ takes the value $1$ and $\Pr$ for the Shakhov
model and the ES-BGK model, respectively. 

In all simulations, the tolerance indicating the achievement of steady
state is set to $\mathit{Tol}=10^{-8}$, and the CFL number to
determine the time step size is set to $0.9$.
For the proposed NMG solver, the smoothing steps $s_{1} = s_{2} = 2$
and $s_{3}=4$ are employed. It has been observed that for most
situations, the efficiency of the NMG solver would be effectively
improved by increasing the total levels of grids, when the coarsest
grid still has a large number of grid cells. To be as efficient as
possible, after a number of preliminary simulations, the number of
total grid levels of the NMG solver is chosen in such a way that the
coarsest grid consists of $8\times 8$ cells.

Additionally, it is easy to show that the boundary conditions utilized
in the simulation could not determine a unique steady-state
solution. To recover the consistent steady-state solution with the
time-integration scheme, the correction adopted in
\cite{hu2014nmg} would be also applied at each NMG iteration in our simulation.

\subsection{Single lid-driven cavity flow}
\label{sec:num-ex-single-lid}

\begin{figure}[!htbp]
  \centering
  \begin{minipage}{1.\linewidth}
  \subfloat[Single lid-driven]{
    \begin{tikzpicture}[scale=0.9]
      \draw[black,line width =1pt] (0,0) rectangle (5,5);
      \draw[latex-latex,line width =0.5pt] (0,4.5)--(5,4.5);
      \draw[latex-latex,line width =0.5pt] (0.5,0)--(0.5,5);
      \draw[-latex] (1,4)--(1.5,4) node [right] {$x$};
      \draw[-latex] (1,4)--(1,3.5) node [below] {$y$};
      \node[] at (-0.25,-0.25) {A};
      \node[] at (-0.25,5.25) {D};
      \node[] at (5.25,5.25) {C};
      \node[] at (5.25,-0.25) {B};
      \node[] at (0.75,2) {$L_{y}$};
      \node[] at (3,4.25) {$L_{x}$};
      \node[] at (2.5,5.5) {$U_W$};
      \draw[-latex] (2,5.125)--(3,5.125);
    \end{tikzpicture}
    \label{fig:single-lid-geo}
  }\hfill
  \subfloat[Four-sided lid-driven]{
    \begin{tikzpicture}[scale=0.9]
      \draw[black,line width =1pt] (0,0) rectangle (5,5);
      \draw[latex-latex,line width =0.5pt] (0,4.5)--(5,4.5);
      \draw[latex-latex,line width =0.5pt] (0.5,0)--(0.5,5);
      \draw[-latex] (1,4)--(1.5,4) node [right] {$x$};
      \draw[-latex] (1,4)--(1,3.5) node [below] {$y$};
      \node[] at (-0.25,-0.25) {A};
      \node[] at (-0.25,5.25) {D};
      \node[] at (5.25,5.25) {C};
      \node[] at (5.25,-0.25) {B};
      \node[] at (0.75,2) {$L_{y}$};
      \node[] at (3,4.25) {$L_{x}$};
      \draw[-latex] (2,5.125)--(3,5.125) node at (2.5,5.5) {$U_W$};
      \draw[latex-] (2,0.125)--(3,0.125) node at (2.5,0.5) {$U_W$};
      \draw[-latex] (-0.125,3)--(-0.125,2) node at (-0.5,2.5) {$U_W$};
      \draw[latex-] (5.125,3)--(5.125,2) node at (5.5,2.5) {$U_W$};
    \end{tikzpicture}
    \label{fig:four-lid-geo}
  }
\end{minipage}
\caption{Configuration of the lid-driven cavity flow problems.}
  \label{fig:cavity_geometry}
\end{figure}
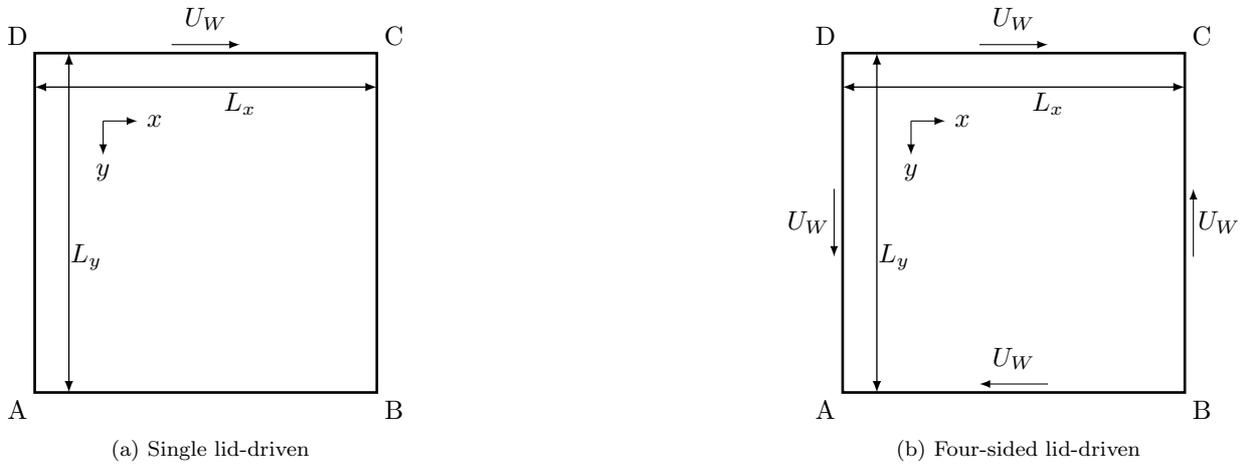

The single lid-driven cavity flow, whose configuration is shown in
\cref{fig:single-lid-geo}, is one of the most frequently used
benchmark tests in the multi-dimensional case. Our setting is the same
as in \cite{Qiao, Cai2018}. To be specific, the length and height of
the cavity are $L_{x} = L_{y} = \SI{9.63e-7}{\m}$. The top lid moves
horizontally to the right with a constant speed
$U_{W}=\SI[per-mode=symbol]{50}{\m\per\s}$ and is maintained at
temperature $\SI{273}{\kelvin}$, while the other sides of the cavity
are stationary and have the same temperature. Initially, the gas is
uniformly distributed and in the Maxwellian with constant density,
mean velocity of $0$, and temperature of $\SI{273}{\kelvin}$. Driven
by the motion of the top lid, the gas would finally reach a steady
state. Two initial densities given by
$\rho = \SI[per-mode=symbol]{0.891}{\kg \per \m^{3}}$ and
$\SI[per-mode=symbol]{0.0891}{\kg \per \m^{3}}$ are considered
below. They correspond to the Knudsen number $0.1$ and $1.0$,
respectively.

\subsubsection{Solution validation}
\label{sec:num-ex-single-lid-sol}

\begin{figure}[!htbp]
  \centering
  \subfloat[Temperature $(\si{\kelvin})$]{
    \label{fig:kn_pt1_shakhov_T}
    \includegraphics[width=0.47\textwidth]{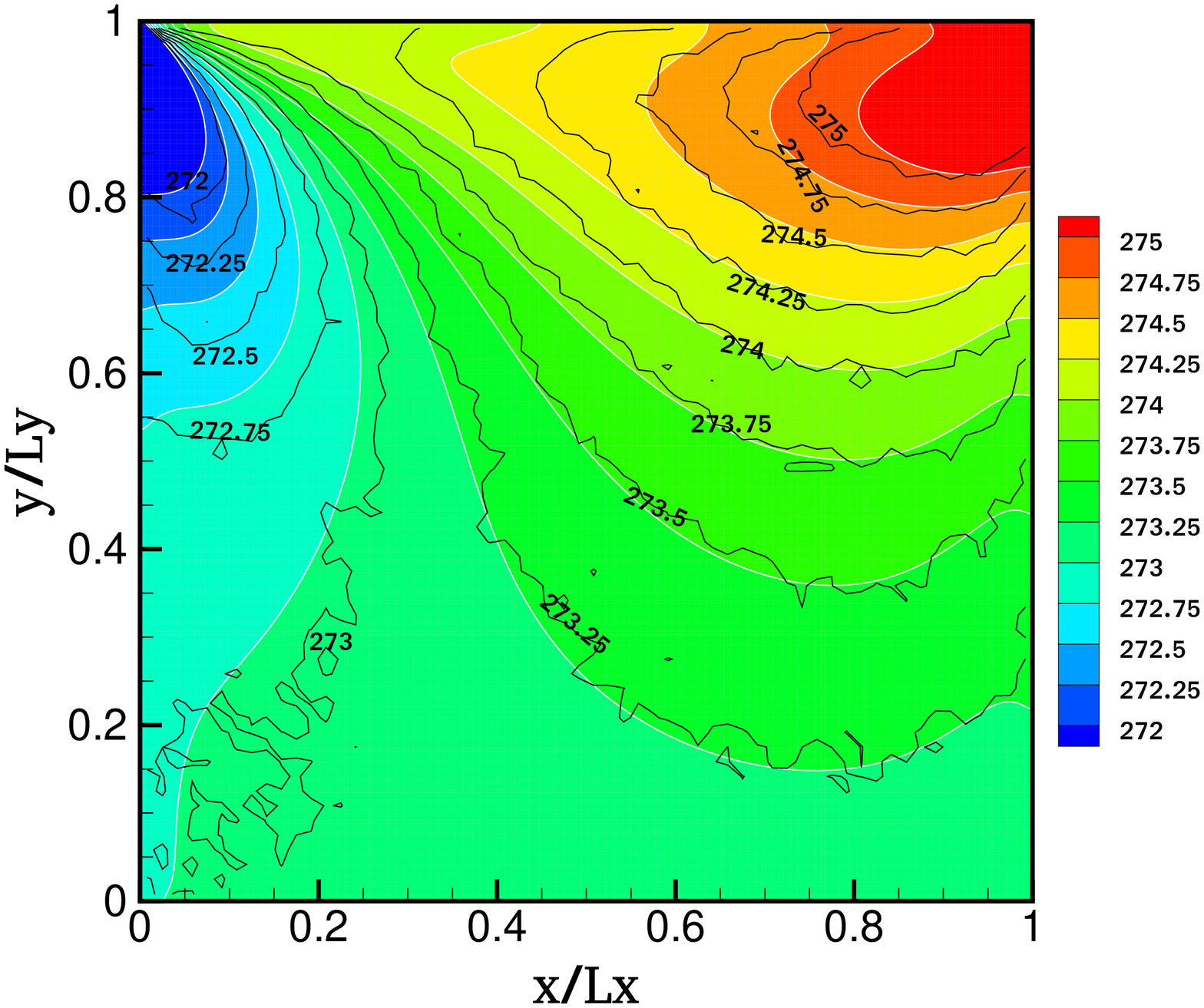}
  }\hfill
  \subfloat[Heat flux streamlines]{
    \label{fig:kn_pt1_shakhov_flux}
    \includegraphics[width=0.47\textwidth]{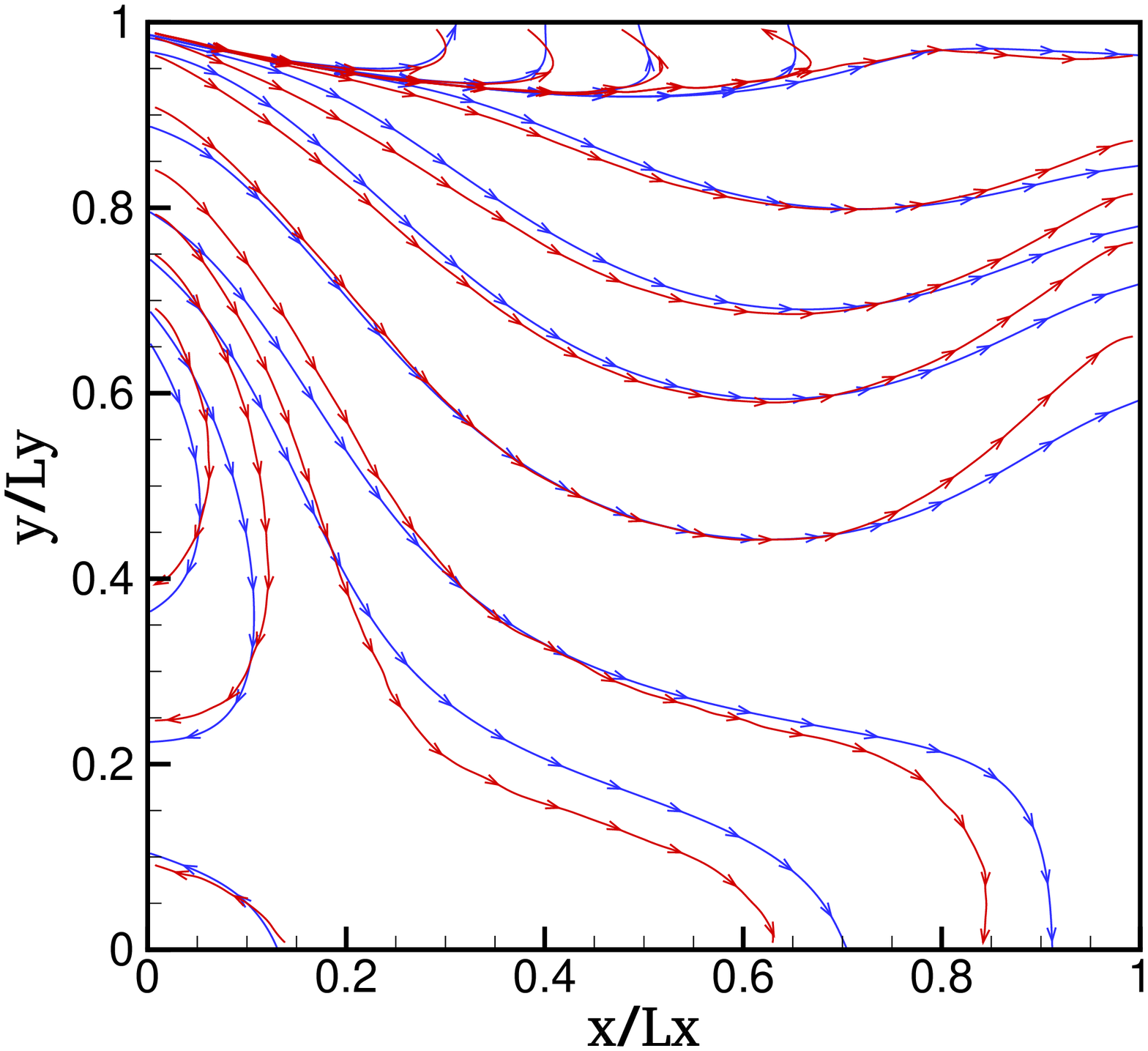}
  }
  \caption{Numerical solutions of the single lid-driven flow for the Shakhov model with  $\Kn=0.1$, $M=20$, and $N_{x}\times N_{y} = 256\times 256$. The black (left) and red (right) lines are the reference DSMC solutions.}
  \label{fig:kn_pt1_shakhov}
\end{figure}

\begin{figure}[!htbp]
  \centering
  \subfloat[Temperature $(\si{\kelvin})$]{
    \label{fig:kn_1_shakhov_T}
    \includegraphics[width=0.47\textwidth]{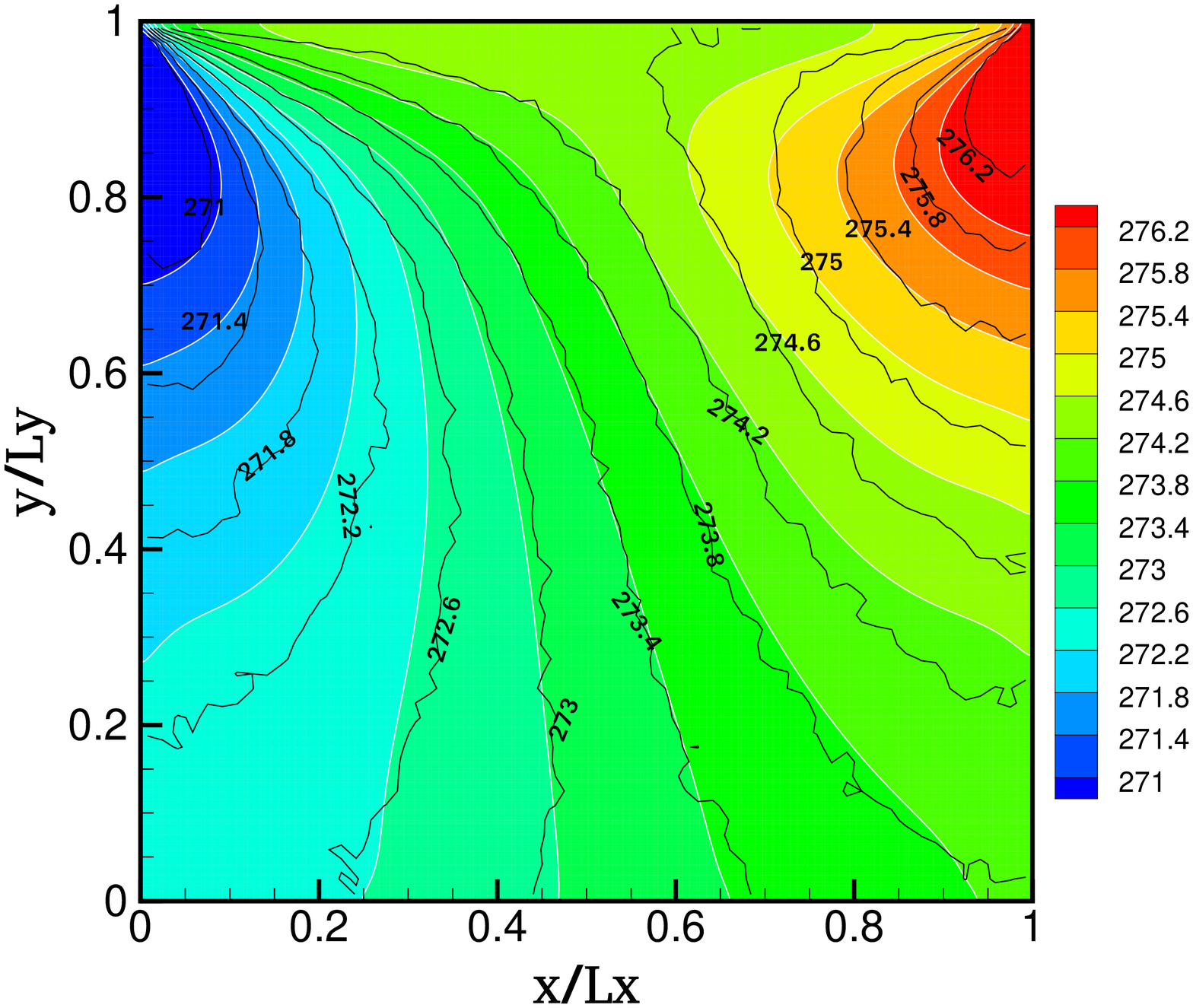}
  }\hfill
  \subfloat[Heat flux streamlines]{
    \label{fig:kn_1_shakhov_flux}
    \includegraphics[width=0.47\textwidth]{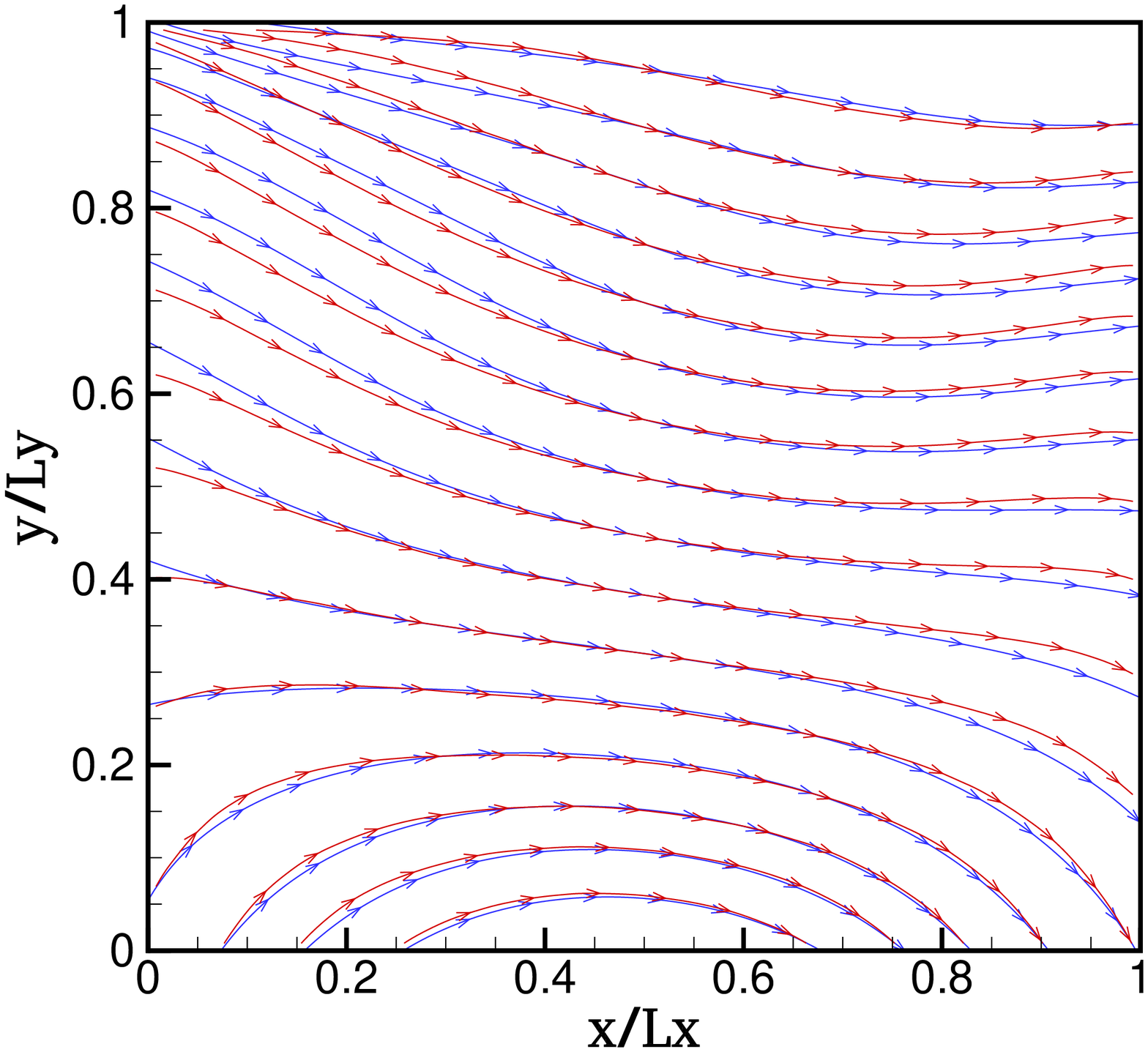}
  }
  \caption{Numerical solutions of the single lid-driven flow for the Shakhov model with  $\Kn=1.0$, $M=25$, and $N_{x}\times N_{y} = 256\times 256$. The black (left) and red (right) lines are the reference DSMC solutions.}
  \label{fig:kn_1_shakhov}
\end{figure}

The Shakhov model is first utilized as the collision model to validate
the solution of the NMG solver. Numerical solutions of the temperature
and heat flux on the uniform grid with $N_{x}=N_{y}=256$ are
presented, respectively, in \cref{fig:kn_pt1_shakhov} for $\Kn=0.1$,
$M=20$, and in \cref{fig:kn_1_shakhov} for $\Kn=1.0$, $M=25$. Therein
the DSMC solutions obtained in \cite{john2010investigation}
are provided as a reference. It can be seen that our results
agree well with the DSMC results at both Knudsen numbers, as revealed
in \cite{Qiao}, where the solution of the moment system
\eqref{eq:moment-system} has been studied in detail. Using the NMG
solver, we are in fact recovering the solution obtained in
\cite{Qiao}. Hence we omit more discussion about the behavior of the
solution here.

\begin{figure}[!htbp]
  \centering
  \subfloat[Temperature $(\si{\kelvin})$]{
    \label{fig:kn_pt1_ESBGK_T}
    \includegraphics[width=0.47\textwidth]{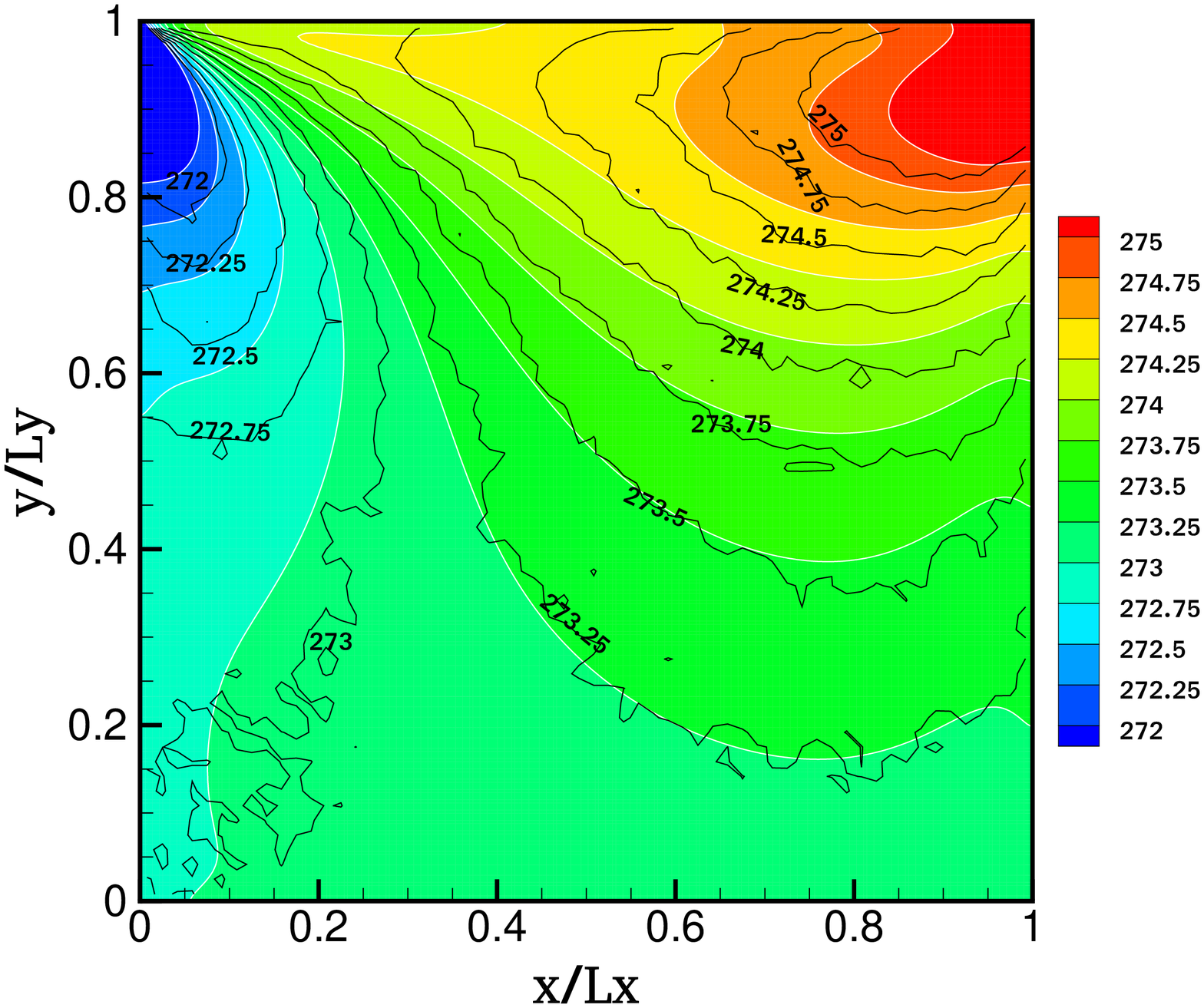}
  }\hfill
  \subfloat[Heat flux streamlines]{
    \label{fig:kn_pt1_ESBGK_flux}
    \includegraphics[width=0.47\textwidth]{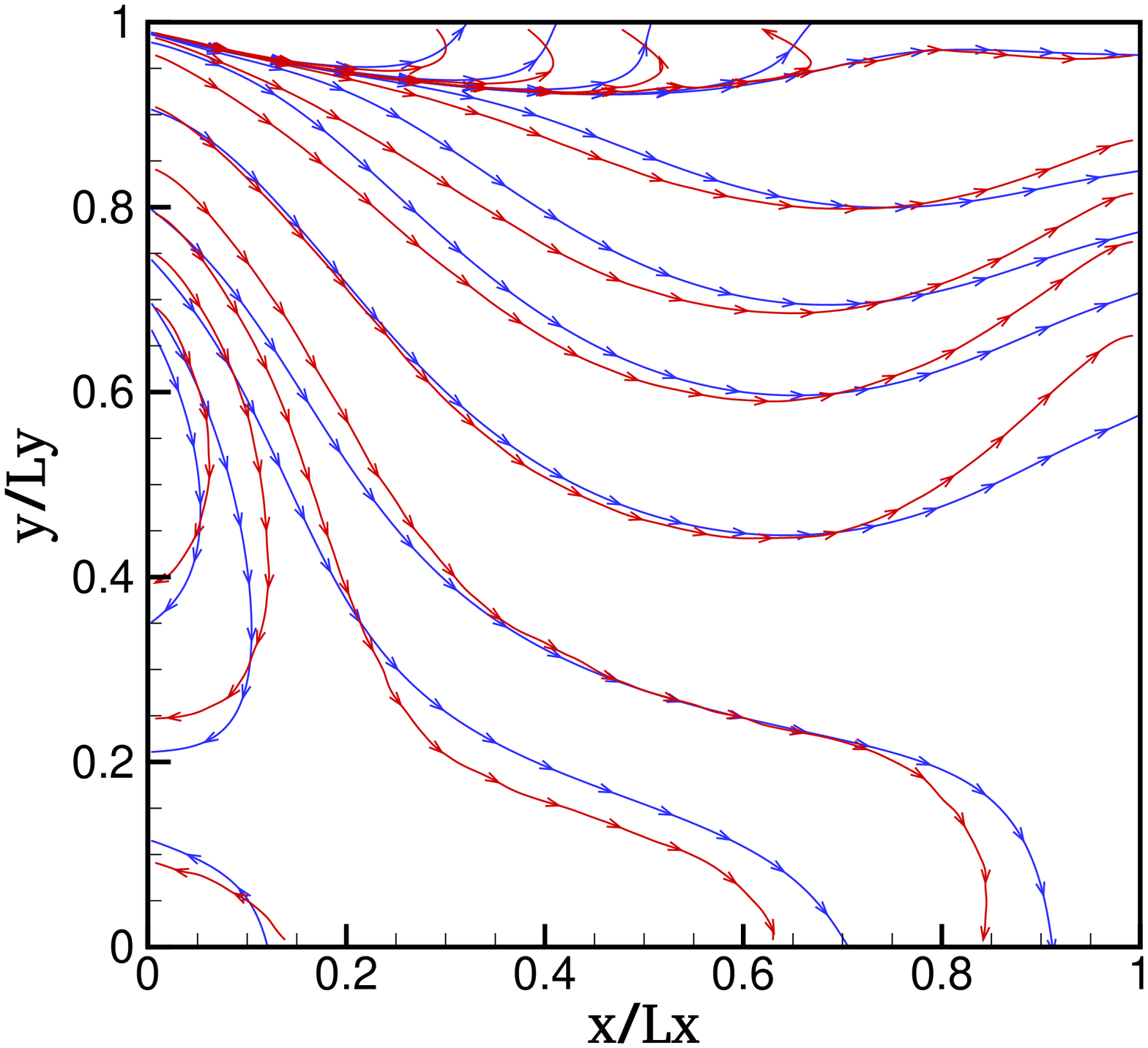}
  }
  \caption{Numerical solutions of the single lid-driven flow for the ES-BGK model with  $\Kn=0.1$, $M=20$, and $N_{x}\times N_{y} = 256\times 256$. The black (left) and red (right) lines are the reference DSMC solutions.}
  \label{fig:kn_pt1_ESBGK}
\end{figure}

\begin{figure}[!htbp]
  \centering
  \subfloat[Temperature $(\si{\kelvin})$]{
    \label{fig:kn_1_ESBGK_T}
    \includegraphics[width=0.47\textwidth]{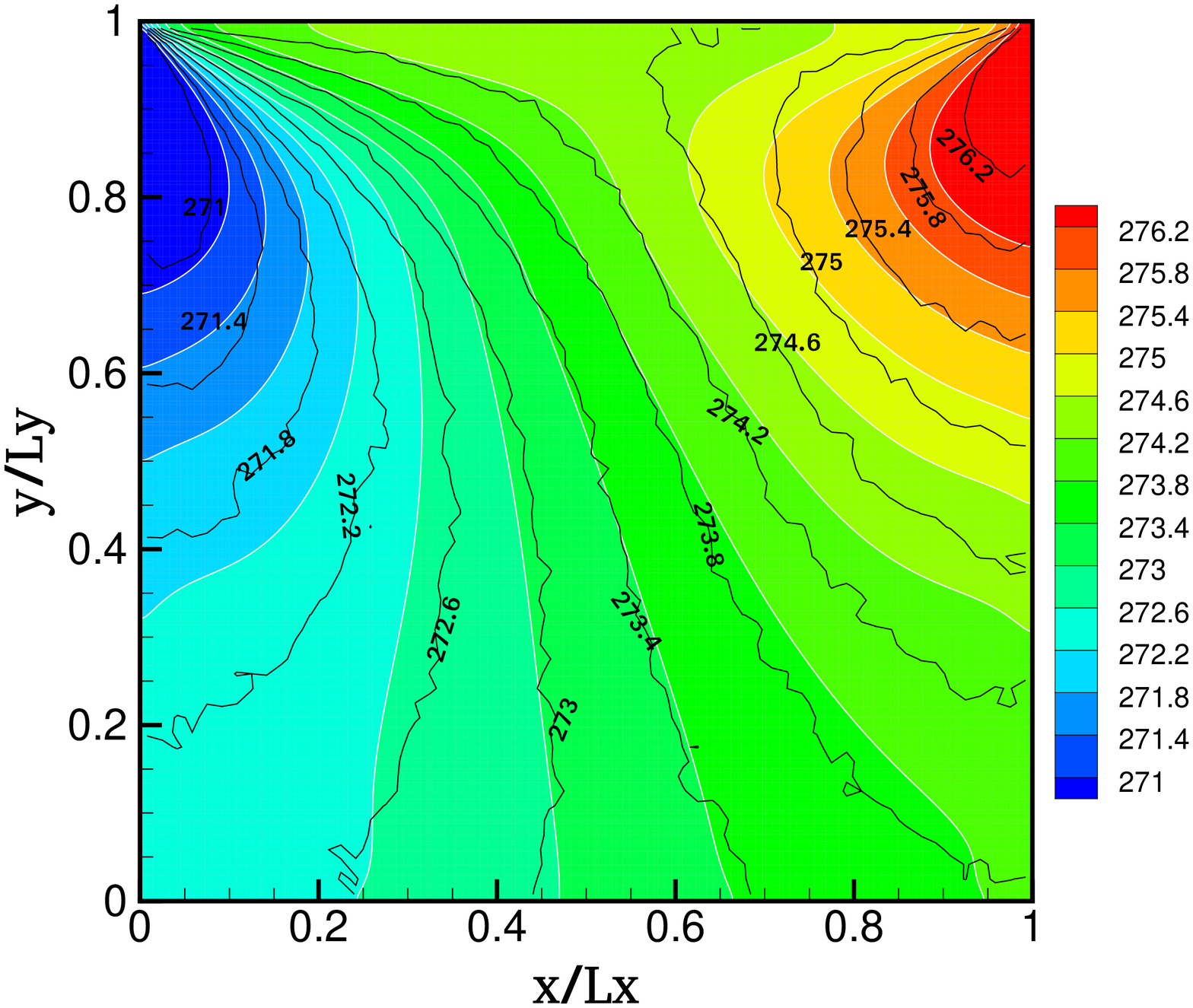}
  }\hfill
  \subfloat[Heat flux streamlines]{
    \label{fig:kn_1_ESBGK_flux}
    \includegraphics[width=0.47\textwidth]{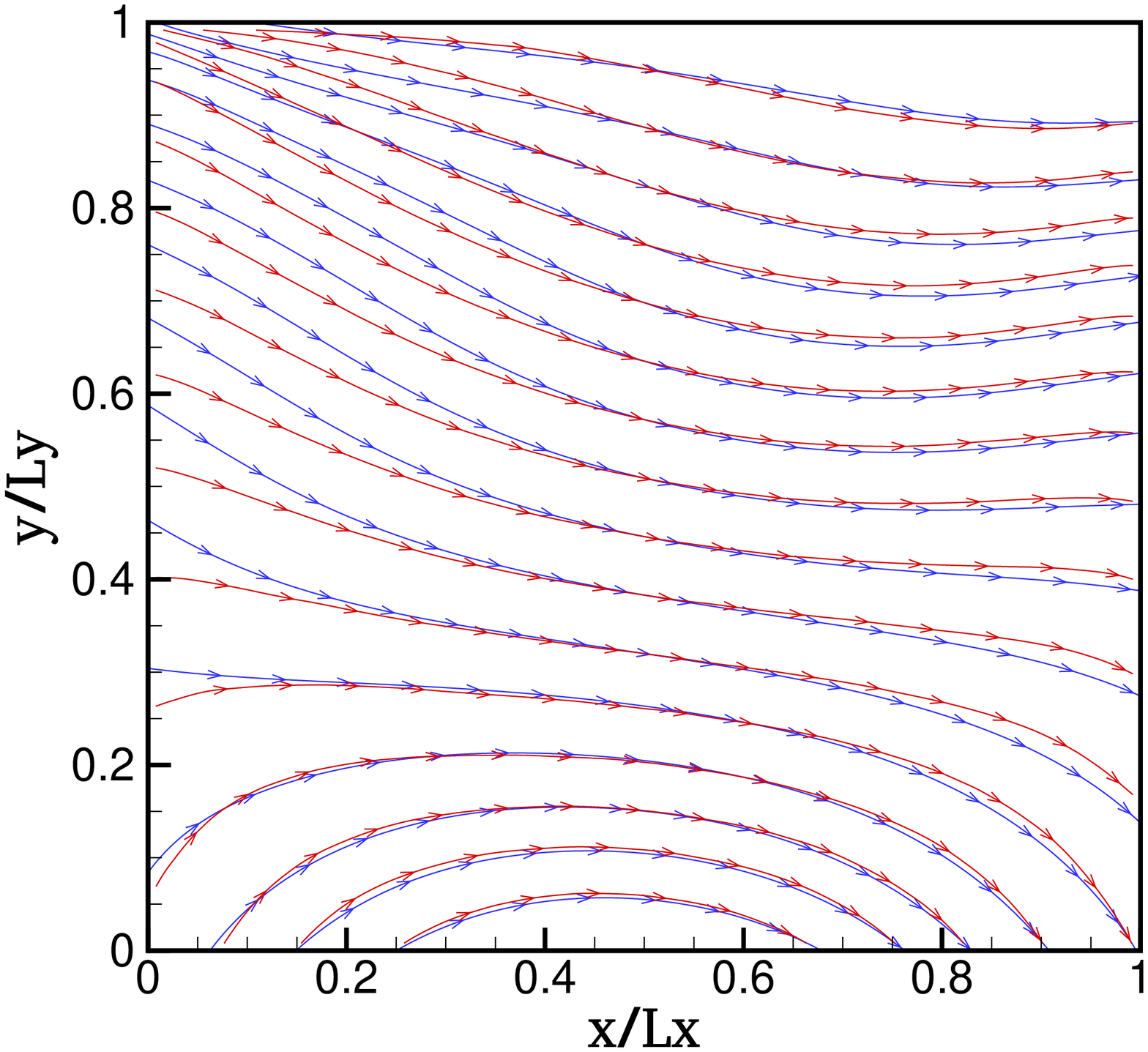}
  }
  \caption{Numerical solutions of the single lid-driven flow for the ES-BGK model with  $\Kn=1.0$, $M=25$, and $N_{x}\times N_{y} = 256\times 256$. The black (left) and red (right) lines are the reference DSMC solutions.}
  \label{fig:kn_1_ESBGK}
\end{figure}

For the ES-BGK model, similar solutions can be obtained by the NMG
solver. As an example, numerical solutions of the temperature and
heat flux for the ES-BGK model, as well as the reference DSMC
solutions, are shown in \cref{fig:kn_pt1_ESBGK} for $\Kn=0.1$ and
\cref{fig:kn_1_ESBGK} for $\Kn=1.0$, respectively. It can be observed
that for the current problem, the numerical solutions of the Shakhov
model agree slightly better with the reference DSMC solutions than the
numerical solutions of the ES-BGK model.

\subsubsection{Numerical efficiency}
\label{sec:num-ex-single-lid-efficiency}

To explore its efficiency and behavior, the NMG solver is performed
with a variety of $M$ on a sequence of uniform grids for both the
Shakhov and ES-BGK models. Since similar features of the NMG solver
are observed for all cases, only partial results are reported in the
present paper.

\begin{figure}[!htbp]
  \centering
    \includegraphics[width=0.49\textwidth]{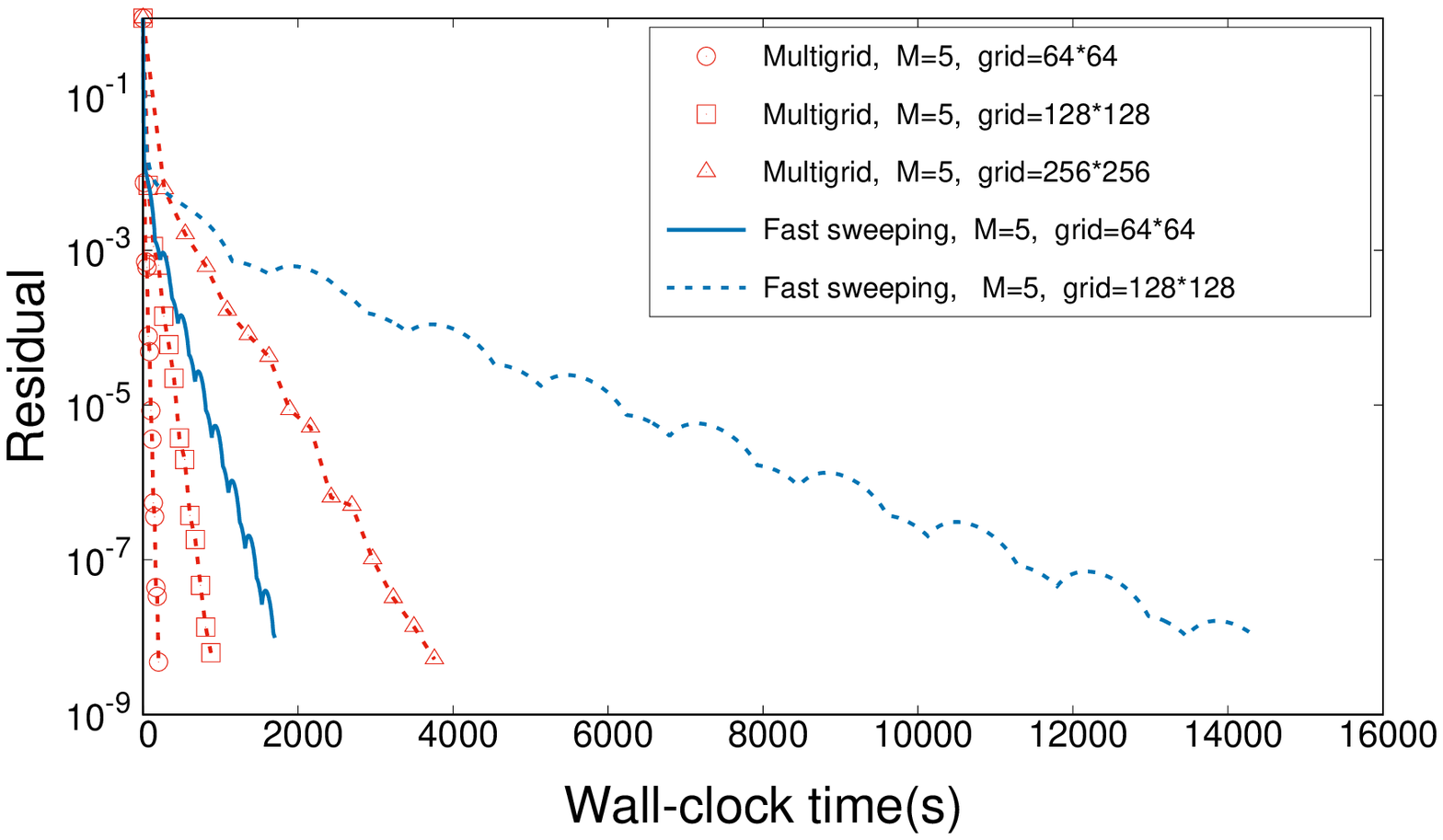}
   \hfill
    \includegraphics[width=0.49\textwidth]{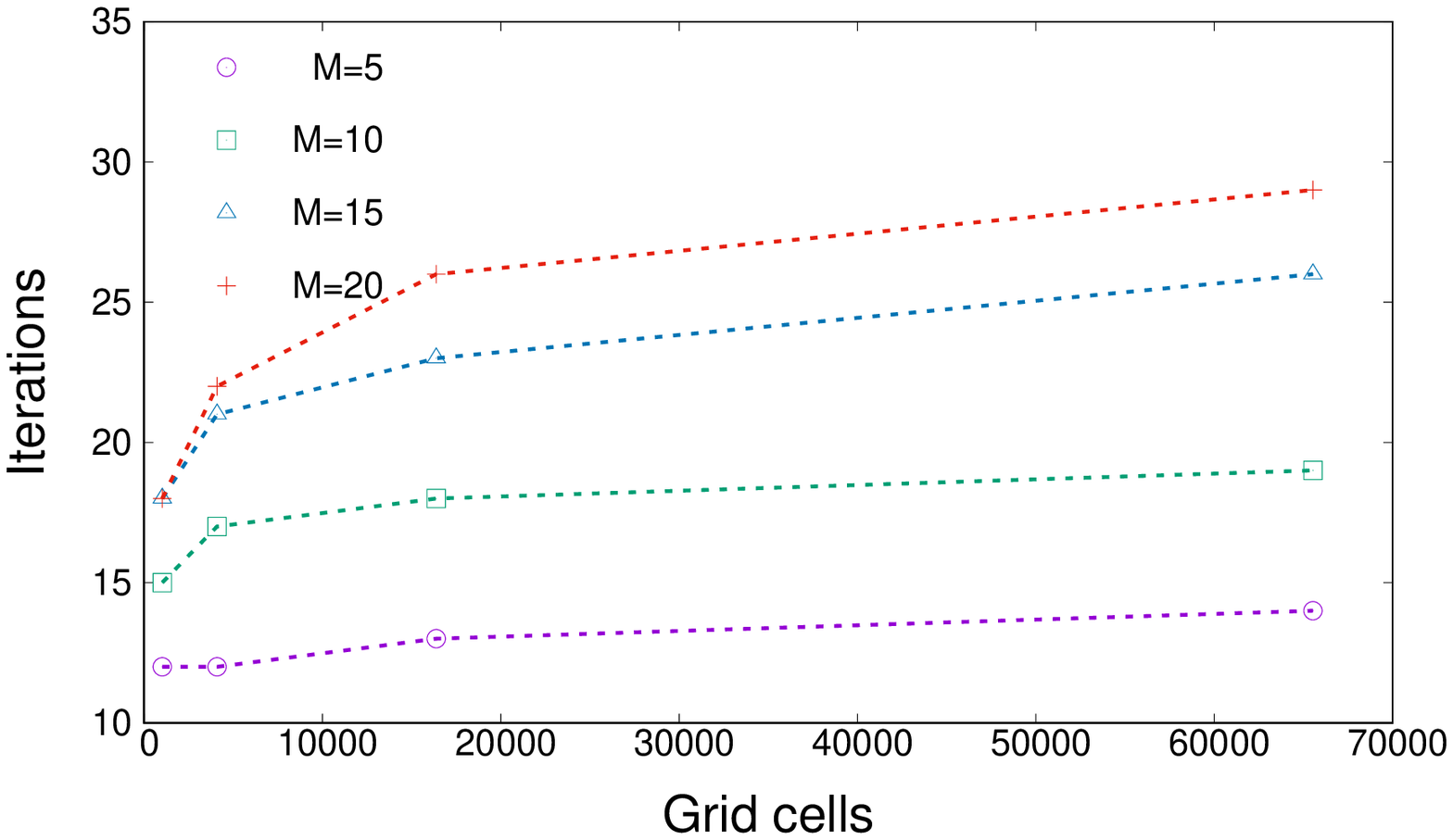}
    \caption{Convergence history (left) and 
      iterations in terms of $N_{x}\times N_{y}$ of the NMG solver
      (right) for the Shakhov model of single lid-driven flow with
      $\Kn=0.1$.}
  \label{fig:kn_pt1_shakhov_cpu_ites}
\end{figure}

In the case of the Shakhov model, the total number of iterations and
the wall-clock time, spent by the NMG solver with $M=5$ on four
different grids, are listed in \cref{tab:Kn_pt1_Shakhov_grid} for
$\Kn=0.1$. For comparison, the corresponding results of two single
level solvers, that is, the forward Euler scheme and the fast sweeping
iteration, are also presented in
\cref{tab:Kn_pt1_Shakhov_grid}. Apparently, the fast sweeping
iteration converges faster than the forward Euler scheme as
expected. It only takes about $1/8$ of iterations of the forward Euler
scheme. Benefiting from this, the wall-clock time of the fast sweeping
iteration is saved to ${1}/{3}$ of that of the forward Euler
scheme. Moreover, for both single level solvers, it can be also seen
that roughly the total number of iterations is doubled and the
wall-clock time is octupled, as the grid is refined each time, for
which $N_{x}$ and $N_{y}$ are doubled. In contrast, the NMG solver
behaves much better than both single level solvers by observing that
it not only converges in a dozen or so iterations, but the total
number of iterations is almost independent of grid size. As a result,
the NMG solver becomes more and more efficient as the grid is
refined. Specifically, the wall-clock time ratio of the NMG solver to
the fast sweeping iteration is about $25.2\,\%$ on the grid composed
of $32\times 32$ cells, whereas it is reduced drastically to $3.5\,\%$
on the grid composed of $256\times 256$ cells.
In addition, some convergence histories of the NMG solver and the fast
sweeping iteration are plotted in \cref{fig:kn_pt1_shakhov_cpu_ites}
(left), which also shows the efficiency of the NMG solver.

\begin{table}[!htbp]
  \centering
  \caption{Performance results for the Shakhov model of single
    lid-driven flow with $\Kn=0.1$ and $M=5$. Euler: the forward Euler
    scheme; FS: the fast sweeping iteration; NMG: the NMG solver;
    $T_{r}$: the wall-clock time ratio of NMG to FS.}
  \label{tab:Kn_pt1_Shakhov_grid}
  \begin{tabular}{cccccccc}
    \toprule
    \multirow{2}*{$N_{x}\times N_{y}$} & \multicolumn{3}{c}{Iterations} & \multicolumn{4}{c}{Wall-clock time $(\si{\s})$} \\
    \cmidrule(r){2-4} \cmidrule(r){5-8}
    & Euler & FS & NMG & Euler & FS & NMG & $T_{r}$ \\ 
    \midrule
    $32 \times 32$  & 1823 & 231 & 12 & 619.84 & 199.63 & 50.40 & 25.2\,\% \\
    $64 \times 64$   & 4100 & 507 & 12 & 5433.38 & 1701.1 & 199.28 & 11.7\,\% \\
    $128 \times 128$   & 8398 & 1084 & 13 & 44270.3 & 14373.2 & 875.00 & 6.1\,\% \\
    $256 \times 256$   & 16364 & 2000 & 14 & 356118 & 105979 & 3760.26 & 3.5\,\% \\
    \bottomrule
  \end{tabular}
\end{table}

\begin{figure}[!htbp]
  \centering
  \includegraphics[width=0.49\textwidth]{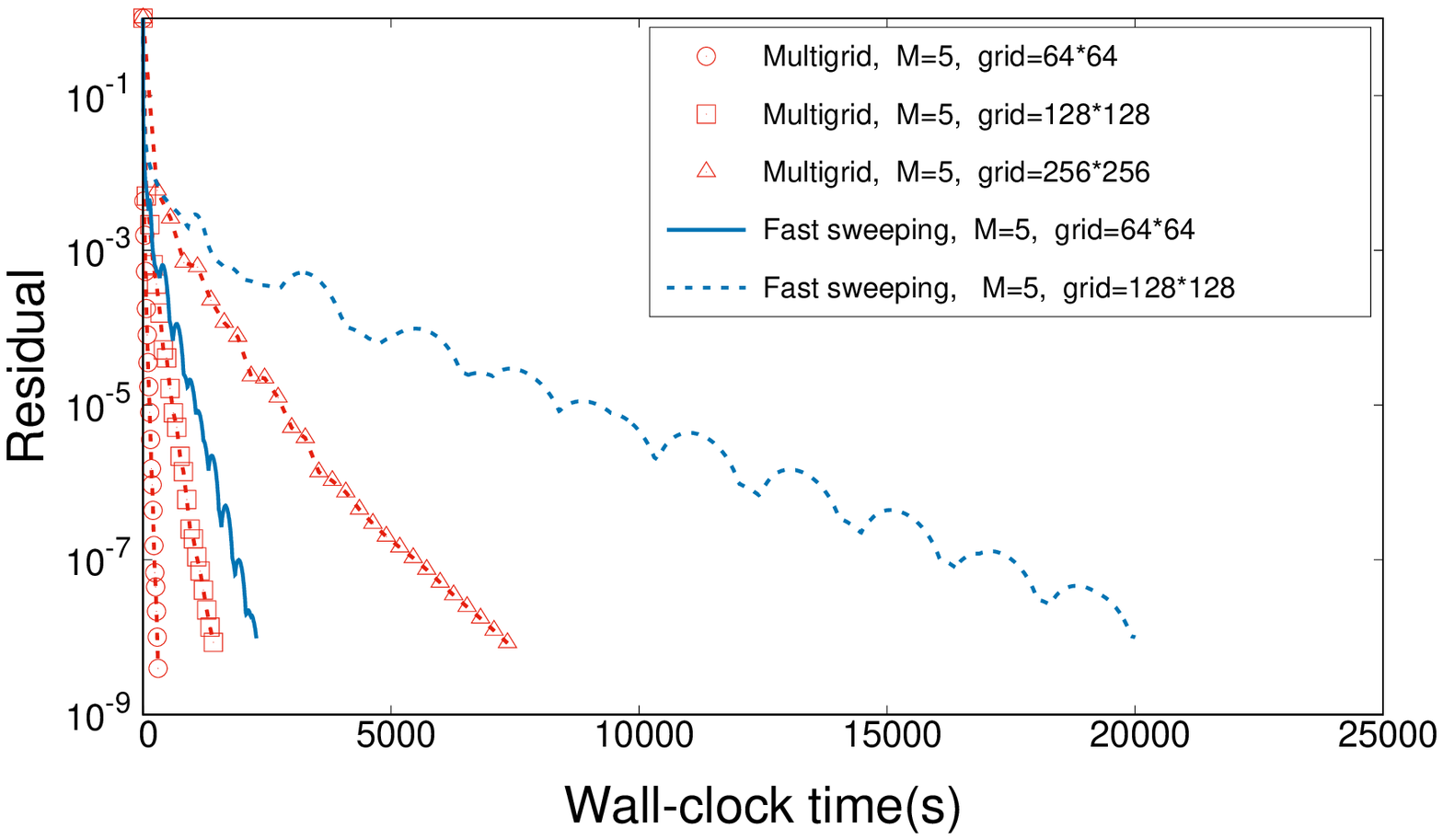}
   \hfill
   \includegraphics[width=0.49\textwidth]{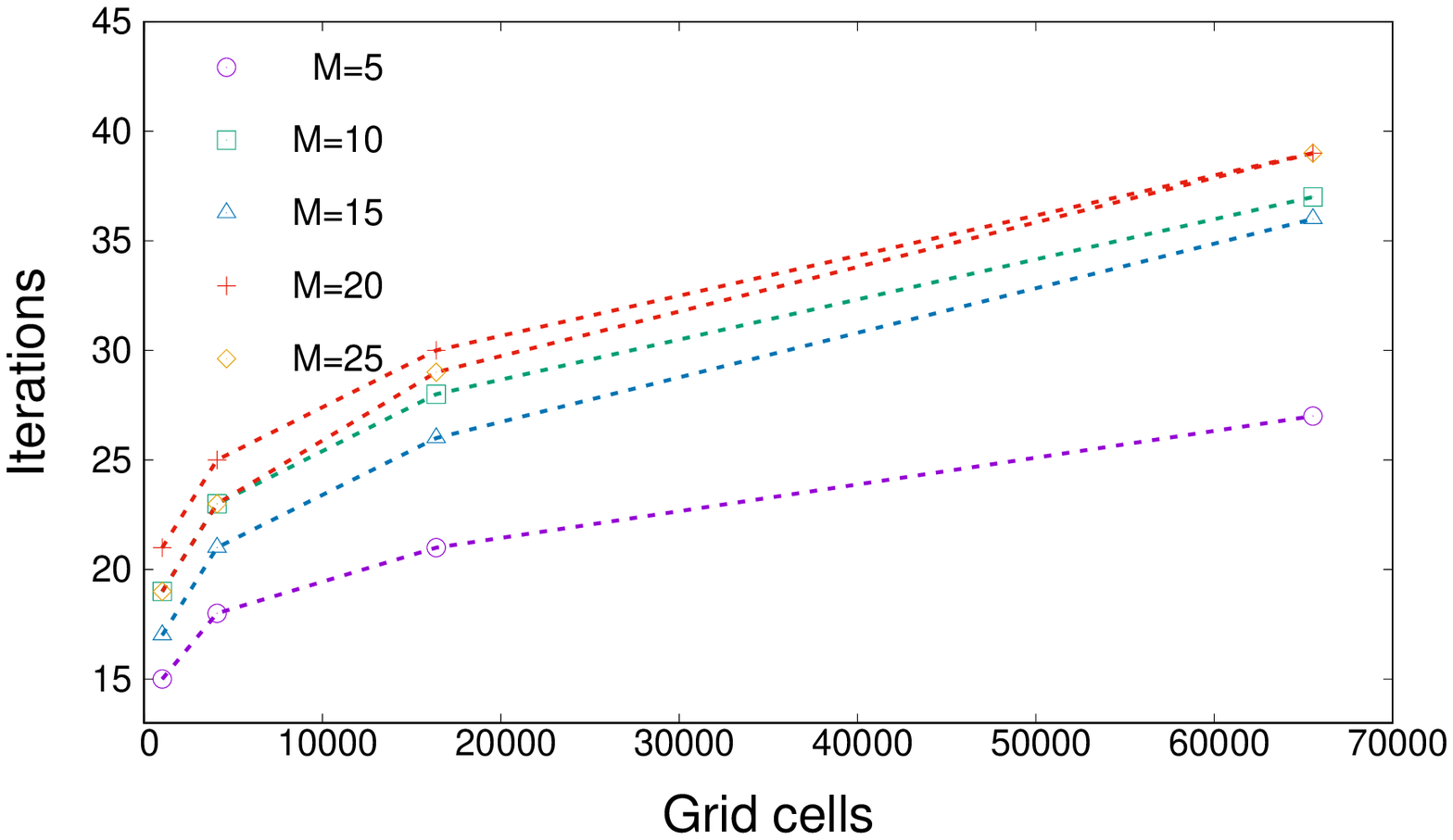}
    \caption{Convergence history (left) and 
      iterations in terms of $N_{x}\times N_{y}$ of the NMG solver
      (right) for the Shakhov model of single lid-driven flow with
      $\Kn=1.0$.}
  \label{fig:kn_1_shakhov_cpu_ites}
\end{figure}

As the Knudsen number increases to $\Kn=1.0$, the performance results
of the NMG solver, as well as two single level solvers, are listed in
\cref{tab:Kn_1_Shakhov_grid} for $M=5$. It is shown that the forward
Euler scheme and the fast sweeping iteration converge more slowly and
their total number of iterations increase slightly faster than the
case of $\Kn=0.1$. Yet the ratios of iterations and wall-clock time
between two single level solvers seem to be preserved well. Since the
NMG solver adopts the fast sweeping iteration as the smoother and the
coarsest grid solver, it converges also more slowly than the case of
$\Kn=0.1$. Besides, the total number of iterations can not be
maintained any more, and a few increments of it are found as the grid
is refined. However, this phenomenon is acceptable in view of that the
parameters such as $s_{\iota}$, $\iota=1,2,3$, are fixed in the NMG
solver, and the convergence histories of the fast sweeping iteration,
shown in \cref{fig:kn_1_shakhov_cpu_ites} (left) for $\Kn=1.0$,
fluctuate a little more in comparison to the results shown in
\cref{fig:kn_pt1_shakhov_cpu_ites} (left) for $\Kn=0.1$. What is more,
the wall-clock time ratios of the NMG solver to the fast sweeping
iteration, ranged from $26.5\,\%$ to $4.2\,\%$ for the grids ranged
from $32\times 32$ to $256\times 256$ cells, are preserved
satisfactorily as the case of $\Kn=0.1$. Therefore, the NMG solver
could be still more and more efficient as the grid is refined.

For the solvers with other $M$, the performance results for the
Shakhov model on the grid composed of $64\times 64$ cells are given in
\cref{tab:single-shakhov-M}. It can be seen, at both Knudsen numbers,
that the total number of iterations for all
solvers increase gradually as $M$ grows, except for the case of $M=15$ and
$\Kn=1.0$, for which they converge faster than the case of even
$M$. Actually, as mentioned in \cite{hu2019efficient, Cai2018},
different behaviors with respect to the parity of $M$ have been
observed for two single level solvers, and so is the NMG solver,
especially in the case of $\Kn=1.0$.
Nevertheless, the acceleration of steady-state computation by using
the fast sweeping iteration or the NMG solver is still prominent for
all cases. 
Concretely, although the ratios of iterations and wall-clock time of
the fast sweeping iteration to the forward Euler scheme increase for
large $M$ in comparison to the case of $M=5$, it is found that the
maximum ratio of wall-clock time is just about $40\,\%$. As for the
NMG solver, it always converges in dozens of iterations, and the
wall-clock time ratios of it to the fast sweeping iteration are around
$12\,\%$ and $13\,\%$ corresponding to $\Kn=0.1$ and $1.0$,
respectively, for all $M$. Furthermore, the variations of total number
of iterations for the NMG solver in terms of the grid size
$N_{x}\times N_{y}$ with several choices of $M$ are plotted separately
in \cref{fig:kn_pt1_shakhov_cpu_ites} (right) for $\Kn=0.1$ and
\cref{fig:kn_1_shakhov_cpu_ites} (right) for $\Kn=1.0$. It is observed
for each $M$ that the total number of iterations grows faster in the
case of $\Kn=1.0$ than in the case of $\Kn=0.1$. However, since the total
number of iterations would be at least doubled as the grid is refined
for two single level solvers, it turns out that the efficiency of the
NMG solver for steady-state computation would be more and more
evident as both $N_{x}$ and $N_{y}$ increase.

\begin{table}[!htbp]
  \centering
  \caption{Performance results for the Shakhov model of single
    lid-driven flow with $\Kn=1.0$ and $M=5$. Euler: the forward Euler
    scheme; FS: the fast sweeping iteration; NMG: the NMG solver;
    $T_{r}$: the wall-clock time ratio of NMG to FS.}
  \label{tab:Kn_1_Shakhov_grid}
  \begin{tabular}{cccccccc}
    \toprule
    \multirow{2}*{$N_{x}\times N_{y}$} & \multicolumn{3}{c}{Iterations} & \multicolumn{4}{c}{Wall-clock time $(\si{\s})$} \\
    \cmidrule(r){2-4} \cmidrule(r){5-8}
    & Euler & FS & NMG & Euler & FS & NMG & $T_{r}$\\
    \midrule
    $32 \times 32$  & 1978 & 277 & 15 & 673.30 & 239.31 & 63.47 & 26.5\,\% \\
    $64 \times 64$   & 4952 & 681 & 18 & 6636.03 & 2284.32 & 303.99 & 13.3\,\% \\
    $128 \times 128$   & 11120 & 1509 & 21 & 58647.5 & 20003.1 & 1417.91 & 7.1\,\% \\
    $256 \times 256$   & 24266 & 3300 & 27 & 525312 & 176582 & 7348 & 4.2\,\% \\
    \bottomrule
  \end{tabular}
\end{table}

\begin{table}[!htbp]
  \centering
  \caption{Performance results for the Shakhov model of single
    lid-driven flow with $N_{x}= N_{y} = 64$. Euler: the
    forward Euler scheme; FS: the fast sweeping iteration; NMG: the
    NMG solver; $T_{r}$: the wall-clock time ratio of NMG to FS.}
  \label{tab:single-shakhov-M}
  \begin{tabular}{ccccccccc}
    \toprule
    \multirow{2}*{$\Kn$} & \multirow{2}*{$M$} & \multicolumn{3}{c}{Iterations} & \multicolumn{4}{c}{Wall-clock time $(\si{\s})$}  \\
    \cmidrule(r){3-5} \cmidrule(r){6-9}
&    & Euler & FS & NMG & Euler & FS & NMG & $T_{r}$\\
    \midrule
    0.1 & 10  & 5244 & 708 & 17 & 33569.5 & 11675.8 & 1401.07	 & 12.0\,\% \\
    & 15  & 6315 & 875 & 21 & 116762 & 43611.8 & 5008.71 & 11.5\,\% \\
    & 20  & 6839 & 936 & 22 & 295431 & 114102 & 12596 & 11.0\,\% \\[0.3em]
    1.0 & 10  & 6051 & 866 & 23 & 38707.2 & 14249.8 & 1892.43 & 13.3\,\% \\
        & 15  & 4805 & 702 & 21 & 88035.9 & 34995.3 & 5027.71 & 14.4\,\% \\
        & 20  & 6987 & 1008 & 25 & 301326 & 122825 & 14533.7 & 11.8\,\% \\
    \bottomrule
  \end{tabular}
\end{table}

\begin{figure}[!htbp]
  \centering
  \includegraphics[width=0.49\textwidth]{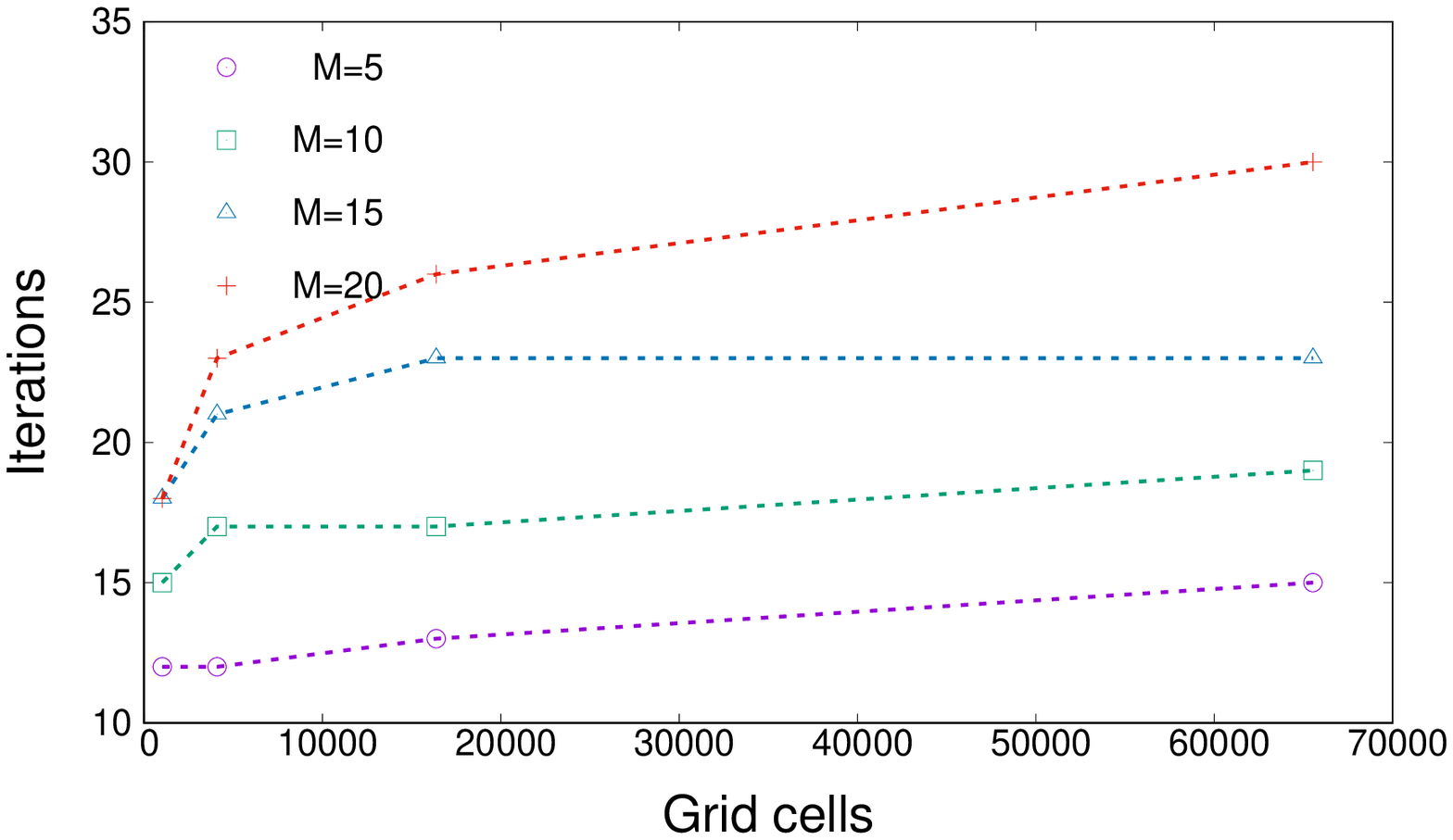}
   \hfill
   \includegraphics[width=0.49\textwidth]{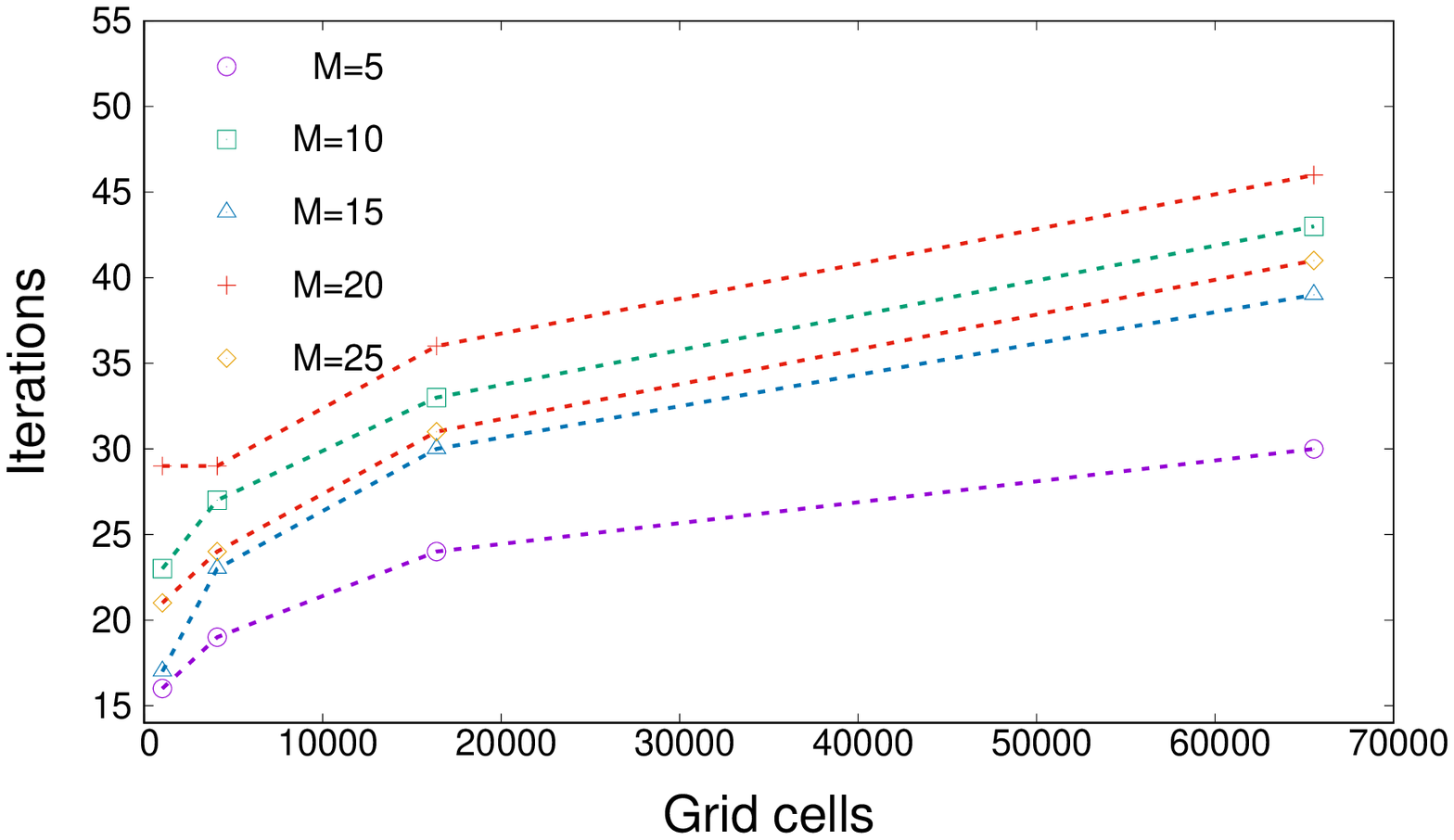}
   \caption{Iterations in terms of $N_{x}\times N_{y}$ of the NMG
     solver for the ES-BGK model of single lid-driven flow with
     $\Kn=0.1$ (left) and $\Kn=1.0$ (right).}
  \label{fig:kn_ESBGK_ites_grid}
\end{figure}

The three solvers, discussed above, also work fairly well for the
ES-BGK model. As examples, the performance results on the grid composed of
$64\times 64$ cells are exhibited in \cref{tab:single-ESBGK-M}, and
the variations of total number of iterations for the NMG solver in
terms of $N_{x}\times N_{y}$ with several choices of $M$ at both
$\Kn=0.1$ and $1.0$ are plotted in \cref{fig:kn_ESBGK_ites_grid}. It
follows that in the case of $\Kn=0.1$, all solvers for the ES-BGK
model behave almost the same as for the Shakhov model, while in the
case of $\Kn=1.0$, the forward Euler scheme requires much more
iterations to guarantee convergence as $M$ increases, compared with
the same case for the Shakhov model. Consequently, in this case, the
fast sweeping iteration and the NMG solver converge more slowly than
the counterpart solver for the Shakhov model. In spite of this, the
speedup ratio of the fast sweeping iteration to the forward Euler
scheme is found to be maintained, and so is the ratio of the NMG
solver to the fast sweeping iteration. It is also observed from
\cref{fig:kn_ESBGK_ites_grid} that the NMG solver converges in dozens
of iterations even when $\Kn=1.0$ and the grid consists of
$256\times 256$. Therefore, a great improvement in efficiency can be
again obtained by using the NMG solver.

At last, four results of the parallel NMG solver on the grid with
$256\times 256$ cells for both Shakhov and ES-BGK models are listed in
\cref{tab:single-lid-parallelNMG}, which shows a nice speedup is able
to be obtained for the NMG solver by using OpenMP-based
parallelization. Specifically, the speedup ratio is around $2$ and
$3.7$ when, respectively, $2$ and $4$ threads are applied in parallel
computation. Besides, it can be seen that the average of total number
of iterations for the parallel NMG solver is almost the same as the
serial NMG solver, although the details of iterations are not exactly
equivalent as mentioned in \cref{sec:solver-openmp}.

\begin{table}[!htbp]
  \centering
  \caption{Performance results for the ES-BGK model of single
    lid-driven flow with $N_{x}= N_{y} = 64$. Euler: the
    forward Euler scheme; FS: the fast sweeping iteration; NMG: the
    NMG solver; $T_{r}$: the wall-clock time ratio of NMG to FS.}
  \label{tab:single-ESBGK-M}
  \begin{tabular}{ccccccccc}
    \toprule
    \multirow{2}*{$\Kn$} & \multirow{2}*{$M$} & \multicolumn{3}{c}{Iterations} & \multicolumn{4}{c}{Wall-clock time $(\si{\s})$}  \\
    \cmidrule(r){3-5} \cmidrule(r){6-9}
    & & Euler & FS & NMG & Euler & FS & NMG & $T_{r}$\\
    \midrule
    0.1 & 5 & 4085 & 503 & 12 & 5138.97 & 1636.6 & 198.05 & 12.1\,\% \\
    & 10  & 5177 & 700 & 17 & 31352 & 10963.7 & 1389.26	& 12.7\,\% \\
    & 15  & 6357 & 867 & 21 & 110388 & 40007.2 & 4915.99 & 12.3\,\% \\
    & 20  & 6764 & 926  & 23 & 270399 & 107218 & 13268.6 & 12.4\,\% \\[0.3em]
    1.0 & 5 & 5053 & 690 & 19 & 6347.05 & 2253.42 & 302.39 & 13.4\,\% \\
        & 10  & 7663 & 1105 & 27 & 46164.6 & 17154.5 & 2086.02 & 12.2\,\% \\
    & 15  & 5541 & 780 & 23 & 96388.4 & 36128.8 & 5154.89	& 14.3\,\% \\
    & 20  & 9177 & 1339 & 29 & 364166 & 158320 & 16619	& 10.5\,\% \\    
    \bottomrule
  \end{tabular}
\end{table}

\begin{table}[!htbp]
  \centering
  \caption{Performance results of the parallel NMG solver with
    $N_{x}= N_{y} = 256$ for the single lid-driven flow. $T_{r,n}$:
    the wall-clock time ratio of single thread to $n$ threads.}
  \label{tab:single-lid-parallelNMG}
  \begin{tabular}{ccccccccc}
    \toprule
    \multirow{2}*{\shortstack{Collision\\ model}} & \multirow{2}*{$\Kn$} & \multirow{2}*{$M$} & \multirow{2}*{\shortstack{Average\\ iterations}}   & \multicolumn{5}{c}{Wall-clock time $(\si{\s})$ for $n$ threads} \\
    \cmidrule(r){5-9}
    &  &  &  & $n=1$ & $n=2$ & $T_{r,2}$ & $n=4$ & $T_{r,4}$\\
    \midrule
    Shakhov & $0.1$ & $20$ & 30 & 215773 & 112842 & 1.91 & 59086.3 & 3.65 \\
    ES-BGK  & $0.1$ & $20$ & 30 & 209509 & 107997 & 1.94 & 56486.8 & 3.71 \\
    Shakhov & $1.0$ & $25$ & 39 & 554987 & 277464 & 2.00 & 147878 & 3.75 \\
    ES-BGK  & $1.0$ & $25$ & 41 & 549046 & 276126 & 1.99 & 146859 & 3.74 \\
    \bottomrule
  \end{tabular}
\end{table}

\subsection{Four-sided lid-driven cavity flow}
\label{sec:num-ex-four-lid}

The four-sided lid-driven cavity flow is also very important for both
benchmarking and application viewpoints, and has been studied in
\cite{nabapure2021dsmc}. The geometry of the anti-parallel wall motion
case considered in \cite{nabapure2021dsmc} is shown in
\cref{fig:four-lid-geo}. That is, the top and right sides move
horizontally to the right and vertically up, respectively, with a
constant speed $U_{W}$, while the other two sides move in the opposite
direction with the same speed relative to the motion of their parallel
sides. We adopt the setting as follows. The constant speed is given by
$U_{W}=\SI[per-mode=symbol]{50}{\m\per\s}$. All the cavity walls are
maintained at temperature $\SI{273}{\kelvin}$. The length and height
of the cavity are $L_{x} = L_{y} = \SI{1}{\m}$. Initially, the gas is
again assumed to be uniformly distributed and in the Maxwellian with
constant density, mean velocity of $0$, and temperature of
$\SI{273}{\kelvin}$.
The initial density taken into account is
$\rho = \SI[per-mode=symbol]{1.1044e-7}{\kg \per \m^{3}}$. Thus the
associated Knudsen number reads $\Kn=0.777$.

\begin{figure}[!htbp]
  \centering
  \subfloat[Velocity streamlines]{
    \label{fig:kn_pt777_shakhov_v}
    \includegraphics[width=0.47\textwidth]{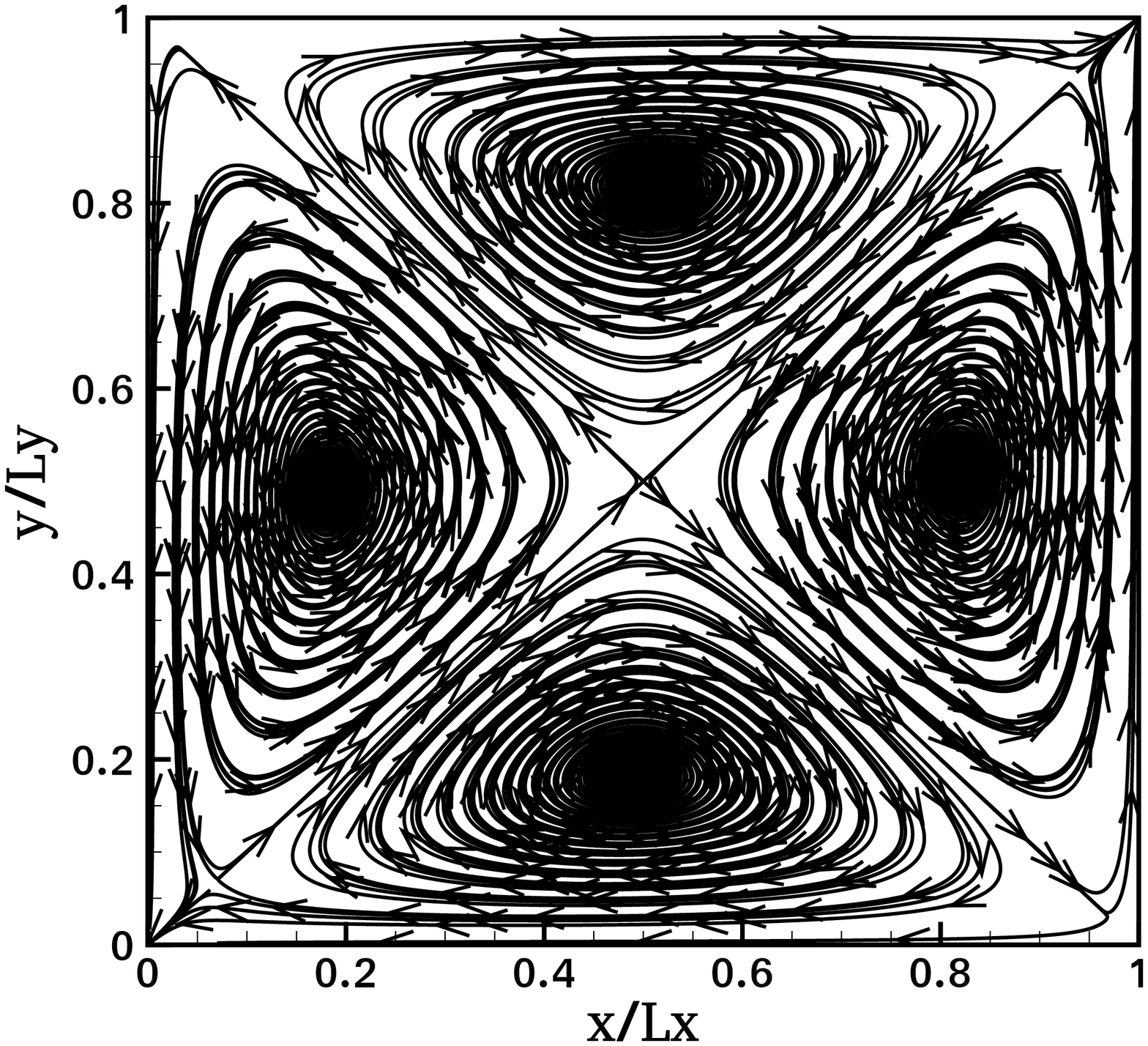}
  }\hfill
  \subfloat[Temperature $(\si{\kelvin})$]{
    \label{fig:kn_pt777_shakhov_T}
    \includegraphics[width=0.47\textwidth]{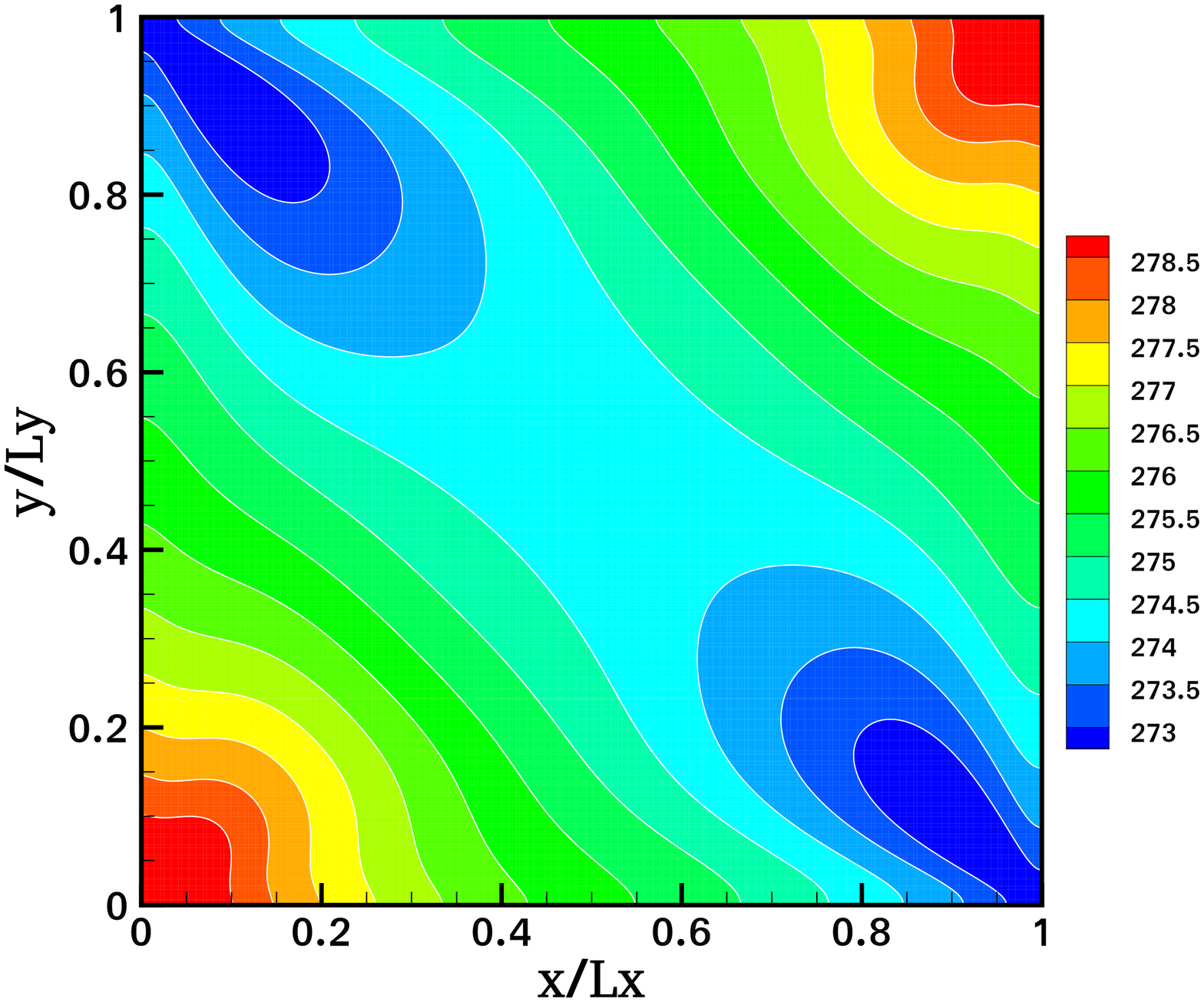}
  }
  \caption{Numerical solutions of the four-sided lid-driven flow for the Shakhov model with $M=25$ and $N_{x}\times N_{y} = 256\times 256$.}
  \label{fig:kn_pt777_shakhov}
\end{figure}

The Shakhov model is employed as the collision model for this
problem. The velocity streamlines and the temperature contours
obtained by the NMG solver with $M=25$ and $N_{x}=N_{y}=256$ are shown
in \cref{fig:kn_pt777_shakhov}. These results exhibit the similar
structure in comparison to the DSMC results presented in
\cite{nabapure2021dsmc}.

\begin{figure}[!htbp]
  \centering
  \includegraphics[width=0.49\textwidth]{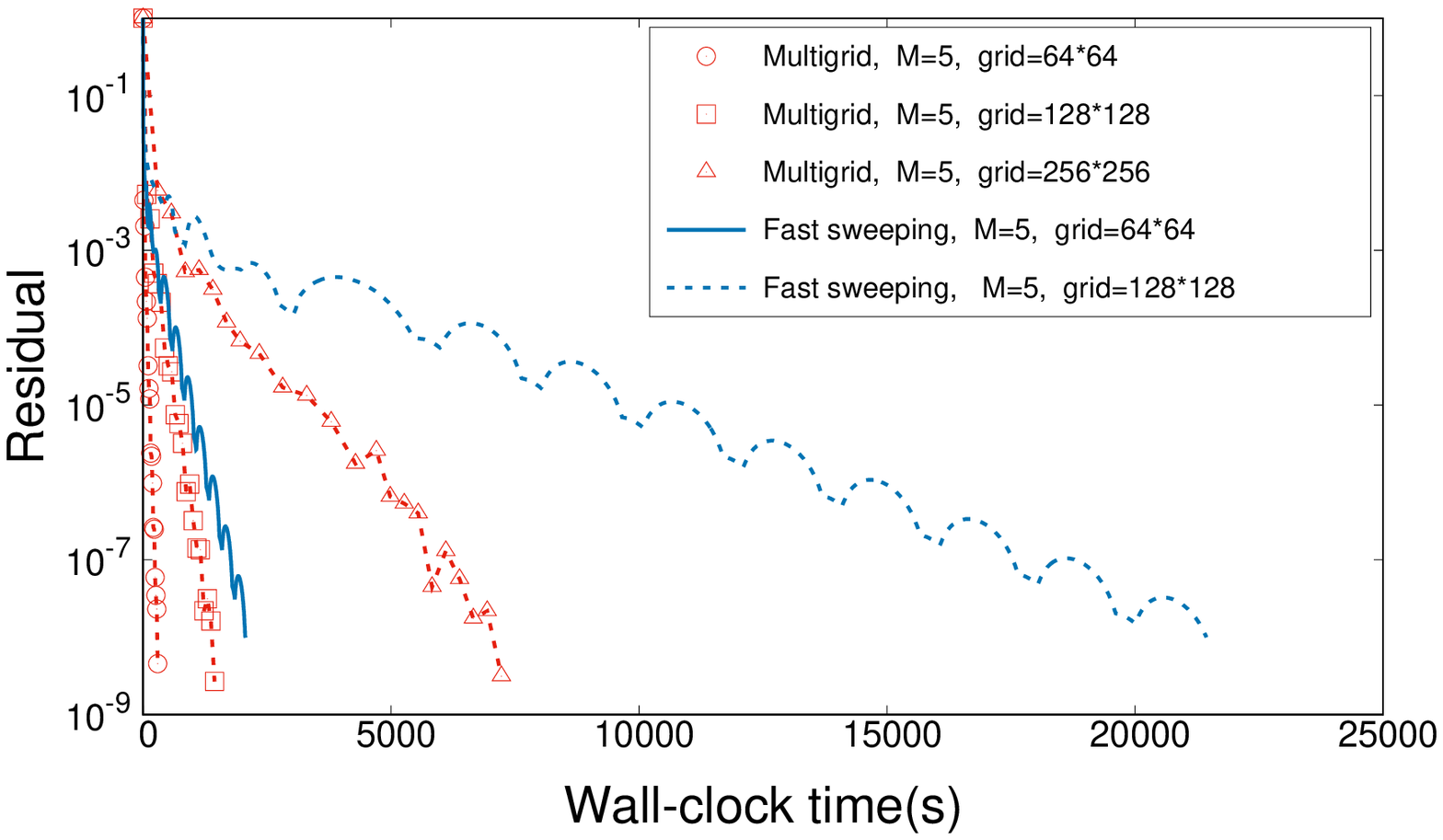}
   \hfill
   \includegraphics[width=0.49\textwidth]{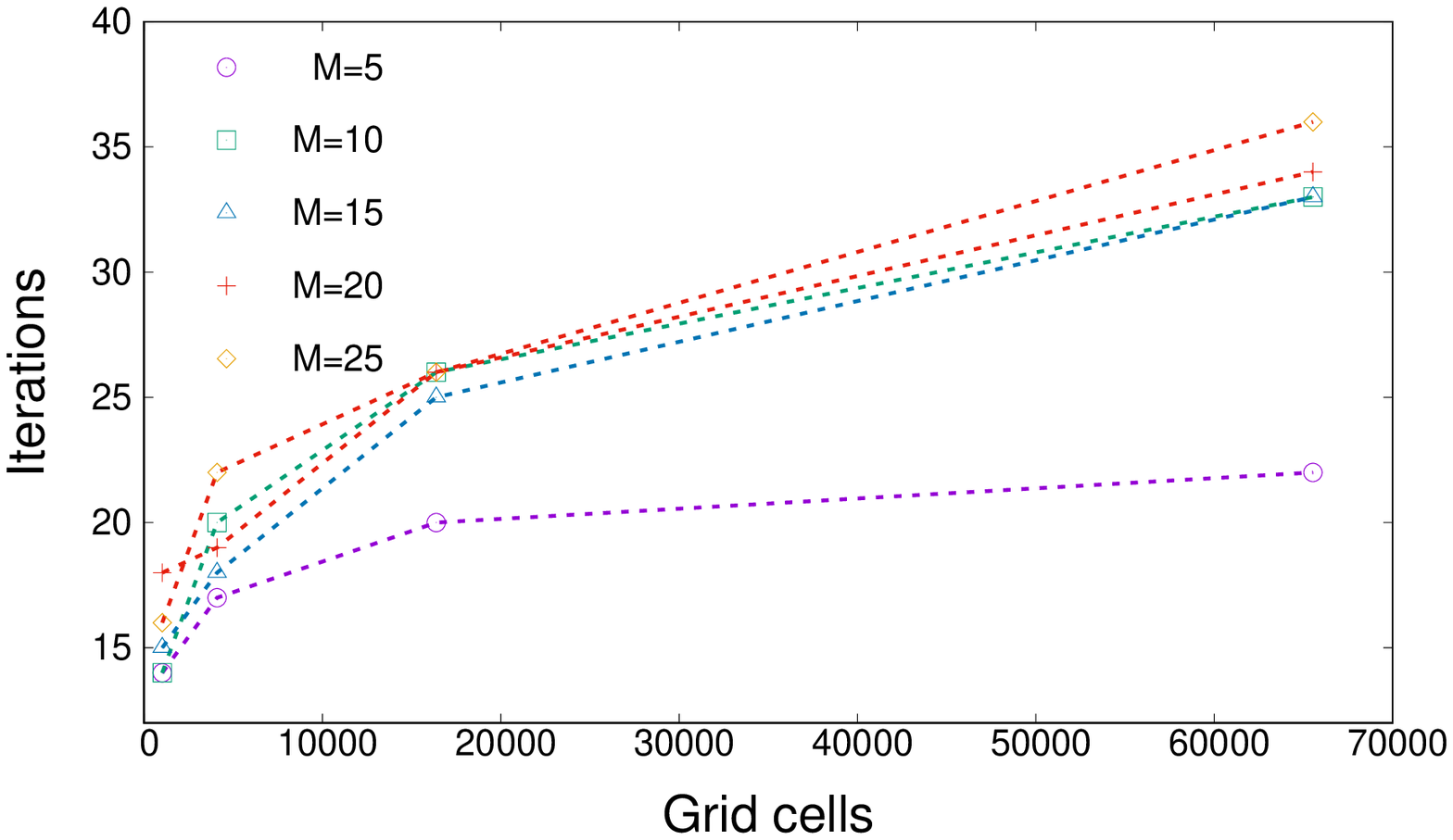}
    \caption{Convergence history (left) and
      iterations in terms of $N_{x}\times N_{y}$ of the NMG solver
      (right) for the Shakhov model of four-sided lid-driven flow.}
  \label{fig:kn_pt777_shakhov_cpu_ites}
\end{figure}

For the efficiency and behavior of the NMG solver, its results
compared with the two single level solvers on the grid composed of
$64\times 64$ cells for a variety of $M$ are listed in
\cref{tab:Kn_pt777_shakhov_M}. As can be seen once more, the fast
sweeping iteration converges faster than the forward Euler scheme for
all $M$, so that the wall-clock time is saved a lot in steady-state
computation. For the NMG solver, it converges within $20$ iterations
for all cases, resulting in the wall-clock time ratios of it to the
fast sweeping iteration being around $13\,\%$. When the grid is
refined, several convergence histories of the NMG solver and the fast
sweeping iteration for $M=5$, and the variations of total number of
iterations for the NMG solver with five values of $M$ are presented in
\cref{fig:kn_pt777_shakhov_cpu_ites}. These results are enough to show
the wonderful efficiency and behavior of the NMG solver. Additionally,
the efficiency in steady-state computation could be further improved
with a nice speedup ratio by using the parallel NMG solver, according
to the results shown in \cref{tab:four-bottom-parallelNMG}.

\begin{table}[!htbp]
  \centering
  \caption{Performance results for the Shakhov model of four-sided
    lid-driven flow with $N_{x}= N_{y} = 64$. Euler: the
    forward Euler scheme; FS: the fast sweeping iteration; NMG: the
    NMG solver; $T_{r}$: the wall-clock time ratio of NMG to FS.}
  \label{tab:Kn_pt777_shakhov_M}
  \begin{tabular}{p{1cm}<{\centering}ccccccc}
    \toprule
    \multirow{2}*{$M$} & \multicolumn{3}{c}{Iterations} & \multicolumn{4}{c}{Wall-clock time $(\si{\s})$}  \\
    \cmidrule(r){2-4} \cmidrule(r){5-8}
    & Euler & FS & NMG & Euler & FS & NMG & $T_{r}$\\
    \midrule
    5 & 4442 & 607 & 17 & 6031.19 & 2068.4 & 295.75 & 14.3\,\% \\
    10  & 4565 & 661 & 20 & 29037.5 & 12436.3 & 1693.24 & 13.6\,\% \\
    15  & 3783 & 621 & 18 & 71309.2 & 32209.3 & 4319.35 & 13.4\,\% \\
    20  & 5731 & 831 & 19 & 238658 & 103063 & 12709 & 12.3\,\% \\
  \bottomrule
  \end{tabular}
\end{table}

\begin{table}[!htbp]
  \centering
  \caption{Performance results of the parallel NMG solver with
    $N_{x}= N_{y} = 256$ in 
    \cref{sec:num-ex-four-lid,sec:num-ex-heat}. $T_{r,n}$: the
    wall-clock time ratio of single thread to $n$ threads.}
  \label{tab:four-bottom-parallelNMG}
  \begin{tabular}{ccccccccc}
    \toprule
    \multirow{2}*{\shortstack{Subsection}} & \multirow{2}*{$\Kn$} & \multirow{2}*{$M$} & \multirow{2}*{\shortstack{Average\\ iterations}}   & \multicolumn{5}{c}{Wall-clock time $(\si{\s})$ for $n$ threads} \\
    \cmidrule(r){5-9}
    &  &  &  & $n=1$ & $n=2$ & $T_{r,2}$ & $n=4$ & $T_{r,4}$\\
    \midrule
    \ref{sec:num-ex-four-lid} & $0.777$ & $25$ & 37 & 522770 & 266148 & 1.96 & 142230 & 3.68 \\
    \ref{sec:num-ex-heat} & $0.3$ & $20$ & 26 & 201206 & 102487 & 1.96 & 53074.4 & 3.79 \\
    \bottomrule
  \end{tabular}
\end{table}

\subsection{Heat transfer in a bottom-heated cavity}
\label{sec:num-ex-heat}

The third example is a heat transfer problem for the rarefied gas
confined in a bottom-heated cavity. As in \cite{2015microcavity,
  Cai2018}, all sides of the cavity are stationary and kept at
temperature $\SI{300}{\kelvin}$, except for the bottom side, which is
kept at temperature $\SI{600}{\kelvin}$. The length and height of the
cavity are taken as $L_{x} = L_{y} = \SI{e-6}{\m}$, and the initial
density, corresponding to the Knudsen number $0.3$, is set to
$\SI[per-mode=symbol]{0.2733}{\kg \per \m^{3}}$.
In addition, all the computations start from the initial global
Maxwellian with the given density, mean velocity of $0$, and
temperature of $\SI{300}{\kelvin}$.

\begin{figure}[!htbp]
  \centering
  \subfloat[Temperature $(\si{\kelvin})$]{
    \label{fig:kn_pt3_shakhov_T}
    \includegraphics[width=0.47\textwidth]{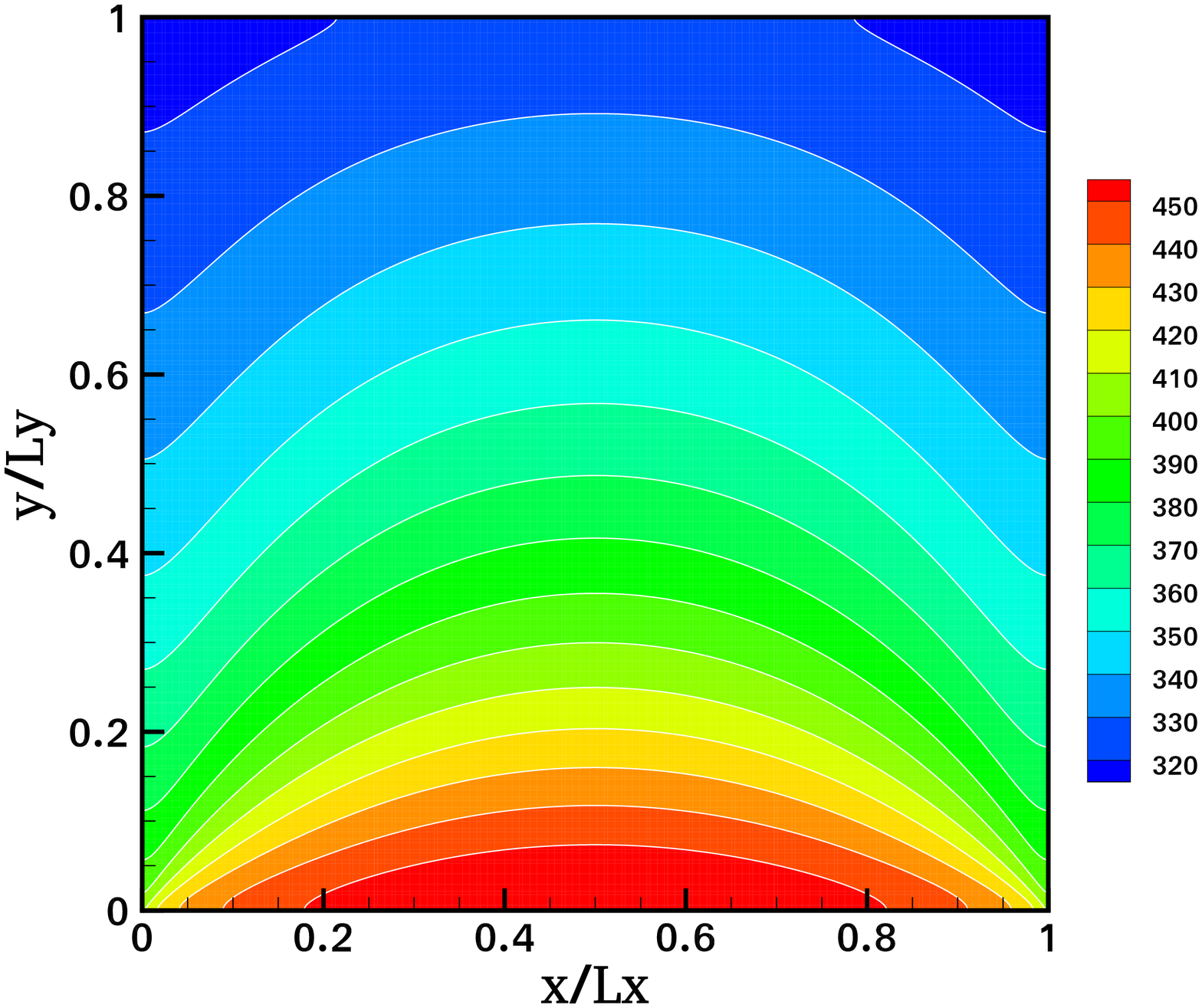}
  }\hfill
  \subfloat[Shear stress $(\si{\pascal})$]{
    \label{fig:kn_pt3_shakhov_sigma}
    \includegraphics[width=0.47\textwidth]{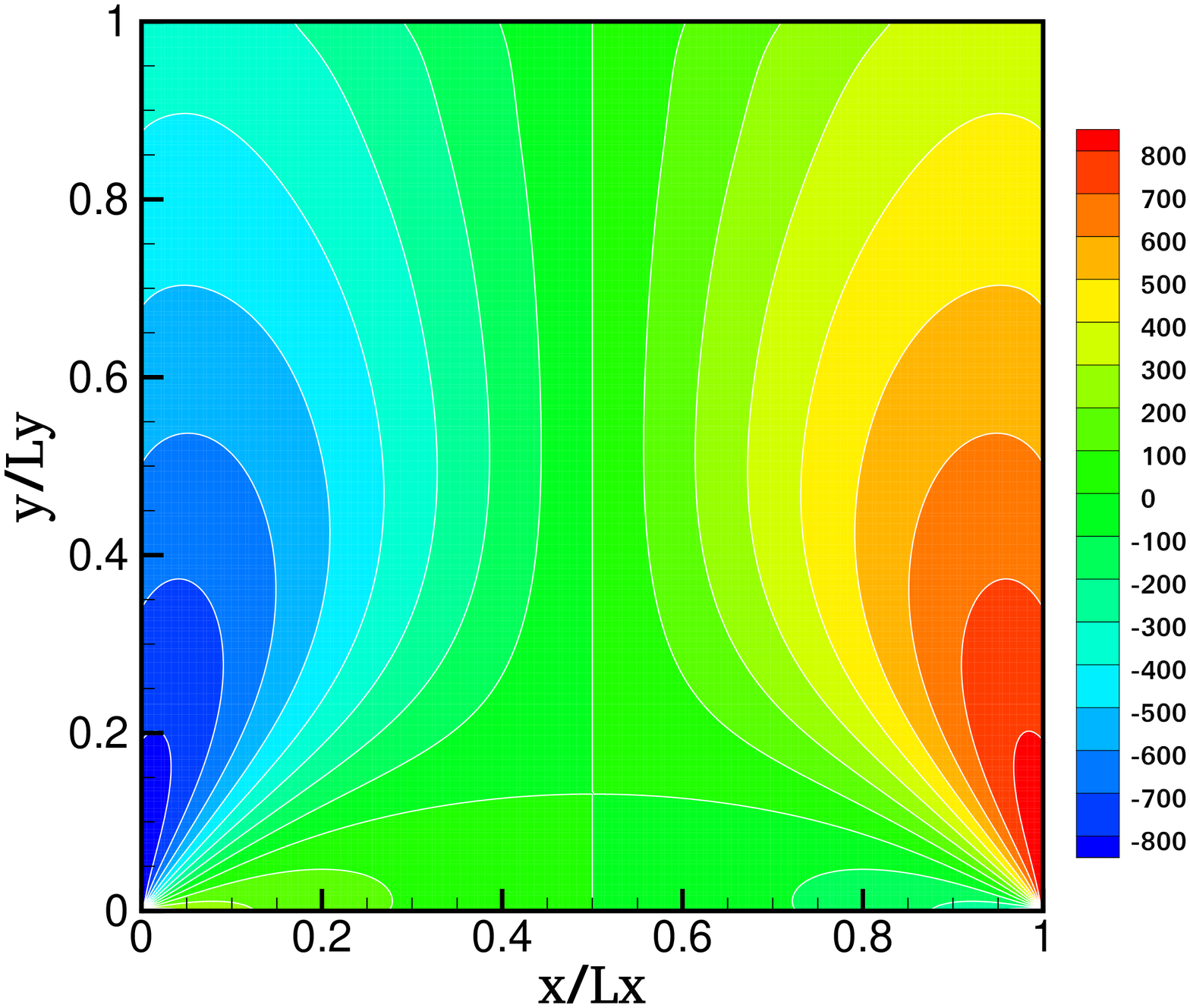}
  }
  \caption{Numerical solutions of the bottom-heated flow for the Shakhov model with $M=20$ and $N_{x}\times N_{y} = 256\times 256$.}
  \label{fig:kn_pt3_shakhov}
\end{figure}

Numerical solutions of the temperature and shear stress obtained by
the NMG solver for the Shakhov model with $M=20$ and
$N_{x}=N_{y} = 256$ are plotted in \cref{fig:kn_pt3_shakhov}. Compared
with the results presented in \cite{Cai2018}, it appears that the NMG
solver for the Shakhov model gives a reasonable distribution of the
temperature and shear stress.

\begin{figure}[!htbp]
  \centering
  \includegraphics[width=0.49\textwidth]{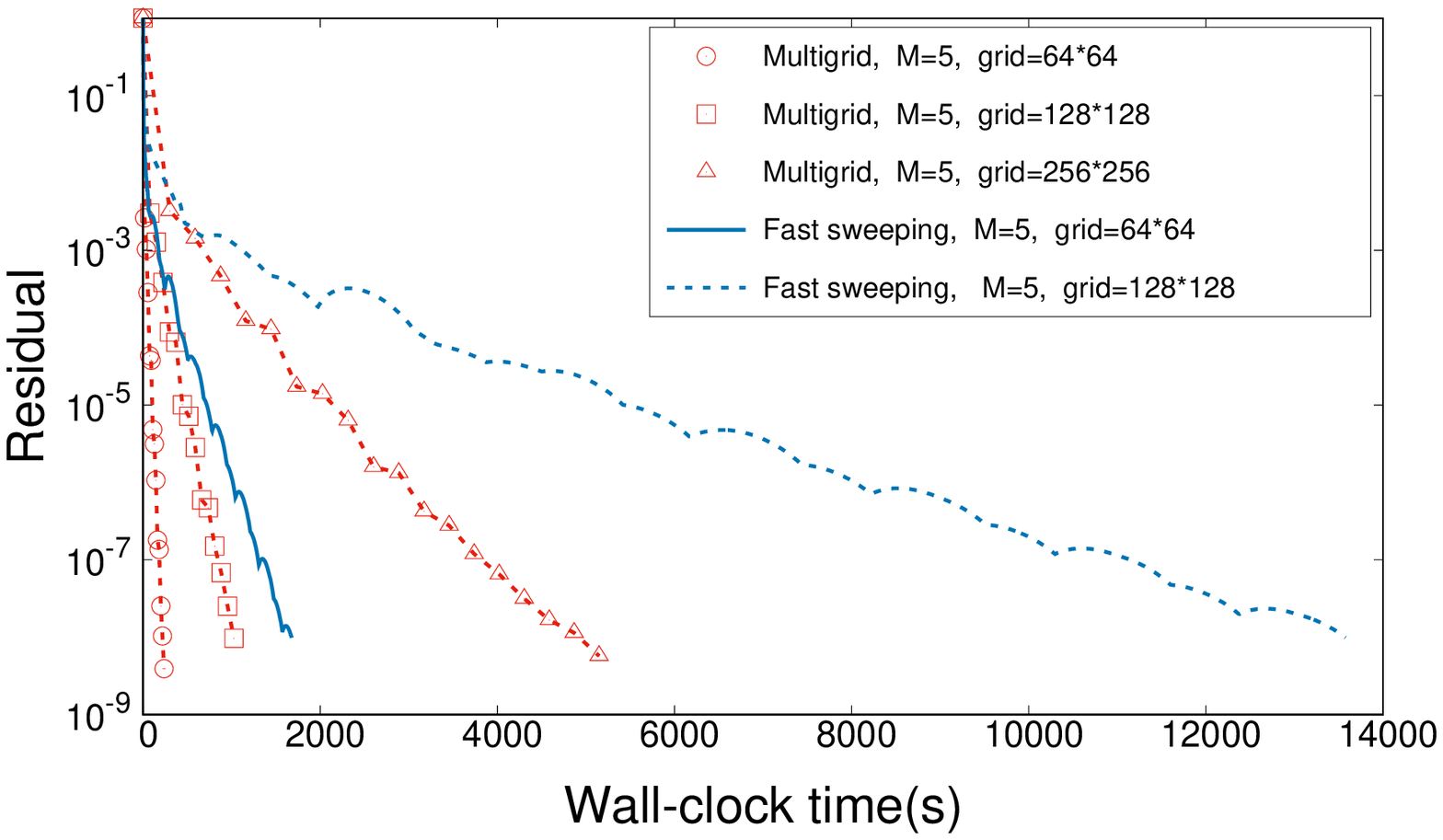}
   \hfill
   \includegraphics[width=0.49\textwidth]{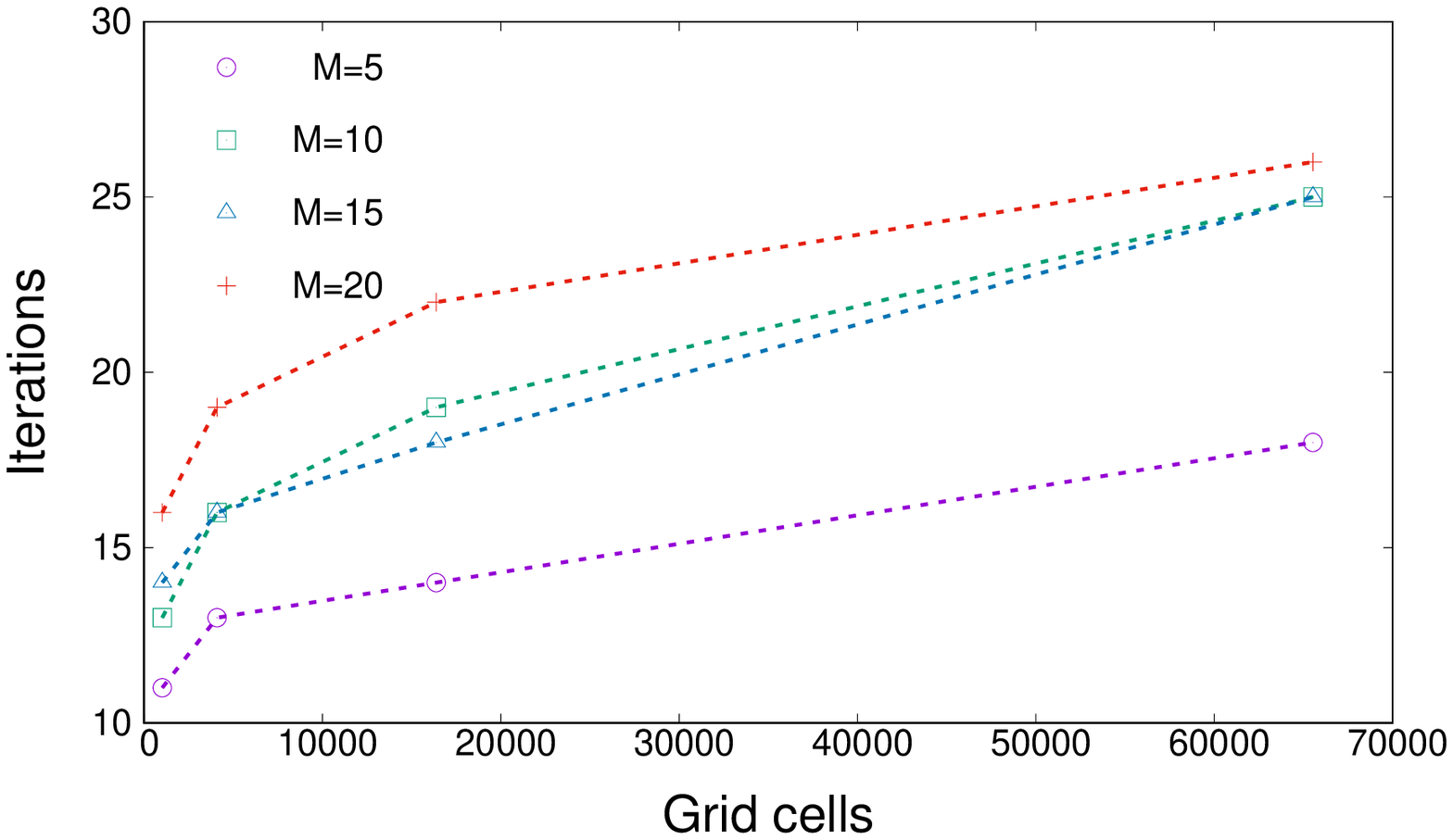}
    \caption{Convergence history (left) and 
      iterations in terms of $N_{x}\times N_{y}$ of the NMG solver
      (right) for the Shakhov model of bottom-heated flow.}
  \label{fig:kn_pt3_shakhov_cpu_ites}
\end{figure}

Then, the performance of the NMG solver for the above problem,
compared with the two single level solvers, is investigated. The
results on the grid with $64\times 64$ cells for several choices of
$M$ are presented in \cref{tab:Kn_pt3_shakhov_M}. With these values of
$M$, the variations of total number of iterations for the NMG solver
with respect to $N_{x}\times N_{y}$, as well as their partial
convergence histories,
are shown in \cref{fig:kn_pt3_shakhov_cpu_ites}. For the speedup of
simulation by using the parallel NMG solver, it can be found in
\cref{tab:four-bottom-parallelNMG}.
From all of these results, we deduce that the NMG solver performs
analogously as in the previous examples. That is, the NMG solver
always converges in dozens of iterations. Consequently, compared with
the fast sweeping iteration, which is already more efficient than the
forward Euler scheme, the NMG solver is able to accelerate the
steady-state computation more significantly, especially in the case
that the number of grid cells is large. When multi-threads are
applied, the parallel NMG solver would further accelerate the
simulation with a pretty speedup ratio.

\begin{table}[!htbp]
  \centering
  \caption{Performance results for the Shakhov model of bottom-heated
    flow with $N_{x}= N_{y} = 64$. Euler: the forward Euler scheme;
    FS: the fast sweeping iteration; NMG: the NMG solver; $T_{r}$: the
    wall-clock time ratio of NMG to FS.}
  \label{tab:Kn_pt3_shakhov_M}
  \begin{tabular}{p{1cm}<{\centering}ccccccc}
    \toprule
    \multirow{2}*{$M$} & \multicolumn{3}{c}{Iterations} & \multicolumn{4}{c}{Wall-clock time $(\si{\s})$}  \\
    \cmidrule(r){2-4} \cmidrule(r){5-8}
    & Euler & FS & NMG & Euler & FS & NMG & $T_{r}$\\
    \midrule
    5 & 4073 & 499 & 13 & 5668.54 & 1677.59 & 236.72 & 15.7\,\% \\
    10  & 4926 & 589 & 16 & 31041.7 & 9739.9 & 1404.25 & 14.4\,\% \\
    15  & 5357 & 664 & 16 & 100517 & 33575.4 & 4039.22 & 12.0\,\% \\
    20  & 5991 & 745 & 19 & 253515 & 93542.3 & 10999.1 & 11.8\,\% \\
  \bottomrule
  \end{tabular}
\end{table}

\subsection{Second-order spatial discretization for single lid-driven cavity flow}
\label{sec:num-ex-second-order}

\begin{figure}[htbp]
  \centering
  \subfloat[Temperature $(\si{\kelvin})$]{
    \label{fig:kn_pt1_shakhov_T_order2}
    \includegraphics[width=0.47\textwidth]{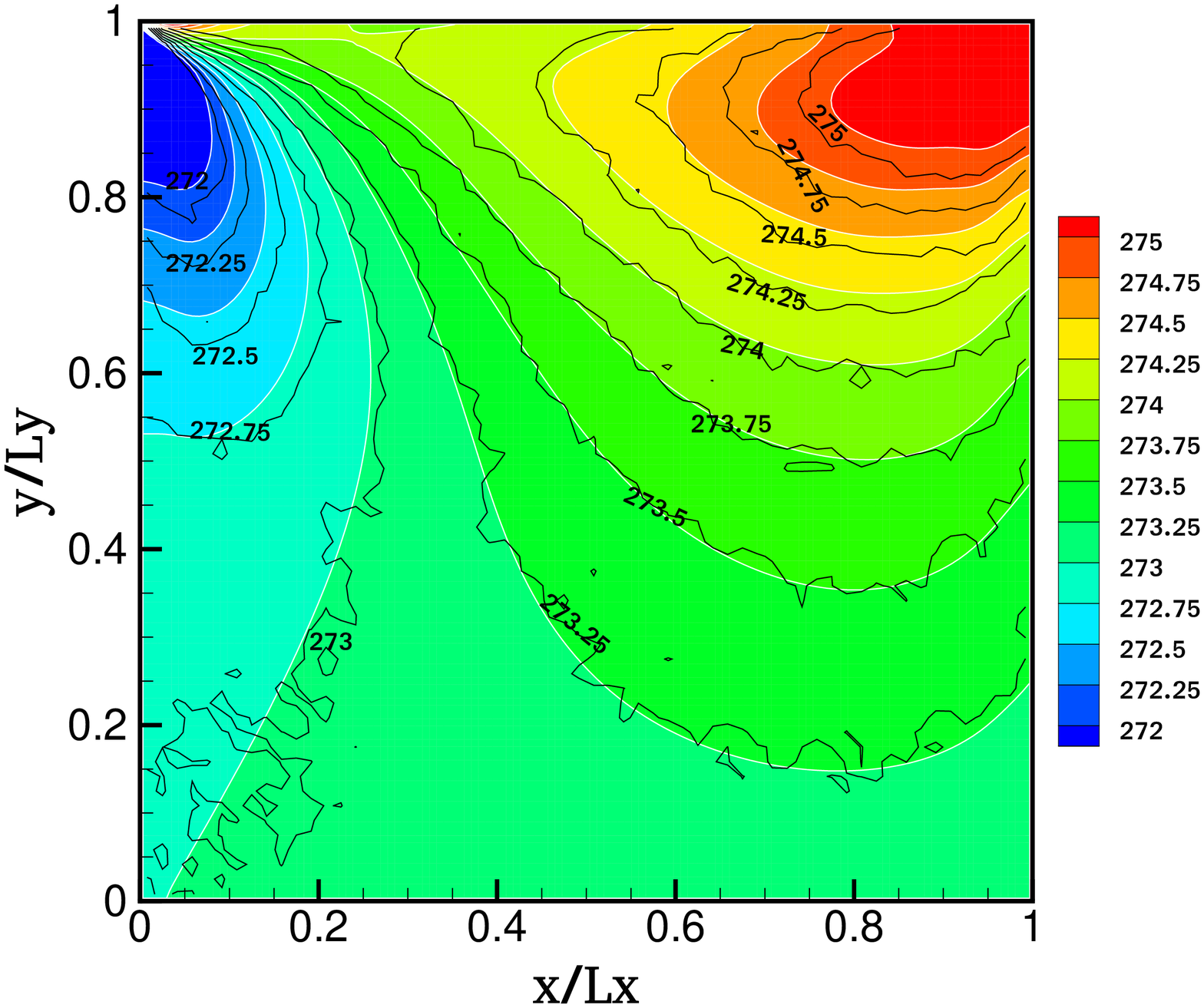}
  }\hfill
  \subfloat[Heat flux streamlines]{
    \label{fig:kn_pt1_shakhov_flux_order2}
    \includegraphics[width=0.47\textwidth]{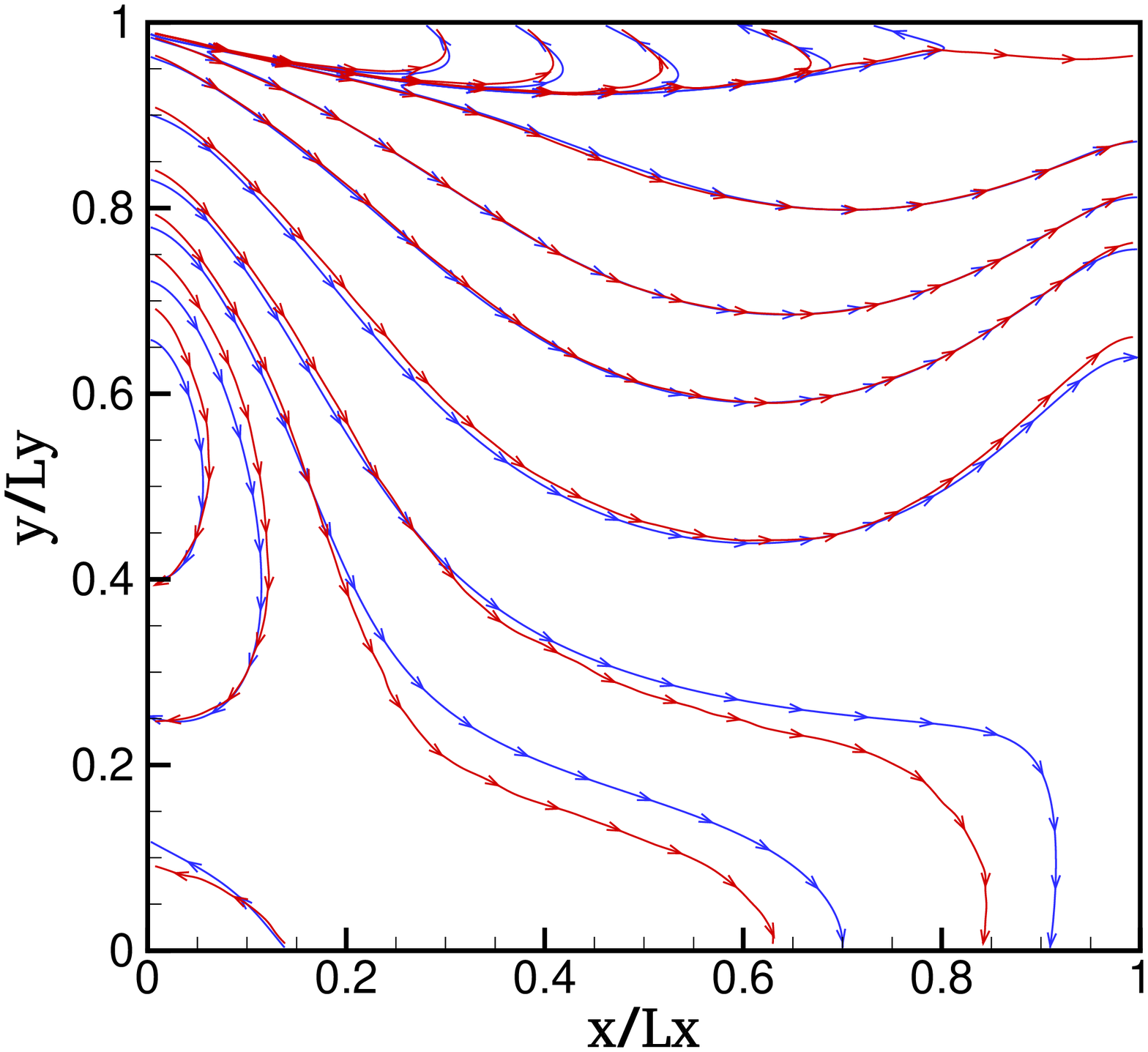}
  }
  \caption{Numerical solutions of the single lid-driven flow for the Shakhov model with  $\Kn=0.1$, $M=20$, $N_{x}\times N_{y} = 128\times 128$, and second-order spatial discretization. The black (left) and red (right) lines are the reference DSMC solutions.}
  \label{fig:kn_pt1_shakhov_order2}
\end{figure}

In the last subsection, the NMG solver for the second-order spatial
discretization is investigated by recomputing the single lid-driven
cavity flow.
The Shakhov model and the Knudsen number $\Kn=0.1$ are
applied. Numerical solutions of the temperature and heat flux for
$M=20$ and $N_{x}=N_{y}=128$, together with the reference DSMC
solutions, are presented in \cref{fig:kn_pt1_shakhov_order2}. As
anticipated, we can observe that the temperature contours and heat
flux streamlines using the second-order spatial discretization on a
coarser grid with $128\times 128$ cells coincide with the DSMC results
much better than the results using the first-order spatial
discretization on the grid with $256\times 256$ cells that are shown
in \cref{fig:kn_pt1_shakhov}, especially in the region near the
boundary.

\begin{figure}[!htbp]
  \centering
    \includegraphics[width=0.49\textwidth]{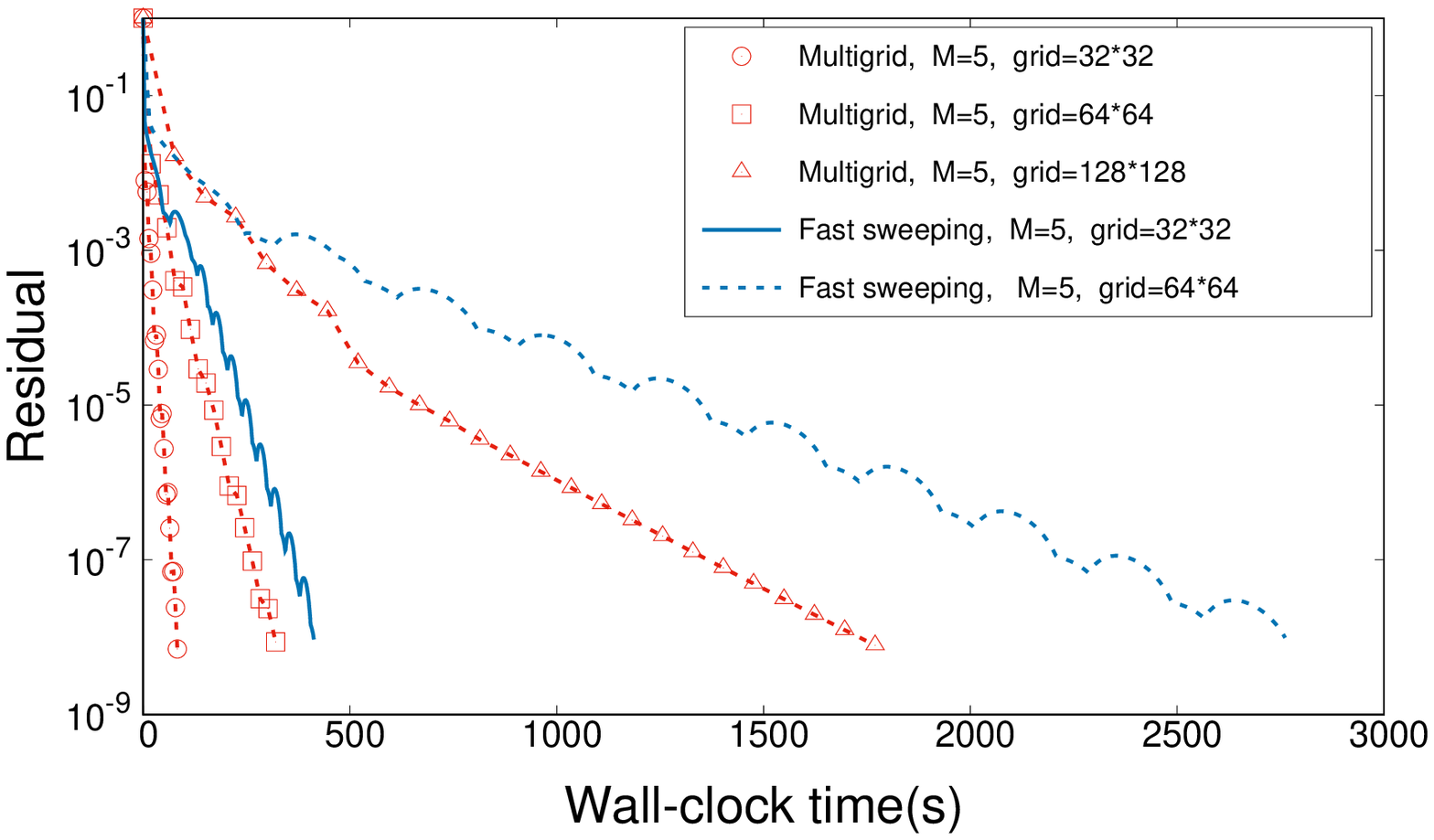}
   \hfill
    \includegraphics[width=0.49\textwidth]{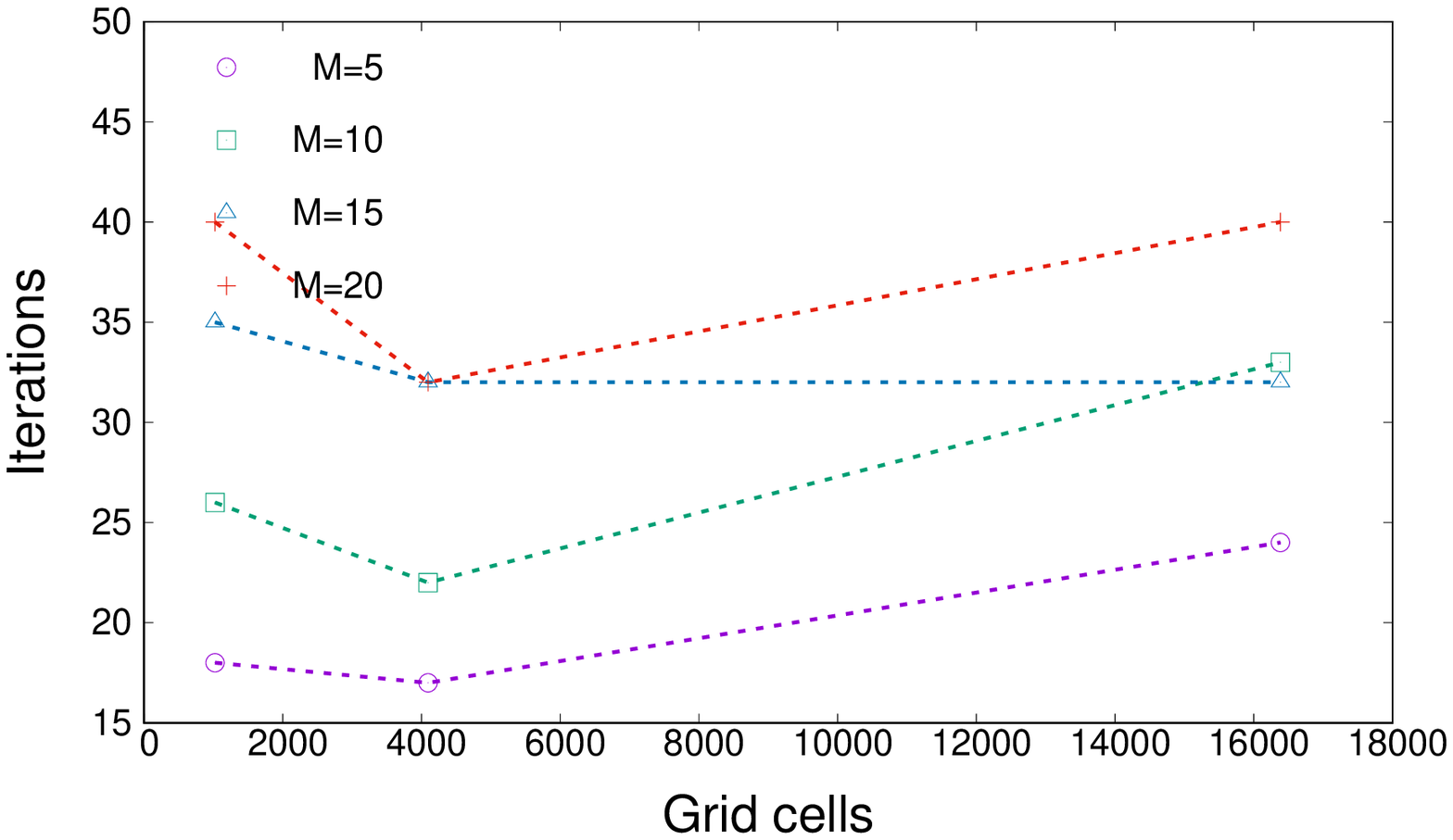}
    \caption{Convergence history (left) and 
      iterations in terms of $N_{x}\times N_{y}$ of the NMG solver
      (right) for the Shakhov model of single lid-driven flow with
      $\Kn=0.1$, and second-order spatial discretization.}
  \label{fig:kn_pt1_shakhov_cpu_ites_order2}
\end{figure}

For the second-order spatial discretization, it has been observed in
our simulation that the forward Euler scheme fails to converge due to
instability, while the fast sweeping iteration still works well so
that the NMG solver using it as the smoother is able to give
satisfactory results too.
In \cref{tab:kn_pt1_grid_second}, the total number of iterations and
the wall-clock time, spent by the NMG solver and the fast sweeping
iteration for the second-order spatial discretization in the case of
$M=5$ on three grids, are listed. Compared with the results shown in
\cref{tab:Kn_pt1_Shakhov_grid} for the first-order spatial
discretization, it turns out that, on the same grid, the fast sweeping
iteration converges more slowly for second-order case, whereas as the
grid is refined, the resulting total number of iterations increases
with a better growth rate for second-order case.
Meanwhile, the NMG solver also behaves somewhat differently from the
first-order case. It can be seen from \cref{tab:kn_pt1_grid_second} as
well as \cref{fig:kn_pt1_shakhov_cpu_ites_order2} (right) that the
total number of iterations of the NMG solver for almost all $M$
decreases first and then increases on the three grids ranged from
$32\times 32$ to $128\times 128$ cells. For more details, it can be
found from \cref{fig:kn_pt1_shakhov_cpu_ites_order2} (left) that, due
to the fluctuation of residuals generated by the smoother, i.e., the
fast sweeping iteration, the convergence history of the NMG solver on
the smallest grid fluctuates a little more than that on the other
grids. Apart from this, the convergence rate of the NMG solver on the
finest grid appears to be degenerate slightly after a few number of
iterations. As a result, a little more iterations to obtain the steady
state are spent by the NMG solver on these grids than on the middle
grid.

Nevertheless, we still have that the NMG solver converges in dozens of
iterations for all simulations, so that it could improve the
efficiency remarkably in steady-state computation.
Additionally, it is clear from
\cref{tab:single-lid-parallelNMG-order2} that the expected speedup
could be further obtained by using the NMG solver with multi-threads.

\begin{table}[!htbp]
  \centering
  \caption{Performance results for the Shakhov model of single
    lid-driven flow with $\Kn=0.1$, $M=5$ and second-order spatial
    discretization. FS: the fast sweeping iteration; NMG: the NMG
    solver; $T_{r}$: the wall-clock time ratio of NMG to FS.}
  \label{tab:kn_pt1_grid_second}
  \begin{tabular}{cccccc}
    \toprule
    \multirow{2}*{$N_{x}\times N_{y}$} & \multicolumn{2}{c}{Iterations} & \multicolumn{3}{c}{Wall-clock time $(\si{\s})$} \\
    \cmidrule(r){2-3} \cmidrule(r){4-6}
    & FS & NMG & FS & NMG & $T_{r}$ \\
    \midrule
    $32 \times 32$  & 384 & 18 & 412.78 & 82.80 & 20.1\,\% \\
    $64 \times 64$   & 728 & 17 & 2761.07 & 320.25 & 11.6\,\% \\
    $128 \times 128$   & 1331 & 24 & 19393.3 & 1769.73 & 9.1\,\% \\
    \bottomrule
  \end{tabular}
\end{table}

\begin{table}[!htbp]
  \centering
  \caption{Performance results of the parallel NMG solver for the
    Shakhov model of single lid-driven flow with $\Kn=0.1$, $M=20$ and
    second-order spatial discretization. $T_{r,n}$: the wall-clock
    time ratio of single thread to $n$ threads.}
  \label{tab:single-lid-parallelNMG-order2}
  \begin{tabular}{ccccc}
    \toprule
    \multirow{2}*{$N_{x}\times N_{y}$} & \multirow{2}*{\shortstack{Average\\ iterations}}   & \multicolumn{3}{c}{Wall-clock time $(\si{\s})$ for $n$ threads} \\
    \cmidrule(r){3-5}
    &  & $n=1$ & $n=2$ & $T_{r,2}$ \\
    \midrule
    $64\times 64$ & 36 & 19307.8 & 11485.8 & 1.68 \\
    $128\times 128$ & 40 & 81154.8 & 41326.2 & 1.96 \\
    \bottomrule
  \end{tabular}
\end{table}

\section{Conclusion}
\label{sec:conclusion}

Aiming at efficient steady-state simulation of rarefied gas cavity
flow described by the Boltzmann equation with BGK-type collision term,
a nonlinear multigrid solver has been successfully developed based on
the following approaches. At first, it adopts the unified framework of
regularized moment method for velocity discretization and finite
volume method for spatial discretization. To solve the resulting
discrete problem, a fast sweeping iteration that converges faster and
more robust than the time-integration scheme is introduced. Then, the
NMG method, which employs the fast sweeping iteration as the smoother,
is proposed to greatly improve the convergence rate. Finally, the
OpenMP-based parallelization is applied for the implementation of the
NMG method, so that the efficiency could be further improved by
multi-threaded parallel computation. Plenty of numerical experiments
have been carried out to investigate the performance of the resulting
NMG solver. All numerical results show the efficiency and robustness
of the solver for both first- and second-order spatial discretization.

Incidentally, it is easy to extend the proposed NMG solver to the
Boltzmann equation with other collision models by providing the
computation of the coefficients $Q_{ij,\alpha}$ in
\cref{eq:flux-force-expansion}. Indeed, combining the algorithm
presented in \cite{Wang2019}, the performance of the NMG solver
for the Boltzmann equation with the quadratic collision term is under
investigation and will be reported elsewhere.

\bibliographystyle{plain}
\bibliography{MG2022.bib}
\end{document}